%% file: main.tex
\documentclass[10pt, oneside]{amsart}
\input{preamble}

\title{Beilinson-Parshin adeles via solid algebraic geometry}
\author{Christopher Brav and Grigorii Konovalov}
\date{July 2025}

\begin{document}

\begin{abstract}
We apply Clausen-Scholze's theory of solid modules to the existence of adelic decompositions for schemes of finite type 
over $\mathbb{Z}$. Specifically, we use the six-functor formalism for solid modules to define the skeletal filtration of a scheme,
and then we show  that decomposing a quasi-coherent sheaf with respect to this filtration gives rise to a new construction of
the Beilinson-Parshin adelic resolution. As an application of the adelic decomposition combined with some nice completeness
properties of the solid tensor product, we prove a version of adelic descent for solid quasi-coherent sheaves.
\end{abstract}
\maketitle

\input{intro}
\input{general_solid}
\input{some_algebra}
\input{skeleton}

\input{main_thm}
\input{global}

\input{canonical_colim_of_affs}
\input{appendixB}

\bibliographystyle{IEEEtran}
\bibliography{refs}{}

\bigskip
\bigskip

\noindent Christopher~Brav, {\sc Shanghai Institute for Mathematics and Interdisciplinary Sciences, Shanghai,}
\href{mailto:brav.k@mipt.ru}{brav@simis.cn}

\smallskip

\noindent Grigorii~Konovalov, {\sc Center for Advanced Studies, Skoltech, Moscow; Centre of Pure Mathematics, MIPT, Moscow;
Faculty of Mathematics, National Research University Higher School of Economics, Moscow;}
\href{mailto:grisha.v.konovalov@gmail.com}{grisha.v.konovalov@gmail.com}

\end{document}

%% file: preamble.tex
\usepackage[hscale = 0.8, top=0.8in]{geometry}
\usepackage{appendix}

\usepackage{xcolor}
\usepackage{float}
\usepackage{amsmath}      
\usepackage{amsthm}        
\usepackage{amssymb}       
\usepackage{empheq}       
\usepackage{stackrel}
\usepackage{braket}
\usepackage[all]{xy}
\usepackage{hyperref}
\hypersetup{
  colorlinks   = true,          
  urlcolor     = pink,          
  linkcolor    = red, 
  citecolor   = blue            
}

\makeatletter
\newcommand{\makesmaller}[2]{#1{\mathpalette\make@smaller@{#2}}}
\newcommand{\make@smaller@}[2]{%
  \vcenter{\hbox{$\m@th\downgrade@style#1#2$}}%
}
\newcommand{\downgrade@style}[1]{%
  \ifx#1\displaystyle\scriptstyle\else
    \ifx#1\textstyle\scriptstyle\else
      \scriptscriptstyle
  \fi\fi
}
\makeatother

\renewcommand{\theequation}{\thesubsubsection.\arabic{equation}}         
\makeatletter                                                      
   \@addtoreset{equation}{subsubsection} 
   \@addtoreset{footnote}{subsubsection}

\def\newprooflikeenvironment#1#2#3#4{%                             
     \newenvironment{#1}[1][]{%                                    
         \refstepcounter{equation}%                                 
         \begin{proof}[{\rm\csname#4\endcsname{#2~\theequation}}]%
         \def\qedsymbol{#3}}%                                      
        {\end{proof}}}                                             
\makeatother                                                   
                                                                   
\newprooflikeenvironment{definition}{Definition}{}{textbf}
\newprooflikeenvironment{example}{Example}{}{textbf}     
\newprooflikeenvironment{remark}{Remark}{}{textbf}       
\newprooflikeenvironment{notation}{Notaion}{}{textbf}
\newprooflikeenvironment{observation}{Observation}{}{}
\newprooflikeenvironment{quasiobservation}{}{}{}
                                                                   
\theoremstyle{plain}                                               
\newtheorem{theorem}[subsubsection]{Theorem}                            
                          
\newtheorem{lemma}[subsubsection]{Lemma}                                
\newtheorem{corollary}[subsubsection]{Corollary}                        
\newtheorem{proposition}[subsubsection]{Proposition}

\renewcommand{\a}{\alpha}
\newcommand{\arr}{\longrightarrow}

\newcommand{\calg}{\mathrm{CAlg}}
\newcommand{\A}{\mathcal{A}}

\newcommand{\bs}{{\makesmaller{\mathrel}{\blacksquare}}}

\newcommand{\CC}{\mathcal{C}}

\newcommand{\comp}{{\makesmaller{\mathrel}{\wedge}}}

\newcommand{\colim}{\mathrm{colim}}

\newcommand{\DD}{\mathcal{D}}

\newcommand{\define}{\text{def}}
\renewcommand{\dim}{\mathrm{dim}}

\newcommand{\f}{\phi}

\newcommand{\fib}{\mathrm{fib}}

\newcommand{\g}{\psi}
\newcommand{\Ga}{\Gamma}
\renewcommand{\hom}{\text{Hom}}

\newcommand{\I}{\mathcal{I}}
\DeclareMathOperator{\id}{id}

\newcommand{\im}{\mathrm{Im}}
\newcommand{\ind}{\mathrm{Ind}}

\newcommand{\leftcomp}{{\stackrel{\comp}{\text{\phantom{a}}}\!\!}}

\newcommand{\map}{\mathrm{Map}}

\renewcommand{\mod}{\mathrm{Mod}}
\newcommand{\N}{\mathbb{N}}
\renewcommand{\O}{\mathcal{O}}

\newcommand{\op}{\mathrm{op}}
\newcommand{\p}{\mathfrak{p}}
\renewcommand{\P}{\mathcal{P}}
\newcommand{\pb}{\rule{.4pt}{5.4pt}\rule[5pt]{5pt}{.4pt}\llap{\hspace{1pt}}}

\newcommand{\projfinlen}{\mathrm{proj,\: fin\text{-}len}}

\newcommand{\perf}{\text{\sf Perf}}
\newcommand{\qcoh}{\text{QCoh}}
\newcommand{\Q}{\mathbb{Q}}
\newcommand{\q}{\mathfrak{q}}

\DeclareMathOperator{\spec}{Spec}
\newcommand{\anspec}{{\mathrm{Spec}^\mathrm{an}}}

\renewcommand{\t}{\theta}

\DeclareMathOperator{\uhom}{\underline{Hom}}
\newcommand{\w}{\omega}
\newcommand{\Z}{\mathbb{Z}}
\newcommand{\ZZ}{\mathcal{Z}}

\newcommand{\inveq}{\mathrel{\rotatebox[origin=c]{90}{$=$}}}
\newcommand{\invsimeq}{\mathrel{\rotatebox[origin=c]{90}{$\simeq$}}}

\newcommand{\invertin}{\mathrel{\rotatebox[origin=c]{90}{$\in$}}}
\newcommand{\invertinn}{\mathrel{\rotatebox[origin=c]{270}{$\in$}}}

\renewcommand{\Psi}{L}
\newcommand{\fakej}{\phi}

\newcommand{\ssec}{\subsection}
\newcommand{\sssec}{\subsubsection}
\def\on{\operatorname}

\def\blox{\bs}
\def\QCoh#1{\on{QCoh}(#1)}
\def\QCohsup#1#2{\on{QCoh}_{#1}(#2)}

\def\Spec#1{\on{Spec}(#1)}
\def\Specan#1{\on{Spec}^{an}(#1)}

\renewcommand{\thicksim}{{\overset{\sim}{\!\!\phantom{\ldots}}}}

\newcommand{\idem}{\mathrm{Idem}}
\newcommand{\cidem}{\mathrm{cIdem}}

\newcommand{\sm}{\mathrm{Sm}}
\newcommand{\gen}{\mathrm{gen}}
\newcommand{\spc}{\mathrm{spc}}
\newcommand{\cmp}{\mathrm{cmp}}
\newcommand{\supp}{\mathrm{supp}}

\newcommand{\cls}{\mathrm{Cls}}
\newcommand{\sd}{\mathrm{sd}}
\newcommand{\pr}{\mathrm{Pr}}
\newcommand{\st}{\mathrm{st}}
\newcommand{\gr}{\mathrm{gr}}

%% file: intro.tex
\section*{Introduction}

Recall the `Hasse fracture square' or `adelic decomposition', which reconstructs the integers from its $p$-completions and its rationalization. Namely, in the derived category of abelian groups there is a homotopy pullback square

\[
\xymatrix{\Z \ar[r] \ar[d] & \Q \ar[d] \\
\prod_{p} \Z_{p} \ar[r] & \Q \otimes \prod_{p} \Z_{p} }
\]
More generally, there is an analogous square for any object of the derived category of abelian groups. Parshin \cite{ParshinAdeles}, then Beilinson \cite{BeilinsonResidues}, constructed a higher dimensional generalization of the adelic decomposition that gives a functorial resolution of a quasi-coherent sheaf on a nice scheme in terms of iterated localisations and completions along all flags of closed integral subschemes.

The goal of this paper is to show how Clausen-Scholze's condensed mathematics \cite{condensed}, more specifically the theory of solid modules, leads to a new construction of Beilinson-Parshin adelic decompositions for schemes of finite type over $\Z$.\footnote{For discussion of more general cases, see comment \ref{rmk-evrth-works-over-a-field} at the end of the introduction.} The basic idea is to use the enhanced locality and functoriality of the category of solid modules to define the skeletal filtration of a scheme, then to decompose a solid quasi-coherent sheaf (in particular, a usual quasi-coherent sheaf) with respect to this filtration. The result is a reconstruction of any (solid) quasi-coherent sheaf via a homotopy limit cube, generalizing the Hasse square above. For a perfect complex, the vertices are given as products over flags of integral closed subschemes of iterated completed localizations, just as in the Beilinson-Parshin resolution, and the general case is reduced to this via commutation with colimits.  

In contrast to the constructions of Parshin \cite{ParshinAdeles}, Beilinson \cite{BeilinsonResidues}, and Balchin-Greenlees \cite{greenlees2020adelic}, which put iterated completed localizations into the construction of the resolution and which require work to check exactness of the resulting complex, in our construction the exactness comes by construction and the iterated completed localizations arise a posteriori, as a calculation of the output of the decomposition of a quasi-coherent sheaf with respect to the skeletal filtration.

In order to motivate the construction, consider a nice topological space $X$ together with a closed subspace $i: Z \hookrightarrow X$ and its open complement $j : U \hookrightarrow X$. Then the six-functor formalism for the (derived category) of abelian sheaves $Sh(X)$ gives a recollement

\[
\xymatrix{
Sh(Z)\ar@{^{(}->}[rr]|{{i_*}} && Sh(X) \ar@<2ex>[ll]^{i^{\:!}}\ar@<-2ex>[ll]_{i^*} \ar[rr]|{j^*} && Sh(U) \ar@<2ex>@{_{(}->}[ll]^{j_*}\ar@<-2ex>@{_{(}->}[ll]_{{j}_!} \;.
}
\]

From this we obtain for any complex of abelian sheaves $F$ two recollement triangles from restriction to $Z$ and to $U$:
\[
j_{!}j^{*}F \arr F \arr i_{*}i^{*}F
\]
\[
i_{*}i^{!}F \arr F \arr j_{*}j^{*}F
\] 
Comparing the two types of restrictions, we obtain a homotopy pullback square
\[
\xymatrix{F \ar[r] \ar[d] & j_{*}j^{*}F \ar[d] \\
i_{*}i^{*}F \ar[r] &  i_{*}i^{*}j_{*}j^{*}F}
\]
Note however that in general, the functor $j_{*}$ does not commute with all direct sums, so the square can be difficult to compute.

Consider now a scheme $X$ of finite type over $\Z$ and its (derived) category of quasi-coherent sheaves $\QCoh{X}$. Given a closed subspace $i: Z \hookrightarrow X$, and its open complement $j : U \hookrightarrow X$, we have continuous adjunctions $j^{*} : \QCoh{X} \leftrightarrow \QCoh{U}: j_{*}$ and $i_{*} : \QCohsup{Z}{X} \leftrightarrow \QCoh{X}:i^{!}$, where $\QCohsup{Z}{X}$ is the full subcategory of quasi-coherent sheaves on $X$ `supported' in $Z$. However, the left adjoint $j_{!}$ of $j^{*}$ usually does not exist, since $j^{*}$ does not preserve arbitrary products. Similarly, the left adjoint $i^{*}$ of $i_{*}$ usually does not exist. We therefore have only one half of a recollement, namely the triangle
\[
i_{*}i^{!}F \arr F \arr j_{*}j^{*}F.
\]
In short, even for some of the simplest maps of schemes, there is not a full six-functor formalism for quasi-coherent sheaves. 

However, using condensed mathematics, Clausen and Scholze (\cite{condensed},\cite{analytic},\cite{clausen2023analytic}) have shown how to enlarge the context of algebraic geometry to an `analytic geometry' built on analytic spectra $\Specan{A}$ of `analytic rings' $A$ and more general spaces glued from such spectra. Quasi-coherent sheaves on analytic spaces then admit a full six-functor formalism for a much wider class of maps. In particular, there is a fully faithful embedding of usual algebraic geometry into analytic geometry that takes the usual spectrum $\Spec{A}$ of a finite type discrete ring $A$ to the analytic spectrum $\Specan{A_{\bs}}$ of the `solid analytic ring structure' $A_{\bs}$ on $A$. Concretely, considering the solid analytic ring structure means specifying a certain category $\mod_{A_\bs}$ of `complete' $A$-modules generated under colimits by products $\prod_{I}A$ over index sets $I$, we refer to section \ref{sec-gen-solid} and references given there for the detailed discussion of this construction. By definition, the quasi-coherent sheaves $\QCoh{\Specan{A_{\bs}}}$ on $\Specan{A_{\bs}}$ is the category $\mod_{A_\bs}$ of `complete' or `solid' $A$-modules, which admits a natural continuous and symmetric monoidal embedding $\mod_A \hookrightarrow \mod_{A_\bs}$. We emphasize however that the embedding does not in general preserve infinite limits. Indeed, the whole point of the category $\mod_{A_\bs}$ is to ensure that infinite limits are somehow appropriately complete.

Now consider $X=\Spec{A}$ where $A$ is of finite type over $\Z$, $Z=\Spec{A/I} \subset X$ a closed subscheme, and $U=X \setminus Z$ the open complement, and consider the embedding of these objects $X_{\bs},U_{\bs},Z_{\bs}$ into analytic geometry. The first new phenomenon that one encounters in applying the six-functor formalism to solid quasi-coherent sheaves is that the kernel of the restriction functor $j^{*}:\QCoh{X_{\bs}} \rightarrow \QCoh{U_{\bs}}$ consists not just of solid quasi-coherent sheaves supported on the complement $Z_{\bs}$, but rather on its tubular neighborhood $\widetilde{Z}_{\bs}=\Specan{A^{\wedge}_{I}}$. Here, $A^{\wedge}_{I}$ denotes the completion of the ring $A$ at the ideal $I$, computed as a limit inside solid modules, together with its induced analytic ring structure over $A_{\bs}$.

Altogether, we obtain a recollement
\[
\xymatrix{
\QCoh{\widetilde{Z}_{\blox}}\ar@{^{(}->}[rr]|{{i_*}} && \QCoh{X_{\blox}} \ar@<2ex>[ll]^{i^{\:!}}\ar@<-2ex>[ll]_{i^*} \ar[rr]|{j^*} && \QCoh{U_{\blox}} \ar@<2ex>@{_{(}->}[ll]^{j_*}\ar@<-2ex>@{_{(}->}[ll]_{{j}_!}
}.
\]
and therefore, as in the topological case, a homotopy pullback square reconstructing any quasi-coherent sheaf $F$ on $X_{\bs}$ from its restrictions to the tubular neighborhood and its complement. Moreover, this recollement is continuous in the sense that all functors preserve (homotopy) colimits, thus solid quasi-coherent sheaves simultaneously enjoy the best properties both of sheaves on nice topological spaces and quasi-coherent sheaves on nice schemes.

The key idea of this paper is that the Beilinson-Parshin adelic decomposition of a quasi-coherent sheaf on a scheme $X$ of finite type over $\Z$ and of dimension $d$ arises by considering a sequence of recollements associated to a skeletal filtration 
\[
{Z}_{0} \subset {Z}_{1} \subset \cdots \subset {Z}_{d-1} \subset X_{\blox}
\]
where ${Z}_{k}$ is in a certain sense the union of all tubular neighbourhoods of closed subschemes of dimension at most $k$. Formally, we define 
$Z_{k}=\Specan{A^{\thicksim}_{S_k}}$, where $A^{\thicksim}_{S_k}$ as a solid ring is defined as the limit of all completions at all ideals of dimension at most $k$.

Given any filtration $Z_{0} \subset Z_{1} \subset \dots \subset Z_{d-1} \subset Z_{d}=X_{\blox}$, define strata $X_{k}=Z_{k} \setminus Z_{k-1}$ with inclusions $s_{k} : X_{k} \hookrightarrow X$. For each stratum we have a localization functor $\Psi_{k}:={s_{k}}_{*}s^{*}_{k} : \QCoh{X_{\blox}} \arr \QCoh{X_{\blox}}$ equipped with a unit of adjunction $\id \arr {s_{k}}_{*}s^{*}_{k}=\Psi_{k}$. For a tuple $d \geq k_{m} > \dots > k_{0} \geq 0$, we can form the iterated localisation functor $\Psi_{k_{0}} \circ \dots \circ \Psi_{k_{m}}$, and using the units of adjunction, the iterated localisations fit together into a commutative `fracture cube' of functors. For example, in the case $d=2$, the cube is of the form

\[
\xymatrix{
& \Psi_0\ar[ld]\ar[dd] && \id \ar[ll]\ar[dd]\ar[dl] \\
\Psi_0\Psi_1\ar[dd] && \Psi_1\ar[ll]\ar[dd]  \\
& \Psi_0\Psi_2\ar[ld] && \Psi_2\ar[ll]\ar[ld] \\
\Psi_0\Psi_1\Psi_2 && \Psi_1\Psi_2 \ar[ll]&.
}
\]
In any dimension, using the theory of total homotopy fibres \cite{2014chromatic}\cite{greenlees2020adelic}, or of non-commutative stratifications \cite{stratified}, the corresponding punctured cube exhibits a sheaf $F$ as the homotopy limit of the rest cube.

The main idea of this paper is to analyse the skeletal filtration $Z_{0} \subset Z_{1} \dots Z_{d-1} \subset Z_{d}=X_{\blox}$ of a $d$-dimensional finite type scheme and to show that the corresponding fracture cube for a solid quasi-coherent sheaf on $X$ gives a Beilinson-Parshin adelic decomposition expressed in terms of filtered colimits of iterated completed localizations over flags of integral closed subschemes. More precisely, our main theorem is Theorem \ref{thm-adeles}, a special case of which we state here for concreteness.

\begin{theorem}[Theorem \ref{thm-adeles}]\label{main-thm-intro}
Let $A$ be a commutative ring of finite type over $\mathbb{Z}$. Then the vertices of the fracture cube applied to $A$ are given as products over nested chains of prime ideals over completed localisations:
\[
L_{k_0}L_{k_1}\ldots L_{k_m}A \;\simeq\; \prod_{\substack{\p_m \\ \dim \p_m = k_m}} \leftcomp\: (\ldots \prod_{\substack{\p_1\supset \p_2\\\dim\p_1 = k_1}}\leftcomp\:(\prod_{\substack{\p_0 \supset \p_1\\\dim\p_0 = k_0}} \leftcomp A_{\p_0})_{\p_1} \ldots )_{\p_m} \;,
\] where the notation $\leftcomp\:(-)_\p$ means the functor that localizes at the prime $\p$ first and then completes at that prime.
Moreover, the fracture cube for a general solid module commutes with colimits and right t-bounded products.
\end{theorem}

In section \ref{ssec-classical-BP}, we discuss in more detail how the formula of Theorem \ref{main-thm-intro}
is actually related to the original formula of Beilinson-Parshin (see \cite{BeilinsonResidues}).

Let us emphasise again that our cubical decomposition works for a general solid quasi-coherent sheaf and is constructed from the geometry of the skeletal filtration, with the flags of closed integral subschemes and the completed localizations appearing a posteriori. The advantage of this point of view is that it suggests various generalizations, by decomposing solid quasi-coherent sheaves with respect to other naturally arising filtrations of schemes, for example the tangent rank stratification for a morphism or the orbit-closure stratification for a scheme with group action. We intend to pursue these ideas in future work.

\medskip
In fact, according to the formalism of \cite{stratified}, part of which we recall in section \ref{ssec-stratifications-and-sm-sp},
our skeletal filtration provides a decomposition of the whole category $\mod_{A_\bs}$ (or more generally
$\qcoh(X_\bs)$ for a finite type scheme $X$) into a cube shaped limit (see Theorem A in \cite{stratified}). However, the categories that appear in the limit do not admit a clear geometric
description and the functors between them are only lax symmetric monoidal.

Nonetheless, it is possible to use the skeletal filtration combined with some additional considerations to prove what we call solid adelic descent, which generalizes Theorem 3.1 of \cite{Groechenig_2017}, and in particular recovers the classical adelic formula of Weil describing the groupoid of vector bundles on a curve. We now would like to explain the statement of solid descent
in the case of the category $\mod_{\Z_\bs}$ of solid abelian groups. On the one hand, it immediately follows from Theorem
\ref{main-thm-intro} and the discussion that preceded it that the square \[
\xymatrix{\Z \ar[r] \ar[d] & \Q \ar[d] \\
\prod_{p} \Z_{p} \ar[r] & \Q \otimes_{\Z_\bs} \prod_{p} \Z_{p} }
\] is a (homotopy) pullback square in $\mod_{\Z_\bs}$. On the other hand, it follows from some basic properties of the
symmetric monoidal structure on $\mod_{\Z_\bs}$, see Lecture 6 in \cite{condensed}, that all of the algebras \[
\prod_p \Z_p \:, \quad \Q \:, \quad \Q \otimes_{\Z_\bs} \prod_p \Z_p \quad \in\; \calg(\mod_{\Z_\bs})
\] are idempotent. Using these two observations, one immediately gets the following statement.
\begin{theorem}
The square \[\xymatrix{
\mod_{\Z_\bs} \ar[rr]\ar[d] && \Q-\mod_{\Z_\bs} \ar[d] \\
\prod_p \Z_p - \mod_{\Z_\bs} \ar[rr] && \Q\otimes_{\Z_\bs}\prod_p \Z_p - \mod_{\Z_\bs} &,
}\] where all of the categories are just categories of modules in $\mod_{\Z_\bs}$ over the corresponding algebra and 
the functors are just the inductions along the maps of algebras,
is a pullback square in $\calg(\Pr^L_\st)$.
\end{theorem}

More generally, in the case of a finite type ring $A$ of dimension $d$, we have the skeletal filtration on $\anspec(A_\bs)$
(or on the category $\mod_{A_\bs}$), and any object $M \in \mod_{A_\bs}$ decomposes as a cubical limit:
\[
M \;\stackrel{\simeq}{\arr}\; \stackrel[\substack{[0 \le k_0 < k_1 < \ldots < k_m \le d] \\ 0\le m}]{}{\lim} L_{k_0}\ldots L_{k_m} M \;,
\]
where $L_k$ denotes the functor that localizes to the $k$-th stratum and embeds back. It is then quite formal and would hold
for any other filtration with similar formal properties that the functor \begin{equation}\label{intro-adelic-descent-functor}
\mod_{A_\bs} \;\arr\; \stackrel[\substack{[0 \le k_0 < k_1 < \ldots < k_m \le d] \\ 0\le m}]{}{\lim} L_{k_0}\ldots L_{k_m} A - \mod_{A_\bs}\;,
\end{equation}
whose components are given by the inductions along the maps of algebras
is fully faithful; the reader will find the details in the proof of Theorem \ref{adelic-cover}.

We then prove that the natural map $L_{k_0}A \otimes_{A_\bs} \ldots \otimes_{A_\bs} L_{k_m} A \rightarrow L_{k_0}\ldots L_{k_m} A$ is an isomorphism (Proposition \ref{factors-are-tensor-products}) and that $L_{k_0}\ldots L_{k_m} A$ is idempotent over $A_\bs$ (Corollary \ref{adelic-rings-are-idempotent}),
which in particular implies that the functor (\ref{intro-adelic-descent-functor}) is essentially surjective.
As opposed to the previous point, this is specific to the fact that we are dealing with the skeletal filtration.
It uses Theorem \ref{main-thm-intro} and some nice completeness properties of the solid tensor product. All that is explained
in detail at the beginning of section \ref{5th_sec} and in Appendix \ref{app-k-comp}.

Summing up, we have the following statement, which generalizes the finite type affine case of Theorem 3.1 of \cite{Groechenig_2017}
(we also treat the non-affine case but do not state it in the introduction).
\begin{theorem}[Solid adelic descent, Theorem \ref{adelic-cover}]\label{intro-adelic-descent-thm}
The functor
\[
\mod_{A_\bs} \;\arr\; \stackrel[\substack{[0 \le k_0 < k_1 < \ldots < k_m \le d] \\ 0\le m}]{}{\lim} L_{k_0}A \otimes_{A_\bs} \ldots \otimes_{A_\bs} L_{k_m} A - \mod_{A_\bs} \;, 
\]
whose components are given by the inductions along the maps of algebras, and which is in particular
a map in the category $\calg(\Pr^L_\st)$, is an equivalence. 
\end{theorem}

\medskip
Although we do not need it for the paper, it is worth mentioning that in terms of the theory of !-descent
introduced by Clausen and Scholze and explained in Lecture 18 of \cite{clausen2023analytic}, Theorem \ref{intro-adelic-descent-thm}
is an example of a closed !-cover. Namely, idempotency of the rings $L_{k}A$ together with the isomorphism
$L_{k_0}\ldots L_{k_m} A \;\simeq\; L_{k_0}A \otimes_{A_\bs} \ldots \otimes_{A_\bs} L_{k_m}A$ implies that the map \[
\anspec (L_0A, A_\bs) \;\bigsqcup\; \ldots \;\bigsqcup\; \anspec (L_d A, A_\bs) \;\arr\; \anspec(A_\bs)
\] is a !-closed cover with intersections given by the spectra of the adelic rings of higher length. 

\medskip
There is also a global version of solid adelic descent, which holds for a scheme $X$ of finite type over $\Z$.
In this case, the category $\qcoh(X_\bs)$ of solid quasi-coherent sheaves is decomposed as a cubical limit
into pieces of the form \[
L_{k_0}\ldots L_{k_m} \O_X - \mod(\qcoh(X_\bs)) \;.
\] Then we also prove---see Theorem \ref{thm-adeles-are-affinoid}---that each of these pieces is in fact
equivalent to a category of modules complete with respect to some analytic ring structure on the ring of global sections \[
\Ga( X_\bs , L_{k_0}\ldots L_{k_m} \O_X) \;\in\; \calg(\mod_{\Z_\bs}) \;.
\]

\medskip
\noindent \textbf{The paper is structured as follows.}

In the first section, we review some general solid algebra, the smashing spectrum, and the language of non-commutative stratifications.
More specifically, in section \ref{sec-gen-solid}, we review solid algebra and the language of analytic rings. In sections
\ref{ssec-completions} and \ref{ssec-completinos2}, we review certain nice properties of local cohomology and completions.
In section \ref{ssec-HN-thm}, we review the Hopkins-Neeman theorem and local cohomology of a specialization closed subset.
In sections \ref{ssec-tub-neig} and \ref{ssec-sm-sp}, we study the properties of the solid ring obtained from our ring $A$
by completion with respect to a specialization closed subset. In particular, we study a decomposition of $\mod_{A_\bs}$ into a recollement
with respect to a specialization closed subset, and also prove that completing the ring $A$ at a specialization closed subset
provides a nice map from the poset opposite to that of specialization closed subsets to the smashing spectrum of $\mod_{A_\bs}$.
In section \ref{ssec-gen-closed-and-ultrasolid}, we study the other side of the recollement, which is given by completing at a
generalization closed subset and is closely related to ultrasolid analytic ring structures.
In section \ref{ssec-stratifications-and-sm-sp}, we review the language of non-commutative
stratifications of \cite{stratified} and prove that a symmetric monoidal stable presentable category
is stratified over its smashing spectrum. Combined with the results of section \ref{ssec-sm-sp}, that tells us that
the category $\mod_{A_\bs}$ is stratified over the poset of specialization closed subsets of $\spec(A)$, which we later use
to define the skeletal filtration and construct the Beilinson-Parshin resolution. In section \ref{ssec-cubical-resolutions}, following \cite{antolíncamarena2014chromaticfracturecubes} and independent of \cite{stratified}, we
give a more elementary construction of cubical resolutions resulting from a finite length filtration by closed subcategories.

In the second section, we construct and analyze the skeletal filtration of a finite type affine scheme and prove the affine case of our main 
theorem---Theorem \ref{thm-adeles}. In section \ref{ssec-classical-BP}, we also briefly review the relation to the classical formula of 
Beilinson-Parshin.

In the third section, we explain how to globalize our construction to non-affine schemes. In particular, we
globalize the results of section \ref{ssec-sm-sp}.

In the final section, we state and prove the solid adelic descent theorem; part of the proof
is postponed until Appendix \ref{app-k-comp}.

\ssec{Notation and conventions}

\sssec{}
In this paper, we work in the setting of $\infty$-categories. We adopt the usual convention that the word category means
$(\infty,1)$-category, and all the standard categorical notions and constructions, such as functors and (co)limits, should be
understood in this context. In particular, we shall typically refer to derived completion simply as completion.

\sssec{} Following \cite{Lurie-HTT}, Definition 5.3.4.5, we call a functor continuous, if preserves filtered colimits.
Since we will be working mostly with stable categories and exact functors, continuous will usually be equivalent 
to preserving all colimits.

\sssec{}
We will denote by $\Pr^L_\st$ the $(\infty,1)$-category of stable presentable categories and exact continuous functors.

\sssec{}\label{rmk-evrth-works-over-a-field}
Throughout the paper, we work with commutative, unital rings $A$ of finite type over $\Z$.

We chose to work over $\Z$ for convenience, but
the first three sections of this paper should go through without change if one chooses to work in finite type over a field $k$ using
the category of ultrasolid $k$-modules, defined and studied in \cite{aparicio2024ultrasolidhomotopicalalgebra}.
Moreover, the arguments of the last section and the appendix
will also go through if one is willing to assume the base field to be countable. 

In the case of a finite dimensional Noetherian scheme, a natural approach would be to first define the category of ultrasolid modules
over a Noetherian ring and then glue those categories over a Noetherian scheme. While the first step doesn't seem very hard,
and thus the arguments of the first three sections of this paper should work for a finite dimensional affine Noetherian scheme,
the second step seems to require developing some more theory that is beyond the scope of this paper.

Another potential approach seems to be to work with filtered colimits of completions to define tubular neighborhoods. However, this would seem to produce a variant of Beilinson-Parshin adeles, in which the local factors are expressed not as completed localizations of general Noetherian rings, but as filtered colimits of completed localizations of finite type subrings.

\sssec{}
Let $Z \subseteq \spec(A)$ be a closed subset defined by an ideal $I \subset A$. We will use $A^\comp_I$
to denote the (derived) completion of $A$ at the ideal $I$. 
The same goes for the completion of a (solid) module over $A$.
We note that the functor of completion $(-)^\comp_I$ at the ideal $I \subset A$ only depends on the closed subset $Z \subseteq \spec(A)$, but
not on the ideal $I$ cutting out $Z$. To emphasize that, and for other reasons, we will often call the functor of completion
at $I$ the functor of completion at $Z$ and denote it $(-)^\thicksim_Z$.

\sssec{}
More generally, for a specialization closed subset $T \subseteq \spec(A)$ (which can be thought of as a subset equal to the union
of all its subsets which are closed in $\spec(A)$), we will define a completion at $T$ and denote it $(-)^\thicksim_T$,
see section \ref{ssec-HN-thm} and its subsection \ref{sssec-completion-at-a-specialization-closed}.

\sssec{}
Similarly, for a generalization closed subset $E \subset \spec(A)$, we will define a completion at $E$ and denote it
$(-)^\thicksim_E$, see section \ref{ssec-gen-closed-and-ultrasolid}. This will turn out to be equivalent to the completion
with respect to the ultrasolid analytic ring structure on the ring \[
\Ga(E, \O) \;\stackrel{\define}{=}\; \stackrel[U \supseteq E]{}{\colim} \Ga(U, \O) \;,
\] where the colimit is taken over all Zariski open subsets $U \subseteq \spec(A)$ containing $E$ and ordered by (opposite) inclusion.
See Definition \ref{def-ultrasolid} for the definition of the ultrasolid analytic ring structure.

We believe that this notation does not clash with our notation for the completion at a specialization closed subset because
a subset of $\spec(A)$ cannot be both specialization and generalization closed except some trivial cases, in which the functors agree.

\sssec{}
For a prime ideal $\p \subset A$, we will denote by $\leftcomp A_\p$ the local ring $A_\p$ completed at the maximal ideal.
Similarly, for any (solid) $A$-module $M$, $\leftcomp\: M_\p$ denotes the localization of $M$ at $\p$ followed by the (derived) completion
at the maximal ideal of $A_\p$. This piece of notation makes formulas involving iterated completed localizations---like
the formula of Theorem \ref{thm-adeles}---less cumbersome.

\sssec{}
We use Lurie's notation for categories of modules over an algebra as opposed to the traditional derived category notation. For example,
$\mod_A$ means the whole stable $(\infty, 1)$-category of modules over the algebra $A$. Likewise, $\mod_{A_\bs}$ denotes
the stable $(\infty, 1)$-category of solid $A$-modules, and we would like to bring it to the reader's attention that
Clausen and Scholze adopt different conventions in \cite{condensed}.

This also applies to categories of quasi-coherent sheaves over a scheme: $\qcoh(X)$ means the stable $(\infty,1)$-category of
quasi-coherent complexes over $X$; $\qcoh(X_\bs)$ means the stable $(\infty,1)$-category of solid quasi-coherent sheaves over $X$.

\sssec{}
Given an analytic ring $\mathcal{B}$, we will use both $\mathcal{B}-\mod$ and $\mod_{\mathcal{B}}$ to denote the category
of complete $\mathcal{B}$-modules. There should usually be no room for confusion.

\sssec{}
To make our notation less cumbersome,
given a commutative algebra object $B \in \calg(\mod_{A_\bs})$ in the symmetric monoidal category $\mod_{A_\bs}$, we will
often use $B - \mod_{A_\bs}$ (as opposed to something like $B-\mod(\mod_{A_\bs})$) to denote the category of $B$-modules inside
$\mod_{A_\bs}$.

\ssec{Acknowledgments}
We would like to thank Sergei Gorchinsky, Dmitry Kaledin, Mikhail Kapranov, Nick Rozenblyum, Denis Osipov, and Artem Prikhodko for conversations about 
adeles and for their interest in this project.
Special thanks are due to A. Beilinson and J. Campbell for pointing out errors and suggesting corrections, especially to J. Campbell
for letting us know that our proof of Proposition 1.1.2 in the previous version was incorrect.
Finally, we would like to thank the anonymous referee for lots corrections and suggested improvements, especially for suggesting
to use the smashing spectrum explicitly and for pointing out that our proofs of Proposition \ref{k-comp-is-monoidal} and Corollary
\ref{appB-sing-ideal-comp-is-monoidal} do not hold without the connectivity assumption.
Part of the research was done in the Centre of Pure Mathematics within MIPT, grant number FSMG-2023-0013.
C. Brav was partially supported by a Junior Leader grant from the Basis Foundation and is grateful to the Sino-Russian Mathematics Center at 
Peking University for hospitality and excellent working conditions while preparing this paper.
G. Konovalov would like to thank the Shanghai Institute for Mathematical and Interdisciplinary Sciences for hospitality while preparing the revised version of this paper.

%% file: general_solid.tex
\section{Solid algebra and tubular neighbourhoods}

\newcommand{\freefinlen}{\mathrm{free,\:fin\text{-}len}}

In this section, we review some basic constructions about solid modules and the language of analytic rings, following \cite{condensed}
and \cite{analytic}; study local cohomology and completion, introduce the recollement of the 
category $\mod_{A_\bs}$ of solid $A$-modules associated to the tubular neighborhood of a specialization closed subset of $\spec(A)$
and its complement; introduce a map from the poset of specialization closed subsets into the smashing spectrum of $\mod_{A_\bs}$;
review the language of non-commutative stratifications by \cite{stratified}, and review a special case of a stratification by
a finite linearly ordered set (which we usually call a filtration), which gives rise to a cubical resolution.

\ssec{Generalities on solid modules}\label{sec-gen-solid}

Consider a static (i.e. sitting in degree 0 with respect to the t-structure)
commutative ring $A$ of finite type over $\Z$. The category of discrete $A$-modules
(i.e. the category of $A$-modules in the stable symmetric monoidal category of usual $\Z$-modules) we denote by $\mod_A$.
Clausen and Scholze (\cite{condensed}, Lectures 5,6 and 8) have defined the category $\mod_{A_\bs}$ of solid 
$A$-modules. In case $A = \Z$ this category is called the category of solid abelian groups (see Lectures 5 and 6 in \cite{condensed})
and serves as our base cateogory in this paper.
$\mod_{A_\bs}$ admits all limits and colimits and a continuous tensor product, denoted $\otimes_{A_\bs}$, along with a fully-faithful embedding
\[
\iota \colon \mod_A \hookrightarrow \mod_{A_\bs} \;,
\]
which is continuous and symmetric monoidal, but in general does not preserve general infinite limits,
in particular does not preserve infinite products.
The category $\mod_{A_\bs}$ admits compact generators given by objects of the form
of the form $\prod_I A$ for $I$ and index set, where $A = \iota(A)$ is the monoidal unit, the product is computed
inside $\mod_{A_\bs}$ and is not necessarily
finite. On the compact generators, the tensor product and the internal $\hom$ are computed by the following formulas:
\begin{equation}\label{gen-sol-basic-properties}
\begin{split}
&\;\prod_I A \otimes_{A_\bs} \prod_J A \;\simeq\; \prod_{I\times J} A \;;\\
&\;\uhom_{A_\bs} (\prod_I A, \prod_J A) \;\simeq\; \prod_J \bigoplus_I A \;.
\end{split}
\end{equation}

\sssec{Synthetic description of $\mod_{A_\bs}$}\label{gen-solid-synth}
Let $\mod_A^{\projfinlen}$ denote the full subcategory of $\mod_A$ spanned by finite complexes of projective $A$-modules of
not necessarily finite rank. There is a lax symmetric monoidal functor \begin{equation}\label{gen-solid-synthetic-functor}\xymatrix{
(\mod_A^\projfinlen)^\op \ar[r] & \mod^\w_{A_\bs} \\
F \ar@{}[u]|\invertin \ar@{|->}[r] & \uhom_{A_\bs} (\iota(F), A) \ar@{}[u]|\invertin &.
}\end{equation}
\begin{lemma}\label{lem-solid-synth-cmp-obj}
The functor (\ref{gen-solid-synthetic-functor}) is fully faithful and symmetric monoidal.
The image of the composite \[\xymatrix{
(\mod_A^\projfinlen)^\op \ar[r] & \mod^\w_{A_\bs} \ar@{^(->}[r] & \mod_{A_\bs}
}\] generates the category $\mod_{A_\bs}$.
\end{lemma}
\begin{proof}
Since the source category is generated under finite limits and retracts by objects of the form $\prod_I A$, both fully faithfulness and 
symmetric monoidality follow from the basic properties (\ref{gen-sol-basic-properties}). The last claim follows from the fact that
the category $\mod_{A_\bs}$ is compactly genetated by products of copies of  $A$.
\end{proof}

\sssec{Recollection on analytic rings}
Here we, following Lecture 12 of \cite{analytic}, briefly recall the notion of an analytic ring. Since working relative to the category
$\mod_{\Z_\bs}$ of solid abelian groups (see Lectures 5 and 6 of \cite{condensed}) is sufficient for the
purposes of this paper, we will give the definition relative to $\mod_{\Z_\bs}$.
\begin{definition}(\cite{analytic}, Proposition 12.20)\label{def-analytic-ring}
Let $R \in \mathrm{CAlg}^{\le 0}(\mod_{\Z_\bs})$ be an animated commutative ring
in the category of solid abelian groups. An analytic ring structure on $R$ is a full subcategory
\[\xymatrix{\DD^{\le 0} \ar@{^{(}->}[rr]^(.4){j_*} && R-\mod^{\le 0}_{\Z_\bs}}\]
such that it's closed under limits and colimits, stable under the functor $\uhom_{\mod^{\le 0}_{\Z_\bs}}(\prod_I \Z, -)$ for any
compact generator $\prod_I \Z$, and the inclusion admits a left adjoint $j^*$ satisfying
\[ R \stackrel{\simeq}{\arr} j_*j^* R \;.\]
\end{definition}
We note that we will be mostly dealing with static analytic rings, i.e. with those whose underlying solid ring satisfies
$R \in \mathrm{CAlg}^{\heartsuit}(\mod_{\Z_\bs})$.

One of the main examples of an analytic ring is the following.
For a commutative finite type algebra (in usual abelian groups) $A$, the map $\Z \arr A$ induces the restriction functor
\begin{equation}\label{local-gen-solid-eq-1}
\mod_{A_\bs} \arr A-\mod_{\Z_\bs} \;.
\end{equation}
\begin{theorem}[8.1 in \cite{condensed}]
The functor (\ref{local-gen-solid-eq-1}) defines an analytic ring structure on the ring $A$.
\end{theorem}

Although we will mostly be dealing with static analytic rings, sometimes non-static and even non-connective analytic rings
will naturally appear as output of some our constructions. We therefore would like to give a more general version of
Definition \ref{def-analytic-ring}, cf. Definition 2.1.1 \cite{camargo2024analyticrhamstackrigid}.
\begin{definition}\label{defn-analytic-Eoo-ring}
Let $R$ be an $\mathbb{E}_\infty$-ring
in the category of solid abelian groups. An analytic ring structure on $R$ is a full subcategory
\[\xymatrix{\DD \ar@{^{(}->}[rr]^(.4){j_*} && R-\mod_{\Z_\bs}}\]
such that it's closed under limits and colimits, stable under the functor $\uhom_{\Z_\bs}(\prod_I \Z, -)$ for any
compact generator $\prod_I \Z$, and the inclusion admits a left adjoint $j^*$ satisfying
\[ R \stackrel{\simeq}{\arr} j_*j^* R \;.\]
\end{definition}

\sssec{Ultrasolid analytic ring structures}
Let $R$ be a solid $\mathbb{E}_\infty$ (or animated) ring. Let \[
\mathcal{D}_R \subset R-\mod_{\Z_\bs} 
\] denote the full subcategory closed under finite colimits and retracts and containing objects $\prod_I R$.
We get a continuous
functor \begin{equation}\label{def-ultrasolid-oblv-functor}
\ind(\mathcal{D}_R) \arr R-\mod_{\Z_\bs} \;.
\end{equation}
\begin{definition}\label{def-ultrasolid}
In case the functor (\ref{def-ultrasolid-oblv-functor}) is fully faithful and satisfies the conditions of
Definition \ref{defn-analytic-Eoo-ring}, and therefore corresponds to an analytic ring structure on the ring $R$,
we call this analytic ring structure ultrasolid (cf. Example 2.6.7 in \cite{camargo2024analyticrhamstackrigid})
and denote $R_\bullet$. The full category
$\ind(\mathcal{D}_R) \subseteq R-\mod_{\Z_\bs}$ is called the category of ultrasolid $R$-modules and, with accordance
to our convections, can denoted by both $R_\bullet-\mod$ and $\mod_{R_{\bullet}}$.
\end{definition}

\sssec{Solid tensor product and limits}
Another technical observation, which we would like to leave here for future reference, concerns the solid tensor product and certain types of limits preserved by 
it. It is one of the basic properties, listed at the beginning of \ref{sec-gen-solid},
that tensoring with a compact object preserves products of copies of $A$, but we will need a stronger version of
that property and we package it into the following lemma.
\begin{lemma}\label{lem-tensor-vs-prods}
For a compact object $C \in \mod^\omega_{A_\bs}$, the functor
\[
C \otimes_{A_\bs} - : \mod_{A_\bs} \longrightarrow \mod_{A_\bs}
\] preserves right t-bounded products. That is, if $\{M_j\}_{j\in J}$ is a collection of $A_\bs$-modules uniformly cohomologically bounded from above, then
the natural map
\[
C \otimes_{A_\bs} \prod_J M_j \stackrel{\simeq}{\longrightarrow} \prod_J \bigl( C\otimes_{A_\bs} M_j \bigr)
\] is an isomorphism.
\end{lemma}
Before we get to the proof of Lemma \ref{lem-tensor-vs-prods}, recall from \cite{condensed}, Theorem 2.2, that the abelian heart of the category $\mod_{\Z_\bs}$
of solid abelian groups satisfies Grothendieck's axiom (AB6). For the sake of clarity, we would like to state and prove a variant of that axiom that deals with
objects inside the derived categories. It will follow from the abelian variant combined with some basic properties of the categories and the t-structures.
\begin{lemma}\label{lem-ab6-for-complexes}
For any index set $J$ and filtered categories $I_j$, $j \in J$, with functors $i_j \mapsto M_{i_j}$ from $I_j$ to $\mod_{A_\bs}$,
the map\[
\stackrel[\prod_j I_j]{}{\colim} \prod_{j\in J} M_{i_j} \arr \prod_{j\in J} \stackrel[I_j]{}{\colim} M_{i_j}
\] is an isomorphism
\end{lemma}
\begin{proof}
Since the standard t-structure on $\mod_{A_\bs}$ is both left and right complete, it suffices to check that the map induces an isomorphism on cohomology.
The latter reduces to the abelian (AB6) by the exactness of products and filtered colimits, which also follows from Theorem 2.2 of \cite{condensed}. 
\end{proof}
\begin{proof}[Proof of Lemma \ref{lem-tensor-vs-prods}]
Writing each $M_j$ as a filtered colimit of compact objects inside $\mod_{A_\bs}^{\le n}$
and pulling these colimits out of the product by Lemma \ref{lem-ab6-for-complexes},
we reduce the statement to the case where each $M_j$ is compact.

In this case, we
present each $M_j$ as a finite length complex whose terms are non-zero only in degrees not greater than $n$ and are given by products of
copies of $A$, thereby
presenting $\prod_J M_j$ as a (possibly infinite to the left) complex whose terms are products of copies of $A$.
Using that presentation, we filter the object $\prod_J M_j$ using the naive truncations. Namely, there is a map \[
\stackrel[m \ge 0]{}{\colim} \prod_J \sigma_{\ge -m} M_j \;\arr\; \prod_J M_j \;,
\] which one easily checks to be an isomorphism.
Indeed, by the completeness of the t-structure, it suffices to see that this map induces an isomorphism
on cohomology, which immediately follows from the t-exactness of products and filtered colimits.
We also present the object $C$ as a bounded above complex whose terms are products of copies of $A$, and filter it
by the naive truncations: \[
C \;\simeq\; \stackrel[n \ge 0]{}{\colim} \sigma_{\ge -n} C \;.
\]

It now suffices to prove that the map \[
\sigma_{\ge -n} C\otimes_{A_\bs} \prod_J \sigma_{\ge -m} M_j \;\arr\; \prod_J \bigl( (\sigma_{\ge -n}C) \otimes_{A_\bs} (\sigma_{\ge -m}M_j) \bigr)
\] is an isomorphism, which is basically the fact that products distribute over the solid tensor product in case both of the arguments are
products of copies of $A$---see (\ref{gen-sol-basic-properties}).
\end{proof}

%% file: some_algebra.tex
\ssec{Completions}\label{ssec-completions}

In this subsection,
we recall how one usually defines completion $A^\comp_I$ of $A$ at an ideal $I \subset A$, and outline some favorable properties of $A^\comp_I$ as a
commutative algebra in $\mod_{A_\bs}$. We are keeping the assumption that the ring $A$ is of finite type over $\Z$.

It might be useful to keep in mind the following intuitive picture:
$A^\comp_I \otimes_{A_\bs} A^\comp_I \simeq A^\comp_I$,
and therefore $\spec^{an}(A^\comp_I) \subseteq \spec^{an}(A_\bs)$ looks like a `Zariski' open subset,
which should be thought of as a kind of a solid tubular neighborhood of $\spec(A/I)$ inside $\spec(A)$. On the other hand,
$A^\comp_I$ is compact as a solid $A$-module, which implies some additional functoriality between sheaves on $\spec(A^\comp_I)$ and $\spec(A_\bs)$---specifically,
the embedding $A^\comp_I - \mod_{A_\bs} \subseteq \mod_{A_\bs}$ admits, besides a left adjoint, a continuous right adjoint. Thus from the point of view of the six-functor formalism, $\spec^{an}(A^\comp_I) \subseteq \spec^{an}(A_\bs)$ behaves like a closed subset, and we shall tend to favor this point of view.

\begin{definition}\label{def-completion}
Let $Z \subset \spec(A)$ be a closed subset, let $I_Z \subset A$ denote its ideal, and let $Z^c = \spec(A) \setminus Z$ denote
the complement. Let
\[
\Gamma_Z \stackrel{\define}{=} \mathrm{fib}\bigl( A \;\to\; \Gamma(Z^c, \O)\bigr) \;,
\]
be local cohomology of $A$ with supports on $Z$.
Recall from \cite[\href{https://stacks.math.columbia.edu/tag/0A6R}{Lemma 0A6R}]{stacks-project}
that it can be represented by an extended alternating \v{C}ech complex:
let $f_1 \ldots f_l \in A$ be a set of generators of the ideal $I_Z$, then
\begin{equation*}
\begin{split}
\Gamma_Z = \Bigl(\Gamma(Z^c, \O) \:/\: A \Bigr)[-1] &\;\simeq\; \Bigl( A\bigl[\frac{1}{f_1}\bigr] / A \Bigr) [-1] \otimes_{A} \ldots \otimes_{A} \Bigl( A\bigl[\frac{1}{f_l}\bigr] / A \Bigr) [-1] \;\simeq\; \\
&\;\simeq\; \Bigl(\:  A \to \bigoplus_i A\bigl[\frac{1}{f_i}\bigr] \to \bigoplus_{i < j} A\bigl[\frac{1}{f_if_j}\bigr] \to \ldots \to A\bigl[\frac{1}{f_1 \ldots f_l}\bigr]  \:\Bigr) \;\;.
\end{split}
\end{equation*}
Completion of $A$ at the ideal $I_Z$, denoted $A^\comp_{I_Z}$, can be expressed by the formula
\[
A^\comp_{I_Z} \:\stackrel{\define}{=}\: \uhom_{A_\bs} \bigl(\Gamma_Z , A\bigr) \;\;.
\]
More generally, completion of a module $M \in \mod_{A_\bs}$ is defined by the formula \[
M^\comp_{I_Z} \:\stackrel{\define}{=}\: \uhom_{A_\bs} \bigl(\Gamma_Z , M\bigr) \;.
\] We note that the functor of completion at $I_Z$ only depends on the closed subset $Z \subseteq \spec(A)$ rather than on a scheme
structure one may put on $Z$. To emphasize that and for other reasons to become apparent later, we will call this functor
the functor of completion at $Z$ and will denote it $(-)^\thicksim_Z$. 
\end{definition}
Let us emphasize that, per our conventions on categories and functors, the above notion of completion is derived, not naive.

\sssec{}\label{sssec-completion-is-idempotent-and-compact}
For future reference, we would like to recall that, for a closed subset $Z \subseteq \spec(A)$, the object \[
A / \Ga_Z \;\simeq\; \Ga(Z^c, \O)
\] is naturally an idempotent algebra inside $\mod_A^\projfinlen$. This follows immediately from
Lemma \ref{lem-local-coh-is-coidempotent}, Proposition \ref{prop-idem-obj-vs-idem-algs}, and Proposition \ref{prop-idem-vs-coidem},
all recorded just below.

Indeed, let $\CC$ denote a stable symmetric monoidal category;
$\CC = \mod^\projfinlen_A$ is the example most relevant to this section.
Recall from \cite{smashing_spectrum}, Definition 2.9, that an object $E \in \CC_{/1_\CC}$ is called a
\textit{coidempotent object} if the map \[
E \otimes E \;\arr\; 1_\CC \otimes E
\] is an isomorphism. The full subcategory of $\CC_{/1_\CC}$ spanned by coidempotent objects is denoted by $\cidem(\CC)$.
\begin{lemma}\label{lem-local-coh-is-coidempotent}
The object $\Ga_Z \in (\mod^\projfinlen_A)_{/A}$ is coidempotent.
\end{lemma}
\begin{proof}
Immediately follows from the presentation \[
\Ga_Z \;\simeq\; \Bigl( A\bigl[\frac{1}{f_1}\bigr] / A \Bigr) [-1] \otimes_{A} \ldots \otimes_{A} \Bigl( A\bigl[\frac{1}{f_l}\bigr] / A \Bigr) [-1]
\] of the local cohomology in terms of the generators of the ideal defining $Z \subseteq \spec(A)$.
\end{proof}
Dually, an object $Q \in \CC_{1_{\CC}/}$ is called an idempotent object if the map \[
1_\CC \otimes Q \;\arr\; Q \otimes Q
\] is an isomorphism. The full subcategory spanned of $\CC_{1_\CC / }$ by idempotent objects is denoted by $\idem(\CC)$. 
Let us also recall that a commutative algebra object $C \in \mathrm{CAlg}(\CC)$ is called
an \textit{idempotent algebra} if the multiplication map \[
C\otimes C \;\arr\; C
\] is an isomorphism.
The following is Proposition 4.8.2.9 of \cite{Lurie-HA}.
\begin{proposition}\label{prop-idem-obj-vs-idem-algs}
The composite \[
\mathrm{CAlg}(\CC) \;\simeq\; \mathrm{CAlg}(\CC)_{1_\CC / } \;\arr\; \CC_{1_\CC /}
\] induces an equivalence of the full subcategory spanned by idempotent algebras with the full subcategory spanned
by idempotent objects.
\end{proposition}
\begin{corollary}\label{cor-idem-is-poset}
The category $\idem(\CC)$ is equivalent to a poset, i.e. the mapping space between a pair of objects is either empty or contractible.
\end{corollary}
\begin{proposition}[Proposition 2.14,\cite{smashing_spectrum}]\label{prop-idem-vs-coidem}
The functor \[\xymatrix{
\cidem(\CC) \ar[rr] && \idem(\CC) \\
E \ar@{}[u]|\invertin \ar@{|->}[rr] && \mathrm{cof}(E \to 1_\CC) \ar@{}[u]|\invertin
}\] is an equivalence.
\end{proposition}

\sssec{}
The coidempotency of the local cohomology $\Ga_Z$ together with Lemma \ref{lem-solid-synth-cmp-obj},
which says that the functor \[\xymatrix{
\mod_A^\projfinlen \ar[rr] && \mod^\omega_{A_\bs} \\ 
P \ar@{|->}[rr]\ar@{}[u]|\invertin && \uhom_{A_\bs}(P, A) \ar@{}[u]|\invertin
}\] is symmetric monoidal, imply the following lemma.
\begin{lemma}\label{lem-comp-is-cmp-and-idmp}
The completion \[
A^\thicksim_Z \;=\; \uhom_{A_\bs}(\Ga_Z, A)
\] is an idempotent algebra inside $\mod_{A_\bs}^\w$.
\end{lemma}
\begin{proof}
By Proposition \ref{prop-idem-obj-vs-idem-algs}, it suffices to prove that the map \[
A \otimes_{A_\bs} A^\thicksim_Z \;\arr\; A^\thicksim_Z \otimes_{A_\bs} A^\thicksim_Z
\] is an isomorphism, which, by Lemma \ref{lem-solid-synth-cmp-obj}, is dual to the map \[
\Ga_Z \otimes_A \Ga_Z \;\arr\; A \otimes_A \Ga_Z \;.
\] The latter is an isomorphism by Lemma \ref{lem-local-coh-is-coidempotent}.
\end{proof}

\sssec{}\label{sssec-local-coh-functorial}
It follows from the discussion in the section \ref{sssec-completion-is-idempotent-and-compact}
(specifically, from Lemma \ref{lem-local-coh-is-coidempotent}, Corollary \ref{cor-idem-is-poset},
and Proposition \ref{prop-idem-vs-coidem})
that, for a pair of closed subsets $Z, Z^\prime$, the space of maps
$\map_{(\mod_A)_{/A}}(\Ga_{Z^\prime}, \Ga_{Z})$ is either empty or contractible. Moreover, it is non-empty if and only if
the map \begin{equation}\label{local-eq-maps-between-local-coh}
\Ga_{Z} \otimes_A \Ga_{Z^\prime} \;\arr\; A \otimes_A \Ga_{Z^\prime}
\end{equation}
is an isomorphism, which happens exactly when there is an inclusion $Z^\prime \subseteq Z$. Indeed, because $\Ga_Z$ splits as a tensor
product over the ideal generators, it suffices to treat the case where the ideal of $Z$ is generated by a single element: \[
\Ga_Z \;\simeq\; \Bigl( A\bigl[\frac{1}{f}\bigr] / A \Bigr) [-1] \;.
\]
In this case, the map (\ref{local-eq-maps-between-local-coh}) is an isomorphism if and only if
$\Ga_{Z^\prime} \otimes_A A\bigl[\frac{1}{f}\bigr] \;\simeq\; 0$, which happens exactly when $Z^\prime$ is contained in $\{f=0\}$.

\ssec{Tensoring with a compact preserves completions}\label{ssec-completinos2}

There is a variant of Lemma \ref{lem-tensor-vs-prods}, which we will need later.
We record it here for future reference.
\begin{lemma}\label{lem-tensor-pres-comp}
For an ideal $I \subset A$ and a pseudo-compact object $C \in \mod_{A_\bs}$,
the functor $C\otimes_{A_\bs}-$ preserves completion at $I$ of a right t-bounded object:
\[
C\otimes_{A_\bs} M^\comp_I \stackrel{\simeq}{\longrightarrow} \bigl( C\otimes_{A_\bs}M \bigr)^\comp_I \;,
\] where $M \in \mod_{A_\bs}^{\le n}$ for some $n$.
\end{lemma}

\begin{proof}
Let us pick a set of generators $f_1, \ldots, f_r \in I$ of the ideal $I$, and let $K_n$ denote
the Koszul complex of $A$ on the sequence $f_1^n, \ldots, f_r^n$. According to
\cite[\href{https://stacks.math.columbia.edu/tag/0913}{Tag 0913}]{stacks-project}, we can rewrite the local cohomology of $A$ with supports in
$V(I) \subseteq \spec(A)$, which is predual to the completion $A^\comp_I$, as a sequential colimit of duals to the Koszul complexes:
\[
\Gamma_{V(I)} \;\simeq\; \stackrel[n]{}{\colim} K^\vee_n \;.
\]
Any sequential colimit, such as $\stackrel[n]{}{\colim} K^\vee_n$, can be rewritten in terms of a cofiber and direct sums:
\[
\stackrel[n]{}{\colim} K^\vee_n \;\simeq\; \mathrm{cof}\Bigl( \oplus_n K_n^\vee \longrightarrow \oplus_n K_n^\vee \Bigr) \;,
\]
which implies that the functor of completion $\uhom_{A_\bs}(\Gamma_I, -)$ can be rewritten in terms of a fiber and products:
\[
\uhom_{A_\bs} (\stackrel[n]{}{\colim} K^\vee_n, - ) \;\simeq\; \mathrm{fib}\Bigl( \prod_n (K_n\otimes_{A_\bs}-) \: \longrightarrow \: \prod_n (K_n\otimes_{A_\bs}-) \Bigr) \;.
\] Now the claim follows from Lemma \ref{lem-tensor-vs-prods}.
\end{proof}

\ssec{Hopkins-Neeman theorem and specialization closed subsets}\label{ssec-HN-thm}

In this section we recall the classification of (co)idempotent (co)algebras in $\mod_A$,
which is originally due to Hopkins and Neeman.
We use \cite{Lurie-SAG}, section 2.6.1, as our reference. As usual, we will be dealing with a finite type commutative ring $A$.

Recall the notion of a coidempotent object from section \ref{sssec-completion-is-idempotent-and-compact}
(or from \cite{smashing_spectrum}, Definition 2.9). We denote by $\cidem(\mod_{A}) \subseteq (\mod_A)_{/A}$ the full subcategory
spanned by coidempotent objects, which is equivalen to a poset by Corollary \ref{cor-idem-is-poset} and
Proposition \ref{prop-idem-vs-coidem}.
There is a classical theorem by Hopkins-Neeman describing that poset, which we now recall.

A subset $T \subseteq \spec(A)$ is called specialization closed if, whenever $T$ contains a point $\p$,
$T$ also contains the whole closure $\overline{\p}$. Each such specialization closed subset is in fact the union of all closed
subsets contained in it. We denote by $\spec(A)^\spc$ the poset of specialization closed subsets of $\spec(A)$
ordered by the inclusion. Dually, the complement $\spec(A) \setminus T$ to a specialization closed subset $T$ is generalization
closed, i.e. with every point $\p \in \spec(A) \setminus T$ it contains every point $\q$ specializing to $\p$. We denote by
$\spec(A)^\gen$ the poset of generalization closed subsets ordered by inclusion.
Note that the posets $\spec(A)^\spc$ and $\spec(A)^\gen$ are anti-equivalent to each other: \[\xymatrix{
(\spec(A)^\spc)^\op \ar[rr]^\simeq && \spec(A)^\gen \\
T \ar@{}[u]|\invertin \ar@{|->}[rr] && T^c \;=\; \spec(A)\setminus T \ar@{}[u]|\invertin \;.
}\]

\begin{proposition}[cf. Theorem 2.6.1.2, \cite{Lurie-SAG}]\label{prop-hopkins-neeman}
\begin{enumerate}
\item The object \[
\Ga_T \;\stackrel{\define}{=}\; \stackrel[Z \subseteq T]{}{\colim} \Ga_Z \;,
\] where the colimit is taken over the poset of all closed subsets contained in $T$ ordered by inclusion
(see section \ref{sssec-local-coh-functorial} for the discussion on functoriality of $Z \mapsto \Ga_Z$),
is a coidempotent object inside $\mod_A$. Moreover, the support of $\Ga_T$ is equal to $T$, and \[\begin{split}
A_\p/\p \otimes_A \Ga_T &\;\simeq\; A_\p/\p \;, \text{ if } \p \in T \;. 
\end{split}\]
\item Dually, the algebra \[
A / \Ga_T \;\simeq\; \stackrel[Z \subseteq T]{}{\colim} \Ga(Z^c, \O) \;=\; \Ga(T^c, \O) 
\] is idempotent. Its support is equal to $T^c = \spec(A)\setminus T$, and \[\begin{split}
A_\p/\p \otimes_A \Ga (T^c, \O) &\;\simeq\; A_\p/\p \;, \text{ if } \p \in T^c \;. 
\end{split}\]
\item The maps \[\xymatrix{
U \ar@{}[d]|\invertinn \ar@{|->}[rr] && \Ga(U, \O) \ar@{}[d]|\invertinn \\
(\spec(A)^\gen)^\op \ar@{}[d]|\invsimeq \ar[rr] && \idem(\mod_A) \ar@{}[d]|\invsimeq \\
\spec(A)^\spc \ar[rr] && \cidem(\mod_A) \\
T \ar@{}[u]|\invertin \ar@{|->}[rr] && \Ga_T \ar@{}[u]|\invertin
}\] are isomorphisms of posets.
\end{enumerate}
\end{proposition}

There are a few immediate corollaries, which we record here for future use.
\begin{lemma}\label{lem-sm-sp-map-preserves-finite-joins}
For a pair of specialization closed subsets $T_1, T_2 \in \spec(A)^\spc$, the natural map
\begin{equation}\label{local-sm-sp-map-preserves-joins}
\Ga_{T_1 \cap T_2} \;\arr\; \Ga_{T_1} \otimes_A \Ga_{T_2}
\end{equation}
is an isomorphism.
\end{lemma}
\begin{proof}
Both the source and the target of the map (\ref{local-sm-sp-map-preserves-joins}) are coidempotent objects in $\mod_{A}$.
Therefore, it suffices to check that their supports agree, which is immediate from the first claim of
Proposition \ref{prop-hopkins-neeman}.
\end{proof}

\begin{lemma}\label{lem-sm-sp-map-preserves-finite-meets}
For a pair of specialization closed subsets $T_1, T_2 \in \spec(A)^\spc$, the square \[\xymatrix{
\Ga_{T_1 \cup T_2} && \Ga_{T_2} \ar[ll] \\
\Ga_{T_1} \ar[u] && \Ga_{T_1 \cap T_2} \ar[u]\ar[ll] 
}\] is a cofiber square.
\end{lemma}
\begin{proof}
There is a canonical map \[
\Ga_{T_1} \bigsqcup_{\Ga_{T_1 \cap T_2}} \Ga_{T_2} \;\arr\; \Ga_{T_1 \cup T_2} \;,
\] both the source and the target of which are coidempotent objects in $\mod_A$. Therefore, it suffices to check that the
supports agree, which is immediate from the first claim of Proposition \ref{prop-hopkins-neeman}.
\end{proof}

\begin{proposition}\label{prop-sm-sp-map-arbitrary-meets}
For a specialization closed subset
$T \in \spec(A)^\spc$ decomposed into a union of other specialization closed subsets---$T = \cup_{i\in I} T_i$, the natural map
\begin{equation}\label{eq1-local-prop-sm-sp-map-arbitrary-meets}
\stackrel[\{i_1, \ldots, i_k \} \subseteq I]{}{\colim} \Ga_{T_{i_1} \cup \ldots \cup T_{i_k}} \;\arr\; \Ga_T
\end{equation} where the colimit is taken over all finite subsets of $I$, is an isomorphism.
\end{proposition}
\begin{proof}
First, we note that the object $\stackrel[\{i_1, \ldots, i_k \} \subseteq I]{}{\colim} \Ga_{T_{i_1} \cup \ldots \cup T_{i_k}}$
is coidempotent because the diagram is filtered. Therefore, it suffices to prove that \[
\supp(\stackrel[\{i_1, \ldots, i_k \} \subseteq I]{}{\colim} \Ga_{T_{i_1} \cup \ldots \cup T_{i_k}}) \;=\; T \;.
\]
Fix a point $\p \in \spec(A)$. In case $\p \notin T$,
we also have $\p \notin T_{i_1} \cup \ldots \cup T_{i_k}$ for any finite subset $\{i_1, \ldots, i_k\} \subseteq I$,
and therefore, \[\begin{split}
A_\p/\p \;\otimes_A \stackrel[\{i_1, \ldots, i_k \} \subseteq I]{}{\colim} \Ga_{T_{i_1} \cup \ldots \cup T_{i_k}} &\;\simeq\; 
\stackrel[\{i_1, \ldots, i_k \} \subseteq I]{}{\colim} A_\p/\p \otimes_A \Ga_{T_{i_1} \cup \ldots \cup T_{i_k}} \;\simeq\; 0 \;.\\
\end{split}\]
On the other hand, if $\p \in T$, we pick an index $i_0 \in I$ such that $\p \in T_{i_0}$. Then we have: \[\begin{split}
A_\p / \p \;\otimes_A \stackrel[\{i_1, \ldots, i_k \} \subseteq I]{}{\colim} \Ga_{T_{i_1} \cup \ldots \cup T_{i_k}} &\;\simeq\;
A_\p / \p \;\otimes_A \stackrel[\{i_1, \ldots, i_k \} \subseteq I \setminus \{i_0\}]{}{\colim} \Ga_{T_{i_0} \cup T_{i_1} \cup \ldots \cup T_{i_k}} \;\simeq\; \\ 
&\;\simeq\; \stackrel[\{i_1, \ldots, i_k \} \subseteq I \setminus \{i_0\}]{}{\colim} A_\p/\p \otimes_A \Ga_{T_{i_0} \cup T_{i_1} \cup \ldots \cup T_{i_k}} \;\simeq\; \\
&\;\simeq\; \stackrel[\{i_1, \ldots, i_k \} \subseteq I \setminus \{i_0\}]{}{\colim} A_\p/\p \;\simeq\; A_\p / \p \;,
\end{split}\] where we first use filteredness and then
distribute the tensor product with the colimit and use the first claim of Proposition \ref{prop-hopkins-neeman}.
\end{proof}

\sssec{Completion at a specialization closed subset}\label{sssec-completion-at-a-specialization-closed}

Similar to Definition \ref{def-completion}, for a specialization closed subset $T \subseteq \spec(A)$,
we can use the coidempotent object $\Ga_T \in \cidem(\mod_A)$ to define a functor of completion at $T$.
Namely, the functor is defined by the formula \[\xymatrix{
\mod_{A_\bs} \ar[rr] && \mod_{A_\bs} \\
M \ar@{}[u]|\invertin \ar@{|->}[rr] && \uhom_{A_\bs}(\Ga_T, M) \ar@{}[u]|\invertin \ar@{}[r]|\simeq & \stackrel[Z \subseteq T]{}{\lim} M^\thicksim_Z &.
}\] It is not hard to see the resulting endofunctor is idempotent and lax symmetric monoidal.

\sssec{} We finish this section with a technical statement saying that, when computing the completion with respect to
a specialization closed subset $T \in \spec(A)^\spc$, we can take the limit only over the irreducible closed subsets contained in $T$.
\begin{proposition}\label{prop-lim-over-primes}
For a specialization closed subset $T \in \spec(A)^\spc$ and a module $M \in \mod_{A_\bs}$, the natural map \[
M^\thicksim_T \;\arr\; \stackrel[\p \in T]{}{\lim} M^\comp_\p \;,
\] where the limit is taken over the poset of all points inside $T$, is an isomorphism.
Equivalently, the map \begin{equation}\label{local-lim-over-primes-map}
\stackrel[\p \in T]{}{\colim} \Ga_{\spec(A/\p)} \;\arr\; \Ga_T
\end{equation}
is an isomorphism.
\end{proposition}
We first prove that the map (\ref{local-lim-over-primes-map}) is an isomorphism
in the case of a single closed subset $Z \subseteq \spec(A)$.
\begin{lemma}\label{lem-lim-over-primes-single-closed}
The map (\ref{local-lim-over-primes-map}) is an isomorphism when $T = Z \subseteq \spec(A)$ is a single closed subset.
In particular, the claim of Proposition \ref{prop-lim-over-primes} holds in this case.
\end{lemma}
\begin{proof}
We prove the lemma by double induction on the dimension of $Z$ and the number of irreducible components of $Z$.
\vskip 5pt
\noindent
\textit{Base}: either $\dim(Z) = 0$ or $Z$ is irreducible. In case of an irreducible closed subset, the poset of primes
has an initial point---the one that corresponds to the closed subset itself, and therefore the claim holds. In the case where
$\dim(Z) = 0$, the poset of primes is just a discrete set (with no arrows), and the claim immediately follows from
Lemma \ref{lem-sm-sp-map-preserves-finite-meets}.
\vskip 5pt
\noindent
\textit{Step}: take a closed subset $Z = Z_1 \cup Z_2$ such that the claim holds for $Z_1$, $Z_2$, and for the intersection
$Z_1 \cap Z_2$. Let $J$ denote the poset having three vertices $j_1, j_2$, and $j_{12}$ and two arrows
$j_{12} \to j_1$ and $j_{12} \to j_2$. Define a map of posets \[
H \colon \{ \p \in Z \}^\op \;\arr\; J
\] by the formula \[
H(\p) \;=\; \left\{\begin{array}{cc}
j_{12} & \text{ if } \p\in Z_{12};\\
j_1 & \text{ if } \p\in Z_1\setminus Z_{12};\\
j_2 & \text{ if } \p\in Z_2\setminus Z_{12}.
\end{array}\right.
\] We also denote by $F$ our diagram---\[\xymatrix{
F\colon \{\p \in Z\}^\op \ar[rr] && \mod_A \\
\p \ar@{}[u]|\invertin \ar@{|->}[rr] && \Ga_{\spec(A/\p)} \ar@{}[u]|\invertin &.
}\] The observation that lets us prove the induction step is that we can compute the left Kan extension of $F$ via $H$ using the 
induction hypothesis: \[\begin{split}
\mathrm{LKE}_H (F) (j_1) &\;\simeq\; \stackrel[\substack{\p \in Z \\ \; H(\p) \to j_1}]{}{\colim} \Ga_{\spec(A/\p)} \;\simeq\; \\
&\;\simeq\; \stackrel[ \p \in Z_1 ]{}{\colim} \Ga_{\spec(A/\p)} \;\simeq\; \Ga_{Z_1} ;\\
\mathrm{LKE}_H (F) (j_2) &\;\simeq\; \ldots \;\simeq\; \Ga_{Z_2} ;\\
\mathrm{LKE}_H (F) (j_{12}) &\;\simeq\; \ldots \;\simeq\; \Ga_{Z_{12}} .\\
\end{split}\]
We therefore get: \[\begin{split}
\colim \;F &\;\simeq\; \colim \;\mathrm{LKE}_H (F) \;\simeq\; \\
&\;\simeq\; \colim
\begin{pmatrix}
\Ga_{Z_1} & \\
\uparrow & \\
\Ga_{Z_{12}} & \to & \Ga_{Z_2} 
\end{pmatrix}
\;\simeq\; \\
&\;\simeq\; \Ga_Z,
\end{split}\] where the last isomorphism follows from Lemma \ref{lem-sm-sp-map-preserves-finite-meets}.
\end{proof}

\begin{proof}[Proof of Proposition \ref{prop-lim-over-primes}]
Let $H$ denote the inclusion \[
H \colon \{ \p \in T \}^\op \;\hookrightarrow\; \{ Z \subseteq T \}^\op
\] of the subposet of irreducible closed subsets of $T$ into the poset of all closed subsets of $T$. Also, let $F$ denote
our diagram---\[\xymatrix{
F\colon \{\p \in T\}^\op \ar[rr] && \mod_A \\
\p \ar@{}[u]|\invertin \ar@{|->}[rr] && \Ga_{\spec(A/\p)} \ar@{}[u]|\invertin &.
}\] Using Lemma \ref{lem-lim-over-primes-single-closed}, we compute the left Kan extension of $F$ via $H$: \[\begin{split}
\mathrm{LKE}_H(F)(Z) &\;\simeq\; \stackrel[\substack{ \p \in T \\ H(\p) \to Z}]{}{\colim} \Ga_{\spec(A/\p)} \;\simeq\; \\
&\;\simeq\; \stackrel[\substack{ \p \in Z}]{}{\colim} \Ga_{\spec(A/\p)} \;\simeq\; \Ga_Z.
\end{split}\] Finally, using this computation, we compute the colimit of $F$:
\[
\colim \; F \;\simeq\; \colim \; \mathrm{LKE}_H F \;\simeq\; \stackrel[Z \subseteq T]{}{\colim} \Ga_Z \;\simeq\; \Ga_T \;,
\] which finishes the proof.
\end{proof}

\ssec{Tubular neighborhoods and recollement}\label{ssec-tub-neig}

Let us fix a closed subset $Z \subset \spec(A)$. In the following couple of paragraphs we make a few observations about
nice solid-algebraic properties of the completion $A^\thicksim_Z$ as an algebra object of $\mod_{A_\bs}$ and describe
a solid analog of the decomposition of the category of sheaves on a topological space with respect to a closed embedding and its complement.

Recall from Lemma \ref{lem-comp-is-cmp-and-idmp} that
the completion $A^\thicksim_Z$ is an idempotent algebra inside the category
$\mod^\w_{A_\bs}$. This immediately implies the following lemma
\begin{lemma}
The forgetful functor \[
i_* \colon A^\thicksim_Z-\mod_{A_\bs} \arr \mod_{A_\bs}
\] is fully faithful and its image is a tensor ideal inside $\mod_{A_\bs}$.
Moreover, besides the left adjoint $i^*$, the functor $i_*$ admits a continuous right adjoint \[
i^! \;=\; \uhom_{A_\bs} (A^\thicksim_Z, -) \::\; \mod_{A_\bs} \arr A^\thicksim_Z-\mod_{A_\bs} \;.
\]
\end{lemma}

Let now $Z^c$ denote the complement $\spec(A) \setminus Z$.
In the following proposition we identify the factor category $\mod_{A_\bs} / \im(i_*)$ with the category
$\mod_{{\Ga(Z^c, \O)}_{\bullet}}$ of ultrasolid modules over the ring ${\Ga(Z^c, \O)}$, see
Definition \ref{def-ultrasolid}.

\begin{proposition}\label{prop-recollement-single-closed}
There is a recollement \[\xymatrix{
A^\thicksim_Z-\mod_{A_\bs} \ar@{^{(}->}[rr]|{i_*} && \mod_{A_\bs} \ar@<2ex>[ll]^{i^{\:!}}\ar@<-2ex>[ll]_{i^*} \ar[rr]|{j^*} && \mod_{A_\bs} / \im(i_*) \ar@<2ex>@{_{(}->}[ll]^{j_*}\ar@<-2ex>@{_{(}->}[ll]_{j_!} \;.
}\] The image \[
\im(j_*) \;\subseteq\; \mod_{A_\bs}
\] is compactly generated by products of copies of the ring $\Ga(Z^c, \O)$, and the functor $j_*$
satisfies the conditions of Definition \ref{defn-analytic-Eoo-ring}.
\end{proposition}

\begin{proof}
It follows formally that the functor $j_*$ satisfies the conditions of Definition \ref{defn-analytic-Eoo-ring} and
that the image $\im(j_*)$ is compactly generated by objects of the form \[
j_* j^* \prod_I A \;\simeq\; \prod_I j_* j^* A \;.
\] In particular, we see that $\im(j_*)$ is equivalent to the category of ultrasolid modules over the ring $j_* j^* A$.
It remains to see that $j_*j^* A \;\simeq\; \Ga(Z^c, \O)$:
\[\xymatrix{
i_*i^!A \ar[rr] && A \ar[rr] && j_*j^*A \\
\uhom_{A_\bs}(A_Z^\thicksim, A)\ar[rr]\ar@{}[u]|\invsimeq && A \ar@{}[u]|\inveq \\
\Gamma_Z \ar@{}[u]|\invsimeq \ar[rr] && A \ar[rr]\ar@{}[u]|\inveq && A / \Gamma_Z \ar@{}[uu]|\invsimeq \ar@{}[r]|(0.3)\simeq & \Ga(Z^c, \O) \;,
}\] where we use reflexivity of the local cohomology $\Ga_Z$ as an $A_\bs$ module.
\end{proof}

In particular, in the case where the closed subset $Z \subseteq \spec(A)$ is a hypersurface $\{f = 0\}$, the category of solid
$A$-modules decomposes as a recollement of $A_f^\comp-\mod_{A_\bs}$ and $\mod_{A[\frac{1}{f}]_\bs}$: \[\xymatrix{
A^\comp_f-\mod_{A_\bs} \ar@{^{(}->}[rr]|{i_*} && \mod_{A_\bs} \ar@<2ex>[ll]^{i^{\:!}}\ar@<-2ex>[ll]_{i^*} \ar[rr]|{j^*} && \mod_{A[\frac{1}{f}]_\bs} \ar@<2ex>@{_{(}->}[ll]^{j_*}\ar@<-2ex>@{_{(}->}[ll]_{j_!} \;.
}\]

\sssec{} More generally, we will see in Proposition \ref{prop-the-smashing-spectrum-map} that,
for any specialization closed subset $T \subseteq \spec(A)$, the category $A^\thicksim_T-\mod_{A_\bs}$ of
modules over the $T$-completed algebra $A^\thicksim_T$ in $\mod_{A_\bs}$ satisfies the same formal properties. In particular,
we will get a full continuous recollement \[\xymatrix{
A^\thicksim_T-\mod_{A_\bs} \ar@{^{(}->}[rr]|{i_*} && \mod_{A_\bs} \ar@<2ex>[ll]^{i^{\:!}}\ar@<-2ex>[ll]_{i^*} \ar[rr]|{j^*} && \mod_{A_\bs} / \im(i_*) \ar@<2ex>@{_{(}->}[ll]^{j_*}\ar@<-2ex>@{_{(}->}[ll]_{j_!} \;.
}\] We will use certain full subcategories of the form
\[
A^\thicksim_T-\mod_{A_\bs} \;\subseteq\; \mod_{A_\bs}
\]
to construct a filtration on the category $\mod_{A_\bs}$, decomposition with respect to which will give us an analog of
Beilinson-Parshin resolution for the solid modules, see section \ref{ssec-skeletal-filtration} and Theorem \ref{thm-adeles}.

\ssec{Unions of tubular neighborhoods and the smashing spectrum}\label{ssec-sm-sp}

Recall the smashing spectrum from \cite{smashing_spectrum}, Definition 3.17. We will simply denote $\sm(A_\bs)$
the smashing spectrum $\sm(\mod_{A_\bs})$
of the category $\mod_{A_\bs}$. We will denote by $\sm(A_\bs)^\cmp \subseteq \sm(A_\bs)$ the subposet spanned by algebras which are
compact as $A_\bs$-modules; we note that this subposet is closed under finite joins and finite meets.

\begin{proposition}\label{prop-the-smashing-spectrum-map}
There is a map of posets \begin{equation}\label{the-smashing-spectrum-map}\xymatrix{
(\spec(A)^\spc)^\op \ar[rr] && \sm(A_\bs)^\cmp \\
T \ar@{}[u]|\invertin \ar@{|->}[rr] && A^\thicksim_T \stackrel{\define}{=} \stackrel[Z \subseteq T]{}{\lim} A^\thicksim_Z \ar@{}[u]|\invertin \;,
}\end{equation} where the limit is taken over all closed subsets $Z \subseteq T$. Moreover, this map preserves finite meets and
finite joins.
\end{proposition}

Once we have proved that $A^\thicksim_T$ is a compact and idempotent algebra,
the fact that the map (\ref{the-smashing-spectrum-map}) preserves finite joins and finite meets follows immediately from Lemma
\ref{lem-sm-sp-map-preserves-finite-joins} and Lemma \ref{lem-sm-sp-map-preserves-finite-meets}.

First, by Lemma \ref{lem-solid-synth-cmp-obj},
the functor \[\xymatrix{
\mod^{\projfinlen}_{A} \ar[rr] && \mod^\omega_{A_\bs} \\
P \ar@{}[u]|\invertin \ar@{|->}[rr] && \uhom_{A_\bs} (P , A) \ar@{}[u]|\invertin
}\] maps coidempotent objects to idempotent algebras, and  
to prove that $A^\thicksim_T$ is a compact and idempotent algebra, it suffices to prove that
the coidempotent object $\Ga_T \in \mod_{A}$ has finite projective amplitude.

Once we are done with that, we recall from \cite{smashing_spectrum}, Theorem 3.8, how finite meets and joins are
computed in the smashing spectrum. Namely, for a pair $B_1, B_2 \in \sm(A_\bs)^\cmp$ of idempotent algebras,
the join $B_1 \vee B_2$ is computed as the tensor product $B_1 \otimes_{A_\bs} B_2$, and the meet $B_1 \wedge B_2$
is computed as the pullback (in $A_\bs$-modules) $B_1 \times_{B_1 \vee B_2} B_2$. Therefore, the fact that the map
(\ref{the-smashing-spectrum-map}) preserves finite meets and finite joins would follow from Lemma
\ref{lem-sm-sp-map-preserves-finite-joins} and Lemma \ref{lem-sm-sp-map-preserves-finite-meets}.
So to finish the proof of Proposition \ref{prop-the-smashing-spectrum-map}, it suffices to prove the following proposition.

\begin{proposition}\label{prop-comp-alg-is-idmpt-and-cmp}
The coidempotent object \[
\Ga_T \;=\; \stackrel[Z \subseteq T]{}{\colim} \Ga_Z
\] has finite projective amplitude as an $A$-module.
\end{proposition}
\begin{proof}
We will use the formula from Proposition \ref{prop-lim-over-primes}, which says that the canonical map
\[
\stackrel[\p \in T]{}{\colim} \Ga_{\spec(A/\p)} \;\arr\; \Ga_T \;,
\]
where the colimit is taken over the poset of primes $\p \in T$ ordered by opposite inclusion, is an isomorphism.
Let us denote that poset by $\I$. Because the spectrum $\spec(A)$ has finite dimension, the nerve $\mathrm{N}(\I)$
has non-degenerate simplices only in degrees $\le d$, where $d$ denotes the Krull dimension of the ring $A$.
In particular, we need at most $d$ steps to build the simplicial set $\mathrm{N}(\I) \simeq \mathrm{N}(\I)_{\le d}$ by iteratively 
adjoining non-degenerate simplicies---\[
\mathrm{N}(\I)_{\le k} \;\simeq\; \mathrm{N}(\I)_{\le k-1} \bigsqcup_{C_{k} \times \partial \Delta^k} C_{k} \times \Delta^k \;,
\] where $C_k$ denotes the set of non-degenerate $k$ simplices. Then, applying Proposition 4.4.2.2 of \cite{Lurie-HTT} iteratively,
we get a variant of Bousfield-Kan formula for our colimit: \[
\stackrel[\p \in T]{}{\colim} \Ga_{\spec(A/\p)} \;\simeq\; \Bigl( \xymatrix{
\oplus_{\p_0 \subset \ldots \subset \p_d} \Ga_{\spec(A/\p_d)} \ar[r] & \ldots \ar[r] & \oplus_{\q \subset \p} \Ga_{\spec(A/\p)} \ar[r] & \oplus_\p \Ga_{\spec(A/\p)} \\
-d & & -1 & 0\:.
} \Bigr)
\] Therefore, it suffices to prove that, for a closed subset $Z \subseteq \spec(A)$, the local cohomology $\Ga_Z$ has
bounded projective amplitude with bounds independent on $Z$. First, one easily sees that an object of the form
$A[\frac{1}{f}]/A[-1]$ can be represented as a two term complex consisting of free $A$-modules; in particular,
it has projective amplitude in $[0,1]$. Therefore, the local cohomology $\Ga_Z$ can be represented by a complex free
$A$-modules in degrees from $0$ to $l$, where $l$ is the number of generators of the ideal of $Z$. On the other hand,
by \cite[\href{https://stacks.math.columbia.edu/tag/0DXC}{Lemma 0DXC}]{stacks-project}, the object $\Ga_Z$ has
non-zero cohomology only in degrees from $0$ to $d$; therefore, in the presentation of $\Ga_Z$ as
complex of free $A$-modules in degrees from $0$ to $l$, we can split off from the right until the complex has
non-zero terms only in degrees from $0$ to $d$.
\end{proof}

\begin{remark}
Let us note that the proof of Proposition \ref{prop-comp-alg-is-idmpt-and-cmp} works for any Noetherian ring of finite
dimension, which is more general than we need in this paper, but could be useful in future works.
In the finite type case, there is shorter argument that uses finiteness of the global dimension of a polynomial ring.
\end{remark}

\begin{corollary}\label{cor-specialization-closed-subset-recollement}
Similar to Proposition \ref{prop-recollement-single-closed}, for a specialization closed subset $T \subseteq \spec(A)$
we get a full continuous recollement \[\xymatrix{
A^\thicksim_T-\mod_{A_\bs} \ar@{^{(}->}[rr]|{i_*} && \mod_{A_\bs} \ar@<2ex>[ll]^{i^{\:!}}\ar@<-2ex>[ll]_{i^*} \ar[rr]|{j^*} && \mod_{\Ga(T^c, \O)_\bullet} \ar@<2ex>@{_{(}->}[ll]^{j_*}\ar@<-2ex>@{_{(}->}[ll]_{j_!} \;,
}\] where \[
\Ga(T^c, \O) \;\simeq\; \stackrel[Z \subseteq T]{}{\colim} \Ga(Z^c, \O) 
\] is the idempotent algebra supported on the generalization closed subset $T^c = \spec(A) \setminus T$.
\end{corollary}

\ssec{Unions of tubular neighborhoods, completions, and ultrasolid analytic ring structures}\label{ssec-gen-closed-and-ultrasolid}

Recall that a subset $E \subseteq \spec(A)$ is called generalization closed if, whenever $E$ contains a point $\p$,
$E$ also contains any point $\q$ specializing to $\p$. We denote by $\spec(A)^\gen$
the poset of generalization closed subsets of $\spec(A)$.
For a generalization closed subset $E \in \spec(A)^\gen$, we introduce the functor of completion at $E$:
\[\xymatrix{
( - )^\thicksim_E : & \mod_{A_\bs} \ar[rr] && \mod_{A_\bs} \\
& M \ar@{}[u]|\invertin \ar@{|->}[rr] && M^\thicksim_{E} \;=\; \uhom_{A_\bs}(A^\thicksim_{E^c} / A [-1], M ) \ar@{}[u]|\invertin &,
}\]
where $E^c = \spec(A) \setminus E$ is the complement, which is a specialization closed subset of $\spec(A)$.
By the results of the previous section, this functor is continuous and idempotent. Moreover, we have a recollement
\[\xymatrix{
A^\thicksim_{E^c}-\mod_{A_\bs} \ar@{^{(}->}[rr]|{} && \mod_{A_\bs} \ar@<2ex>[ll]^{}\ar@<-2ex>[ll]_{} \ar[rr]|(0.35){j^*} && \bigl(\mod_{A_\bs}\bigr)/ \bigl(A^\thicksim_{E^c}-\mod_{A_\bs} \bigr) \ar@<2ex>@{_{(}->}[ll]^(0.65){j_*}\ar@<-2ex>@{_{(}->}[ll]_(0.65){j_!} \;,
}\] where all of the functors are continuous, and, for the completion at $E$, we have $(-)^\thicksim_E \;\simeq\; j_*j^*$.

As our first basic example, we note that the
calculations from section \ref{ssec-tub-neig} imply that, for an element $f \in A$,
the completion at the generalization closed subset $\{ f = 0 \}^c$ is given by $A[\frac{1}{f}]_\bs \otimes_{A_\bs} -$.

\sssec{}
For a prime ideal $\p \subset A$, we define a generalization closed subset $\{\p\}^\gen \in \spec(A)^\gen$
as the intersection of all Zariski open subsets containing the point $\p \in \spec(A)$: \[
\{\p\}^\gen \;\stackrel{\define}{=}\; \bigcap_{\substack{U \subseteq \spec(A) \\ \p \in U}} U \;.
\] We claim that the completion at $\{\p\}^\gen$ is given by ${A_\p}_\bullet \otimes_{A_\bs} -$, where ${A_\p}_\bullet$
denotes the ultrasolid analytic ring structure on the local ring $A_\p$, see Definition \ref{def-ultrasolid}.
\begin{lemma}\label{lem-ultrasolid-on-localized-at-prime}
There exists an ultrasolid analytic ring structure ${A_\p}_\bullet$ on the ring $A_\p$.
Moreover, there is a map of analytic rings $A_\bs \arr {A_\p}_\bullet$, and we have a natural isomorphism \[
(-)^\thicksim_{\{\p\}^\gen} \;\simeq\; {A_\p}_\bullet \otimes_{A_\bs} - \;.
\]
\end{lemma}
\begin{proof}
Indeed, just from looking at the recollement
\[\xymatrix{
A^\thicksim_{(\{\p\}^\gen)^c}-\mod_{A_\bs} \ar@{^{(}->}[rr]|{} && \mod_{A_\bs} \ar@<2ex>[ll]^{}\ar@<-2ex>[ll]_{} \ar[rr]|(0.35){j^*} && \bigl(\mod_{A_\bs}\bigr)/ \bigl(A^\thicksim_{(\{\p\}^\gen)^c}-\mod_{A_\bs} \bigr) \ar@<2ex>@{_{(}->}[ll]^(0.65){j_*}\ar@<-2ex>@{_{(}->}[ll]_(0.65){j_!}
}\] we get:
\begin{enumerate}
\item the image of $j_*$ is closed under both limits and colimits;
\item the factor category is compactly generated by images of the
compact generators under $j^*$, which are products of copies of $j^* A$ because $j^*$ preserves limits;
\item provided that $j_*j^*A \in \mod_{A_\bs}^{\le 0}$, the factor category admits a $t$-structure such that the functor $j_*$ is t-exact.
\end{enumerate}
It remains to compute $j_*j^*A$. Using Proposition \ref{prop-sm-sp-map-arbitrary-meets} and the equality \[
(\{\p\}^\gen)^c \;=\; \bigcup_{\substack{Z \subseteq \spec(A) \\ \p \notin Z}} Z \;=\; \bigcup_{f \in A\setminus\p} \{f = 0\} \;,
\] we get an isomorphism \[
\Ga_{(\{\p\}^\gen)^c} \;\simeq\; \stackrel[f \in A\setminus\p]{}{\colim} A[\frac{1}{f}]/A[-1] \;.
\]
Now since the object $A^\thicksim_{(\{p\}^\gen)^c} \in \mod_{A_\bs}$ is compact by Proposition \ref{prop-the-smashing-spectrum-map}
and any compact $A_\bs$ module is reflexive, we get: \[\begin{split}
j_*j^* A &\;\simeq\; \uhom_{A_\bs} (A^\thicksim_{(\{\p\}^\gen)^c} / A [-1], A) \;\simeq\; \\
&\;\simeq\; \fib( \Ga_{(\{\p\}^\gen)^c} \;\to\; A ) \;\simeq\; \\
&\;\simeq\; \stackrel[f \in A\setminus\p]{}{\colim} A[\frac{1}{f}] \;\simeq\; A_\p \;.
\end{split}\]
\end{proof}

\sssec{} We also remark that, by an identical argument, for an arbitrary generalization closed susbet $E \in \spec(A)^\gen$,
we get a natural isomorphism
\[
(-)^\thicksim_E \;\simeq\; {\Ga(E, \O)}_\bullet \otimes_{A_\bs} - \;,
\] where $\Ga(E, \O) \;\simeq\;\; \stackrel[U \supseteq E]{}{\colim} \Ga(U, \O)$ is the complex of functions on $E$, which may not
be connective.

\ssec{Noncommutative stratifications and the smashing spectrum}\label{ssec-stratifications-and-sm-sp}

In this section, we recall the language of stratification from \cite{stratified}, which provides a nice general framework
for studying resolutions in a category stratified over a poset. We will only the results of this section
in a partial case of a stratification over a finite linearly ordered set, which admits an alternative more elementary
approach, recorded in section \ref{ssec-cubical-resolutions}.
Let us fix a presentable stable category $\CC$. Recall from \cite{stratified}, Definitions 1.3.1 and 2.3.1, that
a full presentable stable subcategory $\ZZ \;\hookrightarrow\; \CC$ is called a closed subcategory the embedding admits a right 
adjoint which itself is required to admit a right adjoint. In loc. cit., the last right adjoint is not required to be continuous
and it is often not, but in most applications discussed in this paper it will be continuous.
\begin{remark}
Fix a finite type ring $A$ and a generalization closed subset $E \in \spec(A)^\gen$.
In section \ref{ssec-gen-closed-and-ultrasolid} we studied the category ${\Ga(E, \O)}_\bullet - \mod$
of ultrasolid modules over the ring $\Ga(E, \O)$. In particular, from that discussion it follows that
the category ${\Ga(E, \O)}_\bullet - \mod$ embedded via the lower-! embedding into $\mod_{A_\bs}$ is an example
of a closed subcatory inside $\mod_{A_\bs}$, and that is a common example from the point of view of this paper.
In particular, it would make more sense to call these subcategories open because that is what they are from
the six functors point of view. Nonetheless, we use the terminology of \cite{stratified} in this section.
\end{remark}
Following \cite{stratified}, we denote by $\cls(\CC)$ the poset of closed subcategories
of $\CC$. Given a poset $P$ together with a map \[\xymatrix{
P \ar[rr] && \cls(\CC) \\
p \ar@{}[u]|\invertin \ar@{|->}[rr] && \CC_p \subseteq \CC \ar@{}[u]|\invertin &,
}\] we denote by $j_p$ the inclusion $\CC_p \hookrightarrow \CC$ and by $j_p^R$ and $j_p^{RR}$ its right adjoint and the further right
adjoint respectively: \[\xymatrix{
\CC_{p} \ar@<2ex>@{^{(}->}[rr]^{j_{p}}\ar@<-2ex>@{^{(}->}[rr]_{j_{p}^{RR}}
&& \CC \ar[ll]|{j_{p}^R} \;.
}\] We denote by $\CC_{< p}$ the closed subcategory generated under colimits by all $\CC_q$ for $q < p$---that does indeed define
a closed subcategory, see Observation 2.3.9 in \cite{stratified}. Similarly, we denote by $\CC_{\gr, p}$ the quotient $\CC_p / \CC_{< p}$;
we get a recollement \begin{equation*}\xymatrix{
\CC_{< p} \ar@<2ex>@{^{(}->}[rr]^{j_{< p \to p}}\ar@<-2ex>@{^{(}->}[rr]_{j_{< p \to p}^{RR}}
&& \CC_p \ar[ll]|{j_{< p \to p}^R} \ar@<2ex>[rr]^{i_{p\to \gr, p}}\ar@<-2ex>[rr]_{i^{RR}_{p\to \gr, p}} && \CC_{\gr, p} \ar@{_(->}[ll]|{i^R_{p\to \gr, p}}\;.
}\end{equation*} We denote by $L_p$ the endofunctor of $\CC$ given by localizing to $\CC_{\gr, p}$ and embedding back
via the right adjoint; concretely, $L_p$ is given by the following composite: \[\xymatrix{
\CC \ar[rr]^{j_p^R} && \CC_p \ar[rr]^{i_{p\to \gr, p}} && \CC_{\gr, p} \ar@{^(->}[rr]^{i^R_{p\to \gr, p}} && \CC_p \ar@{^(->}[rr]^{j_p^{RR}} && \CC \;\;.
}\]

\sssec{} For a symmetric monoidal stable presentable category $\CC$, 
there is a map of posets \begin{equation}\label{smashing-spectrum-to-closed-subcat}\xymatrix{
\sm(\CC) \ar[rr] && \cls_\CC \\
B \ar@{}[u]|\invertin \ar@{|->}[rr] && \CC_B \ar@{}[u]|\invertin &,
}\end{equation}
where we define $\CC_B$ as the closed subcategory given by the quotient $\CC / (B-\mod_\CC)$ embedded into $\CC$ via the left adjoint
to the quotient functor.

\begin{proposition}\label{prop-stratify-over-smashing-spectrum}
The map (\ref{smashing-spectrum-to-closed-subcat}) satisfies the stratification condition of Definitions 2.4.1 and 2.4.3, \cite{stratified}.
\end{proposition}
\begin{proof}
According to loc. cit., we need to prove that, for a pair of idempotent algebras $B_1, B_2 \in \sm(\CC)$, there is dotted arrow
(necessarily unique up to a contractible space of choices) making the diagram
\begin{equation}\label{local-stratify-over-sm-sp-diag}\xymatrix{
\CC_{< B_1 \cap < B_2} \ar@{^(->}[rr]^(0.6){j_{< B_1 \cap < B_2 \to B_1}} && \CC_{B_1} \\
\CC_{B_2} \ar@{-->}[u] \ar@{^(->}[rr]^{j_{B_2}} && \CC \ar[u]^{j^R_{B_1}} &,
}\end{equation}
where $\CC_{< B_1 \cap < B_2}$ denotes the closed subcategory generated under colimits by all $\CC_B$ for $B \in \sm(\CC)$ mapping
to both $B_1$ and $B_2$, commutative. In view of the fact that the poset $\sm(\CC)$ has finite meets,
see Theorem 3.8 in \cite{smashing_spectrum}, we have an equivalence of closed subcategories \[
\CC_{< B_1 \cap < B_2} \;\simeq\; \CC_{B_1 \wedge B_2} \;,
\] where the idempotent algebra $B_1 \wedge B_2$ is described as the pull-back \[\xymatrix{
B_1 \wedge B_2 \ar@{}[rrd]|(0.25)\pb \ar[rr]\ar[d] && B_1 \ar[d] \\
B_2 \ar[rr] && B_1 \otimes_{\CC} B_2 &.
}\]

We now claim that the functor $j^R_{B_1\wedge B_2 \to B_2} \colon \CC_{B_2} \;\arr\; \CC_{B_1 \wedge B_2}$ makes the diagram
(\ref{local-stratify-over-sm-sp-diag}) commutative (via the Beck-Chevalley natural transformation). We first observe that each of the
involved functors is linear over $\CC$; therefore, it suffices to prove that the Beck-Chevalley map
\begin{equation}\label{local-stratify-sm-sp-BC-map}
j_{B_1 \wedge B_2 \to B_1} \: j^R_{B_1\wedge B_2 \to B_2} \;\arr\; j^R_{B_1} \: j_{B_2}
\end{equation} becomes an isomorphism when evaluated at the unit $1_{\CC_{B_2}} \;\simeq\; j^R_{B_2} (1_\CC)$. Moreover,
we can also compose the map (\ref{local-stratify-sm-sp-BC-map}) with $j^R_{B_1}$ since the latter is conservative. Now, using a simple
observation that, for an algebra $B \in \sm(\CC)$, we have \[
j_B \: j_B^R \;\simeq\; B/1_\CC [-1] \otimes_\CC \text{---} \;,
\] we finish the computation: \[\begin{split}
j_{B_1} \: j_{B_1 \wedge B_2 \to B_1} \: j^R_{B_1\wedge B_2 \to B_2} \: j^R_{B_2} (1_\CC) &\;\simeq\; \\
&\;\simeq\; j_{B_1 \wedge B_2} \: j^R_{B_1\wedge B_2} (1_\CC) \;\simeq\; \\
&\;\simeq\; (B_1 \wedge B_2) / 1_\CC [-1] \;\simeq\; \\
&\;\simeq\; (B_1 \times_{B_1 \otimes_\CC B_2} B_2) / 1_\CC [-1] \;,
\end{split}\] and \[\begin{split}
j_{B_1} \: j^R_{B_1} \: j_{B_2} \: j^R_{B_2} (1_\CC) &\;\simeq\; \\
&\;\simeq\; B_1 / 1_\CC \otimes_\CC B_2 / 1_\CC [-2] \;.
\end{split}\]
\end{proof}

\sssec{} Given a poset $P$, one can form the subdivision category $\sd(P)$, see Definition A.3.2 in \cite{stratified}; its objects are
given by injective functors $[n] \arr P$, its morphisms are given by order preserving maps $[m] \arr [n]$ commuting with the maps to $P$.
Proposition \ref{prop-stratify-over-smashing-spectrum} together with Theorem 2.5.14 of \cite{stratified} imply that, for a poset $P$
mapping to $\sm(\CC)$, we get a diagram \[\xymatrix{
\sd(P) \ar[rr] && \mathrm{End}_{\pr^L_\st}(\CC) \\
[ p_0 < \ldots < p_n ] \ar@{}[u]|\invertin \ar@{|->}[rr] && L_{p_n} \ldots L_{p_0} \ar@{}[u]|\invertin
}\] and a map \begin{equation}\label{resolution-map-sm-sp}
\id \;\arr\; \stackrel[{[} p_0 < \ldots < p_n {]} \in \sd(P)]{}{\lim} L_{p_n} \ldots L_{p_0} \;.
\end{equation} The following statement follows immediately from Theorem 2.5.14, \cite{stratified}.
\begin{proposition}\label{resolution-stratified-over-poset}
The map (\ref{resolution-map-sm-sp}) is an isomorphism in case the poset $P$ is down-finite
(as defined by Definition 2.1.6, \cite{stratified}) and the closed subcategories $\CC_p$ generate the whole $\CC$ under colimits.
\end{proposition}

\ssec{Cubical resolutions}\label{ssec-cubical-resolutions}

We will only need the discussion from the previous sections in the case of stratifications of the form \[
[0, d]^\op \;\arr\; \cls(\CC) \;,
\] where $[0, d]$ denotes the poset $0 < 1 < \ldots < d$. (That $\op$ does not bear any actual meaning, but is there rather to make
sure that certain notation in this section agrees with notation that will be used later.) In this case, the subdivision category
\[
\sd([0,d]) \;=\; \{ 0 \le k_0 < k_1 < \ldots < k_m \le d \:, \; m \ge 0\}
\]
is equivalent to a $d$-dimensinal cube without one of its vertices---without the one that would correspond to the empty subdivision.
According to Proposition \ref{resolution-stratified-over-poset}, we get a cubical resolution: \begin{equation}\label{local-cubical-res-eq1}
\id \;\stackrel{\simeq}{\arr}\; \stackrel[\substack{[0 \le k_0 < k_1 < \ldots < k_m \le d] \\ m \ge 0}]{}{\lim} \; L_{k_0} \ldots L_{k_m} \;.
\end{equation}
In this section, we would like to give a simple, independent of \cite{stratified},
construction of the cubical diagram plus a simple independent argument that the map (\ref{local-cubical-res-eq1})
is an isomorphism. This was essentially borrowed from section 2 of \cite{antolíncamarena2014chromaticfracturecubes}.
\begin{proposition}\label{resolution-filtration}
There is a diagram \[\xymatrix{
\{ [0 \le k_0 < k_1 < \ldots < k_m \le d] , \; m \ge -1\} \ar[rr] && \mathrm{End}_{\pr^L_\st}(\CC) \\
[ 0 \le k_0 < \ldots < k_m \le d] \ar@{}[u]|\invertin \ar@{|->}[rr] && L_{k_0} \ldots L_{k_m} \ar@{}[u]|\invertin &,
}\] where the empty tuple maps to the identity endofunctor. The induced map \[
\id \;\stackrel{\simeq}{\arr}\; \stackrel[\substack{[0 \le k_0 < k_1 < \ldots < k_m \le d] \\ m \ge 0}]{}{\lim} \; L_{k_0} \ldots L_{k_m}
\] is an isomorphism.
\end{proposition}
\begin{proof}
The proof is done by induction on $d$. The base case $d = 0$ is trivial. Let's assume that $d \ge 1$. We have a recollement
\begin{equation}\label{local-direct-cube-construction-eq0}\xymatrix{
\CC_{1} \ar@<2ex>@{^{(}->}[rr]^{j_{1}}\ar@<-2ex>@{^{(}->}[rr]_{j_{1}^{RR}}
&& \CC \ar[ll]|{j_{1}^R} \ar@<2ex>[rr]^{i_{0 \to \gr, 0}}\ar@<-2ex>[rr]_{i^{RR}_{0 \to \gr, 0}} && \CC_{\gr, 0} \ar@{_(->}[ll]|{i^R_{0\to \gr, 0}}\;.
}\end{equation}
We apply the induction hypothesis to $\CC_1$ and get a diagram
\begin{equation}\label{local-direct-cube-construction-eq1}\xymatrix{
\{ [1 \le k_0 < k_1 < \ldots < k_m \le d] , \; m \ge -1\} \ar[rr] && \mathrm{End}_{\pr^L_\st}(\CC_1) \\
[ 1 \le k_0 < \ldots < k_m \le d] \ar@{}[u]|\invertin \ar@{|->}[rr] && L_{\CC_1, k_0} \ldots L_{\CC_1, k_m} \ar@{}[u]|\invertin &,
}\end{equation}
for which the claim holds. By construction, we have isomorphisms \begin{equation}\label{local-direct-cube-construction-eq1.5}
L_k \;\simeq\; j_1^{RR} \: L_{\CC_1, k} \: j_1^R \;,
\end{equation}
for $k \ge 1$, and \[
L_0 \;\simeq\; i_{0 \to \gr, 0}^R \: i_{0 \to \gr, 0} \;.
\] Now, in order to get a cubical diagram of the form
\begin{equation}\label{local-direct-cube-construction-eq2}\xymatrix{
\{ [0 \le k_0 < k_1 < \ldots < k_m \le d] , \; m \ge -1\} \ar[rr] && \mathrm{End}_{\pr^L_\st}(\CC) \\
[ 0 \le k_0 < \ldots < k_m \le d] \ar@{}[u]|\invertin \ar@{|->}[rr] && L_{k_0} \ldots L_{k_m} \ar@{}[u]|\invertin &,
}\end{equation} we first compose the diagram (\ref{local-direct-cube-construction-eq1}) with $j_1^{RR}$ on the left and with
$j_1^R$ on the right, and then also compose with the unit map $\id \arr L_0$ on the left.

Finally, by construction, in order to prove that the map \[
\id \;\stackrel{}{\arr}\; \stackrel[\substack{[0 \le k_0 < k_1 < \ldots < k_m \le d] \\ m \ge 0}]{}{\lim} \; L_{k_0} \ldots L_{k_m}
\] induced by (\ref{local-direct-cube-construction-eq2}) is an isomorphism, it suffices to prove that the map \[
\fib\bigl(\id \to L_0\bigr) \;\arr\; \stackrel[\substack{[1 \le k_0 < k_1 < \ldots < k_m \le d] \\ m \ge 0}]{}{\lim} \; \fib \bigl(\id \to L_0\bigr) L_{k_0} \ldots L_{k_m}
\] is an isomorphism. Using the recollement (\ref{local-direct-cube-construction-eq0}), we compute the fiber: \[
\fib\bigl(\id \to L_0\bigr) \;\simeq\; j_1  j_1^R \;.
\] We also use that, the expression (\ref{local-direct-cube-construction-eq1.5}), and fully faithfulness of $j^{RR}_1$,
to rewrite each term appearing in the limit: \[\begin{split}
\fib \bigl(\id \to L_0\bigr) L_{k_0} \ldots L_{k_m} &\;\simeq\; j_1  j_1^R L_{k_0} \ldots L_{k_m} \;\simeq\; \\
&\;\simeq\; j_1  j_1^R j^{RR}_1 L_{\CC_1, k_0} \ldots L_{\CC_1, k_m} j_1^R \;\simeq\;\\
&\;\simeq\; j_1 L_{\CC_1, k_0} \ldots L_{\CC_1, k_m} j_1^R \;.
\end{split}\] We are now left to prove that the map \[
j_1  j_1^R \;\arr\; \stackrel[\substack{[1 \le k_0 < k_1 < \ldots < k_m \le d] \\ m \ge 0}]{}{\lim} \; j_1 L_{\CC_1, k_0} \ldots L_{\CC_1, k_m} j_1^R 
\] is an isomorphism, which immediately follows from the induction hypothesis.
\end{proof}

%% file: skeleton.tex
\section{The skeletal filtration}

In this section, we introduce the skeletal filtration of the analytic spectrum of the analytic ring $A_\bs$,
where, as usual, $A$ denotes a commutative ring of finite type over $\Z$. 
Heuristically, the  skeletal filtration at level $k$ is given by the union of the tubular neighborhoods of closed
subschemes of dimension at most $k$. After defining the filtration, we analyze the strata and prove the main theorem in the affine case.

\ssec{The skeletal filtration}\label{ssec-skeletal-filtration}

\newcommand{\sk}{\mathrm{S}}

As usual, $A$ denotes a ring of finite type over $\mathbb{Z}$. We denote by $d$ its Krull dimension.
Let $[0,d]$ denote the poset $0 < 1 < 2 < \ldots < d$; we define a map of posets
\begin{equation}\label{eq-skeleton-as-spc-closed-subsets}\xymatrix{
[0,d] \ar[rr] && \spec(A)^\spc \\
k \ar@{}[u]|\invertin \ar@{|->}[rr] && \sk_{k-1} \stackrel{\define}{=} \stackrel[\substack{Z \subseteq \spec(A) \\ \dim(Z) \le k-1}]{}{\cup} Z \ar@{}[u]|\invertin &.
}\end{equation} The \textit{skeletal filtration} \[\xymatrix{
[0,d]^\op \ar[rr] && \cls(\mod_{A_\bs}) \\
k \ar@{}[u]|\invertin \ar[rr] && \bigl(\mod_{A_\bs}\bigr)_k \ar@{}[u]|\invertin
}\] on $\mod_{A_\bs}$ is defined by composing the map (\ref{eq-skeleton-as-spc-closed-subsets})
with maps (\ref{the-smashing-spectrum-map}) and (\ref{smashing-spectrum-to-closed-subcat}).
Note that, by construction, the closed subcategory $\bigl(\mod_{A_\bs} \bigr)_{k}$ fits into the recollement
\[\xymatrix{
A^\thicksim_{\sk_{k-1}}-\mod_{A_\bs} \ar@{^{(}->}[rr]|{} && \mod_{A_\bs} \ar@<2ex>[ll]^{}\ar@<-2ex>[ll]_{} \ar[rr]|{} && \bigl(\mod_{A_\bs}\bigr)_k \ar@<2ex>@{_{(}->}[ll]^{}\ar@<-2ex>@{_{(}->}[ll]_{} \;,
}\] where \[
A^\thicksim_{\sk_{k-1}} \;\simeq\; \stackrel[\substack{Z \subseteq \spec(A)\\\dim(Z) \le k-1}]{}{\lim} A^\thicksim_Z \;\;.
\] is the compact idempotent commutative algebra in $\mod_{A_\bs}$ given by $A$ completed at the specialization closed subset
$S_{k-1} \subseteq \spec(A)$, see section \ref{sssec-completion-at-a-specialization-closed} and
Proposition \ref{prop-the-smashing-spectrum-map}. Also see Corollary \ref{cor-specialization-closed-subset-recollement}
and section \ref{ssec-gen-closed-and-ultrasolid} for the discussion about the recollement of $\mod_{A_\bs}$
associated with a specialization closed subset.

\ssec{Strata of the skeletal filtration}

The goal of this subsection is to compute the graded pieces $\bigl(\mod_{A_\bs}\bigr)_{\gr, k}$ as full subcategories inside
$\mod_{A_\bs}$ together with the localization functors $\mod_{A_\bs} \arr \bigl(\mod_{A_\bs}\bigr)_{\gr, k}$. We first note
that, by construction, the category $\bigl(\mod_{A_\bs}\bigr)_{\gr, k}$ fits into the recollement
\begin{equation}\label{eq-rec-for-z}\xymatrix{
A^\thicksim_{S_{k-1}}-\mod_{A_\bs} \ar@{^{(}->}[rr]|{{i_{k-1 k}}_*} && A^\thicksim_{S_k}-\mod_{A_\bs} \ar@<2ex>[ll]^{i^{\:!}_{k-1 k}}\ar@<-2ex>[ll]_{i^*_{k-1 k}} \ar[rr]|{t_k^*} && \bigl(\mod_{A_\bs}\bigr)_{\gr, k} \ar@<2ex>@{_{(}->}[ll]^{{t_k}_*}\ar@<-2ex>@{_{(}->}[ll]_{{t_k}_!} \;,
}\end{equation} and the task amounts to computing the right part of the recollement.

\begin{lemma}\label{lem-single-irr-comp-mod-div-skeleton}
Let $Y \subset \spec(A)$ be an irreducible closed subset of dimension not greater than $k$.
Let $\p \subset A$ denote the ideal of $Y$. We claim that we have a natural isomorphism
\begin{equation*}
{t_k}_* t_k^* A^\thicksim_Y \;\simeq\; \left\{ \begin{array}{cc}
\leftcomp A_\p, & \text{ if } \dim(Y) = k ; \\
0, & \text{ otherwise.}
\end{array}\right.
\end{equation*}
\end{lemma}

\begin{proof}
First, we note that, if $\dim(Y) < k$, then $A^\thicksim_Y \in A^\thicksim_{S_{k-1}} - \mod_{A_\bs}$, and we immediately have
$t^*_k A^\thicksim_Y = 0$.
In the case where $\dim(Y) = k$, it follows from the recollemt that we have an isomorphism \[
\uhom_{A_\bs} (A^\thicksim_{S_{k-1}} / A[-1], A^\thicksim_Y) \;\simeq\; {t_k}_*t_k^* A^\thicksim_Y \:.
\] By Lemma \ref{lem-sm-sp-map-preserves-finite-joins}, we have an isomorphism \[
A^\thicksim_{S_{k-1}} \otimes_{A_\bs} A^\thicksim_Y \;\simeq\; A^\thicksim_{Y\cap S_{k-1}} \;;
\] therefore, using the adjunction, we get an isomorphism \[
\uhom_{A_\bs} (A^\thicksim_{Y\cap S_{k-1}} / A^\thicksim_Y[-1], A^\thicksim_Y) \;\simeq\; {t_k}_*t_k^* A^\thicksim_Y \;.
\]
We now recall from Lemma \ref{lem-ultrasolid-on-localized-at-prime} that \[
\leftcomp A_\p \;\simeq\; {A_\p}_\bullet \otimes_{A_\bs} A^\thicksim_Y \;\simeq\; \uhom_{A_\bs} (A^\thicksim_{Y\cap (\{\p\}^\gen)^c} / A^\thicksim_Y[-1], A^\thicksim_Y) \;,
\] where $(\{\p\}^\gen)^c$ denotes the complement to the minimal generalization closed subset of $\spec(A)$ containing $\p$.
So it suffices to prove that we have an equality
\[
S_{k-1} \cap Y \;=\; (\{\p\}^\gen)^c \cap Y \;.
\] We obviously have an inclusion \[
(\{\p\}^\gen)^c \cap Y \;=\; \bigcup_{f \in A \setminus \p} \{f=0\} \cap Y \;\subseteq\; S_{k-1} \cap Y \;,
\] so it suffices to prove that, for each closed subset $Z \subseteq Y$, $\dim(Z) < k$, there exists an element $f \in A\setminus\p$ such
that $Z \subseteq \{f=0\}$. Any element from $I_Z \setminus \p \ne \emptyset$ satisfies this property.
\end{proof}

\begin{corollary}\label{cor-single-comp-mod-div-skeleton}
Let $T \in \spec(A)^\spc$ be a specialization closed subset such that, for any closed subset $Z \subseteq T$,
we necessarily have $\dim(Z) \le k$. We claim that \begin{equation*}
{t_k}_* t_k^* A^\thicksim_T \;\simeq\; \prod_{\substack{\p \in T \\ \dim(A/\p) = k}}\leftcomp A_{\p} \;.
\end{equation*}
\end{corollary}

\begin{proof}
Using Proposition \ref{prop-lim-over-primes}, we write \[
A^\thicksim_T \;\simeq\; \stackrel[\p \in T]{}{\lim} A^\comp_\p \;.
\] Now we apply ${t_k}_* t^*_k$ to the limit and use Lemma \ref{lem-single-irr-comp-mod-div-skeleton}.
\end{proof}

\begin{proposition}\label{prop-strata-as-analytic-rings}
The category $\bigl(\mod_{A_\bs}\bigr)_{\gr, k}$ is equivalent to the category \[
\Biggl( \prod_{\substack{\p \in \spec(A) \\ \dim(A/\p) = k}} \leftcomp A_\p \Biggr)_\bullet - \mod
\] of ultrasolid modules over the ring
$\prod_{\substack{\p \in \spec(A) \\ \dim(A/\p) = k}} \leftcomp A_\p$. The adjunction \[\xymatrix{
A^\thicksim_{S_k}-\mod_{A_\bs} \ar[rr]|{t_k^*} && \bigl(\mod_{A_\bs}\bigr)_{\gr, k} \ar@<2ex>@{_{(}->}[ll]^{{t_k}_*}
}\] corresponds to the map of analytic rings \[
\Bigl( A^\thicksim_{S_k}, A_\bs \Bigr) \;\arr\; \Biggl( \prod_{\substack{\p \in \spec(A) \\ \dim(A/\p) = k}} \leftcomp A_\p \Biggr)_\bullet
\] inducing the natural map on the underlying rings.
\end{proposition}

\begin{proof}
Applying Corollary \ref{cor-single-comp-mod-div-skeleton} to $T = S_k$, we get \[
{t_k}_* t_k^* A^\thicksim_{S_k} \;\simeq\; \prod_{\p, \;\dim(A/\p) = k} \leftcomp A_\p \;.
\] The rest follows formally from the recollement.
\end{proof}

\begin{proposition}\label{prop-ultrasolid-stratum-t-exactness}
In addition to the previous proposition, the ring \[
\prod_{\substack{\p \in \spec(A) \\ \dim(A/\p) = k}} \leftcomp A_\p
\] is static, i.e. as an $A_\bs$ module it sits in the heart of the t-structure,
and the functor \begin{equation}\label{prop-ultrasolid-stratum-t-exactness-loc-functor}
\Bigl( \prod_{\substack{\p \in \spec(A) \\ \dim(A/\p) = k}} \leftcomp A_\p \Bigr)_\bullet \otimes_{A_\bs} - \;\colon\;\; \mod_{A_\bs} \;\arr\; \mod_{A_\bs} \;,
\end{equation}
given by localizing and forgetting back, is right t-exact. In particular, this analytic ring structure falls under the more restrictive
Definition \ref{def-analytic-ring}.
\end{proposition}

\begin{proof}
We first prove that the ring in question is static. By the t-exactness of products,
it suffices to prove that, for a finite type commutative ring $A$ and
a prime ideal $\q \in A$, the completed local ring $\leftcomp A_\q$ sits in the heart of the t-structure---\[
\leftcomp A_\q \in \mod^\heartsuit_{A_\bs} \;.
\] By Lemma \cite[\href{https://stacks.math.columbia.edu/tag/0921}{Tag 0921}]{stacks-project}, we can express the completed local
ring as follows:
\[
\leftcomp A_\q \;\simeq\; \stackrel[N \ge 0]{}{\lim} A_\q / \q_N \;,
\] where $\q_0 = \q = (f_1, \ldots, f_r)$ and $\q_N = (f_1^N, \ldots, f_r^N)$. It remains to observe that there is no $R^1\lim$
term because the Mittag-Leffler condition is satisfied.

We now prove that the functor (\ref{prop-ultrasolid-stratum-t-exactness-loc-functor}) is right t-exact.
Because the connective part $\mod^{\le 0}_{A_\bs}$ is compactly generated by products of copies of $A$ and filtered colimits are
t-exact, it suffices to prove that the object \[
\Bigl( \prod_{\substack{\p \in \spec(A) \\ \dim(A/\p) = k}} \leftcomp A_\p \Bigr)_\bullet \otimes_{A_\bs} \prod_J A \;\simeq\; \prod_J \prod_{\substack{\p \in \spec(A) \\ \dim(A/\p) = k}} \leftcomp A_\p
\] is connective, which it is since the ring $\prod_{\substack{\p \in \spec(A) \\ \dim(A/\p) = k}} \leftcomp A_\p$ is static.
\end{proof}

%% file: main_thm.tex
\ssec{Main theorem: the Beilinson-Parshin fracture cube}
Recall that, by Proposition \ref{prop-stratify-over-smashing-spectrum} and Proposition \ref{resolution-stratified-over-poset},
there is a limit cube of endofunctors of $\CC$ associated with a map \[
[0,d] \;\arr\; \sm(\CC) \;,
\] mapping $0$ to the trivial algebra (whose underlying object is zero).
In particular, the skeletal filtration 
produces a $(d+1)$-dimensional limit cube
of functors whose vertices are given
by composites of the localization functors: \[
[0 \le k_0 < k_1 < \ldots < k_m \le d] \;\mapsto\; L_{k_0} L_{k_1} \ldots L_{k_m} \;.
\] For a fixed $k \in \{0, 1, \ldots, d\}$,
the functor $L_k$ is given by localizing from $\mod_{A_\bs}$ to $(\mod_{A_\bs})_{\gr, k}$ and embedding back via the right adjoint,
which is naturally isomorphic to \[
\Biggl( \prod_{\substack{\p \in \spec(A) \\ \dim(A/\p) = k}} \leftcomp A_\p \Biggr)_\bullet \otimes_{A_\bs} - \;\;\;\simeq\;\; \bigl( A^\thicksim_{S_k} \otimes_{A_\bs} - \bigr)^\thicksim_{S^c_{k-1}}
\] by Proposition \ref{prop-strata-as-analytic-rings}.

\begin{theorem}\label{thm-adeles}
The vertices of the cube are given by Beilinson-Parshin adeles. Specifically, each vertex is a continuous endofunctor that
is additionally right t-exact and preserves right t-bounded products and completions, and its value on $A$ is given by the formula
\begin{equation}\label{eq-main-thm-BP-adeles}
L_{k_0} L_{k_1} \ldots L_{k_m} (A) \;\simeq\; \prod_{\substack{\p_m \\ \dim A/\p_m = k_m}} \leftcomp\: (\ldots \prod_{\substack{\p_1\supset \p_2\\\dim A/\p_1 = k_1}}\leftcomp\:(\prod_{\substack{\p_0 \supset \p_1\\\dim A/\p_0 = k_0}} \leftcomp A_{\p_0})_{\p_1} \ldots )_{\p_m} \;.
\end{equation}
\end{theorem}

\sssec{Proof of the main theorem}\label{sssec-proof-of-the-main-thm}

\begin{definition}
For easier book-keeping, for an integer $k = 0, \ldots, d$ and a specialization closed subset $T \in \spec(A)^\spc$,
we define an ednodufunctor \[\xymatrix{
L_k^T \colon \mod_{A_\bs} \ar[rr] && \mod_{A_\bs} \\
M\ar@{}[u]|\invertin \ar@{|->}[rr] && \bigl( A^\thicksim_{S_k \cap T} \otimes_{A_\bs} M \bigr)^\thicksim_{S^c_{k-1}} \ar@{}[u]|\invertin &.
}\] Note that by putting $T = \spec(A)$ we recover the functor $L_k$: \[
L_k^{\spec(A)} \;\simeq\; L_k \;.
\]
\end{definition}

\begin{lemma}\label{lem-proof-of-the-main-thm-properties-of-L}
For any $T$ and $k$, the functor $L_k^T$ is continuous, right t-exact, and preserves right t-bounded products and completions.
\end{lemma}
\begin{proof}
First, continuity immediately follows from the formula \[
L_k^T(M) \;=\; \bigl( A^\thicksim_{S_k \cap T} \otimes_{A_\bs} M \bigr)^\thicksim_{S^c_{k-1}} \;,
\] where we use that completion at a generalization closed subset is always continuous, see the discussion at the beginning of
section \ref{ssec-gen-closed-and-ultrasolid}.
Moreover, the object $A^\thicksim_{S_k \cap T} \in \mod_{A_\bs}$ is compact by Proposition \ref{prop-comp-alg-is-idmpt-and-cmp},
and therefore the functor \[
A^\thicksim_{S_k \cap T} \otimes_{A_\bs} - : \mod_{A_\bs} \;\arr\; \mod_{A_\bs}
\]
preserves right t-bounded products and completions by Lemma \ref{lem-tensor-vs-prods}
and Lemma \ref{lem-tensor-pres-comp}.

Finally, because connective part $\mod^{\le 0}_{A_\bs}$ is compactly generated by products of copies of $A$,
to prove that the functor $L_k^T$ is right t-exact, it suffices to prove that its value on a generator is connective.
Moreover, because the functor $L_k^T$ preserves right t-bounded products and those are t-exact, it suffices to prove that
the value of our functor on $A$ is connective:
\[\begin{split}
L_k^T(A) &\;\simeq\; L_k(A^\thicksim_T) \;\simeq\;  \prod_{\substack{\p \in T \\ \dim(A/\p) = k}} \leftcomp A_\p \;,
\end{split}\]
where we use Corollary \ref{cor-single-comp-mod-div-skeleton} to compute the value
of $L_k( A^\thicksim_T)$.
Now the right t-exactness follows immediately from the first claim of Proposition \ref{prop-ultrasolid-stratum-t-exactness}.
\end{proof}

\begin{lemma}\label{lem-proof-of-the-main-thm-sup-condition}
An inclusion $T_1 \subseteq T_2$ of specialization closed subsets induces a map $A^\thicksim_{T_2} \arr A^\thicksim_{T_1}$, that in turn
induces a natural transformation \[
L_k^{T_2} \;\arr\; L_k^{T_1} \;.
\] For any module $M \in A^\thicksim_{T_1} - \mod_{A_\bs} \subseteq \mod_{A_\bs}$, the induced map \[
L_k^{T_2}(M) \;\arr\; L_k^{T_1}(M)
\] is an isomorphism.
\end{lemma}
\begin{proof}
Immediately follows from the commutativity of the diagram \[\xymatrix{
A^\thicksim_{S_k \cap T_2} \otimes_{A_\bs} M \ar[rr]\ar@{}[d]|\invsimeq_{\text{(Lem \ref{lem-sm-sp-map-preserves-finite-joins}) }} && A^\thicksim_{S_k \cap T_1} \otimes_{A_\bs} M \ar@{}[d]|\invsimeq^{\text{ (Lem \ref{lem-sm-sp-map-preserves-finite-joins})}}\\
A^\thicksim_{S_k} \otimes_{A_\bs} A^\thicksim_{T_2} \otimes_{A_\bs} M \ar[rr]\ar@{}[d]|\invsimeq_{(M \in A^\thicksim_{T_2}-\mod_{A_\bs})\;\:} && A^\thicksim_{S_k} \otimes_{A_\bs} A^\thicksim_{T_1} \otimes_{A_\bs} M \ar@{}[d]|\invsimeq^{ \:\;(M \in A^\thicksim_{T_1}-\mod_{A_\bs})} \\
A^\thicksim_{S_k} \otimes_{A_\bs} M \ar[rr]^= && A^\thicksim_{S_k} \otimes_{A_\bs} M &.
}\]
\end{proof}

\begin{proof}[Proof of Theorem \ref{thm-adeles}]
The continuity, the t-exactness, and the product and completion preservation are all contained in
Lemma \ref{lem-proof-of-the-main-thm-properties-of-L}. As for the
formula (\ref{eq-main-thm-BP-adeles}), it will be convenient to prove the following slightly more general version of it:
for any prime $\q \subset A$,
\[
L^{\spec(A/\q)}_{k_0} L^{\spec(A/\q)}_{k_1} \ldots L^{\spec(A/\q)}_{k_m} (A) \;\simeq\; \prod_{\substack{\p_m \supseteq \q \\ \dim A/\p_m = k_m}} \leftcomp\: (\ldots \prod_{\substack{\p_1\supseteq \p_2\\\dim A/\p_1 = k_1}}\leftcomp\:(\prod_{\substack{\p_0 \supseteq \p_1\\\dim A/\p_0 = k_0}} \leftcomp A_{\p_0})_{\p_1} \ldots )_{\p_m} \;.
\]
We will use induction on $m$. In the base case $m = 0$ we have: \[\begin{split}
L^{\spec(A/\q)}_k (A) &\;\simeq\; L_k (A^\comp_\q) \;\simeq\; \\
&\;\simeq\; \prod_{\substack{\p \in \spec(A/\q) \\ \dim A/\p = k}} \leftcomp A_\p \;\simeq\; \\
&\;\simeq\; \prod_{\substack{\p \supseteq \q \\ \dim A/\p = k}} \leftcomp A_\p \;,
\end{split}\] where we use Corollary \ref{cor-single-comp-mod-div-skeleton} to compute the value of $L_k(A^\comp_\q)$.
We will now prove the inductive step.
\begin{equation*}
\begin{split}
L^{\spec(A/\q)}_{k_0} L^{\spec(A/\q)}_{k_1} \ldots L^{\spec(A/\q)}_{k_m} (A) \;&\stackrel{\text{Base}}{\simeq}\; L^{\spec(A/\q)}_{k_0} L^{\spec(A/\q)}_{k_1} \ldots L^{\spec(A/\q)}_{k_{m-1}} \Bigl( \prod_{\substack{\p_m \supseteq \q \\ \dim A/\p_m  = k_m}} \leftcomp A_{\p_m} \Bigr) \;\simeq\\
&\stackrel{\text{Lemma \ref{lem-proof-of-the-main-thm-properties-of-L}}}{\simeq}\; \prod_{\substack{\p_m \supseteq \q \\ \dim A/\p_m = k_m}} L^{\spec(A/\q)}_{k_0} L^{\spec(A/\q)}_{k_1} \ldots L^{\spec(A/\q)}_{k_{m-1}} (\:\leftcomp A_{\p_m}) \;\simeq\\
&\stackrel{\text{Lemma \ref{lem-proof-of-the-main-thm-sup-condition}}}{\simeq}\; \prod_{\substack{\p_m \supseteq \q \\ \dim A/\p_m = k_m}} L^{{\spec(A/\p_m)}}_{k_0} L^{\spec(A/\p_m)}_{k_1} \ldots L^{\spec(A/\p_m)}_{k_{m-1}} (\:\leftcomp A_{\p_m}) \;\simeq\\
&\stackrel{\text{Lemma \ref{lem-proof-of-the-main-thm-properties-of-L}}}{\simeq}\; \prod_{\substack{\p_m \supseteq \q \\ \dim A/\p_m = k_m}} \leftcomp\:\Bigl( L^{\spec(A/\p_m)}_{k_0} L^{\spec(A/\p_m)}_{k_1} \ldots L^{\spec(A/\p_m)}_{k_{m-1}} (A) \Bigr)_{\p_m} \;\simeq\\
\stackrel{\substack{\text{induction}\\\text{hypothesis}}}{\simeq}\; &\prod_{\substack{\p_m \supseteq \q \\ \dim A/\p_m = k_m}} \leftcomp\: ( \prod_{\substack{\p_{m-1} \supseteq \p_m \\ \dim A/\p_{m-1} = k_{m-1}}} \leftcomp\:(\ldots \prod_{\substack{\p_1\supseteq \p_2\\\dim A/\p_1 = k_1}}\leftcomp\:(\prod_{\substack{\p_0 \supseteq \p_1\\\dim A/\p_0 = k_0}} \leftcomp A_{\p_0})_{\p_1} \ldots )_{\p_{m-1}} )_{\p_m} \;.
\end{split}
\end{equation*}
\end{proof}

\sssec{} Let us also make a parenthetical remark that the proof of Theorem \ref{thm-adeles} almost without change applies to
$L_{k_0} L_{k_1} \ldots L_{k_m} (A^\thicksim_T)$ with $T$ being a specialization closed subset, and gives a formula \[
L_{k_0} L_{k_1} \ldots L_{k_m} (A^\thicksim_T) \;\simeq\; \prod_{\substack{\p_m \in T \\ \dim A/\p_m = k_m}} \leftcomp\: (\ldots \prod_{\substack{\p_1\supset \p_2\\\dim A/\p_1 = k_1}}\leftcomp\:(\prod_{\substack{\p_0 \supset \p_1\\\dim A/\p_0 = k_0}} \leftcomp A_{\p_0})_{\p_1} \ldots )_{\p_m} \;.
\]

\ssec{Comparison to the classical definition of Beilinson-Parshin}\label{ssec-classical-BP}
In this parenthetical section we explain how the formula of our Theorem \ref{thm-adeles} agrees with the classical cosimplicial
inductive formula of Beilinson and Parshin, originally introduced in \cite{ParshinAdeles} and \cite{BeilinsonResidues}.
First, we recall the original definition. As usual, $A$ denotes a
commutative ring of finite type over $\Z$. Let $S^{\mathrm{red}}_\cdot$ denote the semisimplicial set with $r$-simplices defined
by the formula \[
S_r^{\mathrm{red}} \;=\; \{ \p_0, \ldots, \p_r \in \spec(A) : \p_i \supset \p_{i+1} \} \;,
\] where by $\p \supset \q$ we mean that the prime $\p$ contains $\q$ but is not equal to it.
The face maps defined by skipping a point in the flag. For a subset $T \subseteq S_r^{\mathrm{red}}$, one can define a 
functor \[
\mathbb{A}(T, -) \colon \mod_{A_\bs}^{< +\infty} \arr \mod_{A_\bs}^{< +\infty} 
\] that
\begin{enumerate}
\item is exact, continuous, and satisfies \[
\mathbb{A}(T, \prod_J A) \;\simeq\; \prod_J \mathbb{A}(T, A) \;;
\]
\item in case $r = 0$, satisfies \[
\mathbb{A}(T, A) \;=\; \prod_{\q \in T} \leftcomp A_\q \;;
\]
\item in case $r > 0$, satisfies \[
\mathbb{A}(T, A) \;=\; \prod_{\q \in T} \leftcomp\: \Bigl(\mathbb{A}(T_\q, A) \Bigr)_\q \;,
\] where \[
T_\q \subseteq S^\mathrm{red}_{r-1}
\] consists of flags $\p_0 \supset \ldots \supset \p_{r-1}$ satisfying
$(\p_0 \supset \ldots \supset \p_{r-1} \supset \q) \in T$.
\end{enumerate}
A direct analog of the Beilinson-Parshin resolution (see \cite{BeilinsonResidues} or Corollary 8.13 in
\cite{morrow2012introductionhigherdimensionallocal}) would be the following statement.
\begin{proposition}\label{classical-adeles-prop}
For an object $M \in \mod^{< +\infty}_{A_\bs}$, the canonical map \[
M \;\arr\; \lim_{r \in \Delta_s^{\op}} \mathbb{A}(S_r^{\mathrm{red}}, M)
\] is an isomorphism.
\end{proposition}
We first prove the following lemma.
\begin{lemma}\label{classical-adeles-local-lem}
For a prime ideal $\q \subset A$, we have an isomorphism
\[
\mathbb{A}((S_{r+1}^{\mathrm{red}})_\q, A) \;\simeq\; \prod_{\substack{\p_r \supset \q}} \leftcomp\: (\ldots \prod_{\substack{\p_1\supset \p_2}}\leftcomp\:(\prod_{\substack{\p_0 \supset \p_1}} \leftcomp A_{\p_0})_{\p_1} \ldots )_{\p_m} \;.
\]
\end{lemma}
\begin{proof}
This is easily done by induction on $r$. The base is clear, assume that $r > 0$. Using the definition, we get
\[
\mathbb{A}((S_{r+1}^{\mathrm{red}})_\q, A) \;\simeq\; \prod_{\substack{\p \supset \q}} \leftcomp\: \Bigl[ \mathbb{A}\bigl(((S_{r+1}^{\mathrm{red}})_\q)_\p,A\bigr) \Bigr]_\q \;.
\] Then we note that, for a prime ideal $\q$ which is contained in $\p$ but not equal to $\p$, we have an equality \[
((S_{r+1}^{\mathrm{red}})_\q)_\p \;=\; (S_r^\mathrm{red})_\p \;.
\] Now the statement immediately follows from the induction hypothesis.
\end{proof}
\begin{proof}[Proof of Proposition \ref{classical-adeles-prop}]
It suffices to prove the claim in the case $M = A$.
By Lemma \ref{classical-adeles-local-lem}, we have an isomorphism \[
\mathbb{A}(S_r^{\mathrm{red}}, A) \;\simeq\; \prod_{[0 \le k_0 < \ldots < k_r \le d]} \prod_{\substack{\p_r \\ \dim A/\p_r = k_r}} \leftcomp\: (\ldots \prod_{\substack{\p_1\supset \p_2\\\dim A/\p_1 = k_1}}\leftcomp\:(\prod_{\substack{\p_0 \supset \p_1\\\dim A/\p_0 = k_0}} \leftcomp A_{\p_0})_{\p_1} \ldots )_{\p_m} \;.
\] Then the claim immediately follows from our Theorem \ref{thm-adeles} combined with
the comparison between cubical and semi-cosimplicial limits, which is recorded in
Lemma 2.2.3 and Proposition 2.2.4 of \cite{kim2021adelicdescentktheory}.
\end{proof}

%% file: global.tex
\section{The main theorem: globalization}

Let us fix a scheme $X$ of finite type over $\Z$. Clausen-Scholze define the category $\qcoh(X_\bs)$
of solid quasi-coherent sheaves on $X$ and prove that it satisfies Zariski descent
(see Theorem 9.8 in \cite{condensed}), which we now use to construct the skeletal
filtration on $\qcoh(X_\bs)$ and prove a version of the main theorem for it.

\ssec{Unions of tubular neighborhoods and the smashing spectrum: globalization}

Recall that, for an open embedding $j \colon U \hookrightarrow X$, one has a symmetric monoidal functor \[
j^* \colon \qcoh(X_\bs) \;\arr\; \qcoh(U_\bs) \;.
\] It admits a fully faithful $\qcoh(X_\bs)$-linear left adjoint $j_!$ and a fully faithful continuous right adjoint $j_*$.
In particular, the natural transformation \[
j^* \uhom_{X_\bs}(- , -) \;\arr\; \uhom_{U_\bs} (j^*-, j^*-) \;,
\] which one gets from the adjunction, is an isomorphism.

We will denote by $\sm(X_\bs)$ the smashing spectrum of the category $\qcoh(X_\bs)$, and likewise for $\sm(X_\bs)^\cmp$.
\begin{lemma}
For an open embedding $j \colon U \hookrightarrow X$, there is an induced map of posets \[\begin{split}
&\sm(X_\bs)^\cmp \;\arr\; \sm(U_\bs)^\cmp 
\end{split}\] taking an idempotent algebra into its image under $j^*$; this map preserves finite meets and finite joins.
\end{lemma}

\sssec{} We denote by $X^\spc$ the poset of specialization closed subsets in $X$. Given a specialization closed subset $T \in X^\spc$,
as in the affine case, one defines the functor of completion at $T$ by the formula \[
(-)^\thicksim_T \;\stackrel{\define}{=}\; \uhom_{X_\bs} (\Ga_T, -) \;,
\] where $\Ga_T \;\simeq\; \stackrel[Z \subseteq T]{}{\colim}\Ga_Z$ is the filtered colimit of local cohomologies with
supports in closed subsets $Z$ taken over the poset of all $Z$ contained in $T$.

\begin{lemma}\label{lem-pull-back-preserves-comp}
For an open embedding $j\colon U \hookrightarrow X$ and a specialization closed subset $T \in X^\spc$, there is a natural isomorphism
\[
j^* (-)^\thicksim_T \;\simeq\; (j^*-)^\thicksim_{T\cap U} \;.
\]
\end{lemma}
\begin{proof}
Using the definition and the fact that $j^*$ respects the inner hom, we get: \[\begin{split}
j^* (-)^\thicksim_T &\;\simeq\; j^*\uhom_{X_\bs} (\Ga_T, -) \;\simeq\; \uhom_{X_\bs} (j^*\Ga_T, j^*-) \;.
\end{split}\] It remains to observe that $j^* \Ga_T \;\simeq\; \Ga_{T\cap U}$, which is manifest
in view of the fact that $j^*$ restricts to the usual pull-back functor on the full subcategory
of discrete objects $\qcoh(X) \subset \qcoh(X_\bs)$.
\end{proof}

\begin{proposition}
There is a map of posets \begin{equation}\label{the-sm-sp-map-global}
(X^\spc)^\op \;\arr\; \sm(X_\bs)^\cmp
\end{equation}
which is uniquely determined by the condition that it restricts to the map (\ref{the-smashing-spectrum-map}) for any affine open
$j \colon \spec(A) \hookrightarrow X$: \[\xymatrix{
(X^\spc)^\op \ar[rr]\ar[d]^{j^*} && \sm(X_\bs)^\cmp \ar[d]^{j^*} \\
(\spec(A)^\spc)^\op \ar[rr] && \sm(A_\bs)^\cmp &.
}\]
\end{proposition}
\begin{proof}
As in the affine case, the map is given by \[\xymatrix{
T \ar@{|->}[rr] && (\O_X)^\thicksim_T \;\simeq\; \uhom_{X_\bs} (\Ga_T, \O_X) \;.
}\] Using Lemma \ref{lem-pull-back-preserves-comp} and Proposition \ref{prop-the-smashing-spectrum-map},
one can check locally in Zariski topology on $X$ that the RHS is compact and idempotent.
The uniqueness follows from the Zariski descent for $\qcoh(X_\bs)$.
\end{proof}

\sssec{} We denote by $X^\gen$ the poset of generalization closed subsets of $X$. As in the affine case,
for a generalization closed subset $E \in X^\gen$, we define the functor of completion at $E$ by the formula \[
(-)^\thicksim_E \;\stackrel{\define}{=}\; \uhom_{X_\bs} ((\O_X)^\thicksim_{E^c} / \O_X[-1], -) \;,
\] where $E^c = X \setminus E \in X^\spc$ is the complement.
\begin{lemma}\label{lem-pull-back-preserves-comp-gen}
For an open embedding $j\colon U \hookrightarrow X$ and a generalization closed subset $E \in X^\gen$, there is a natural isomorphism
\[
j^* (-)^\thicksim_E \;\simeq\; (j^*-)^\thicksim_{E\cap U} \;.
\]
\end{lemma}
\begin{proof}
Using the definition, the fact that $j^*$ respects the inner $\hom$, and Lemma \ref{lem-pull-back-preserves-comp},
we get: \[\begin{split}
j^* (-)^\thicksim_E &\;\simeq\; j^*\uhom_{X_\bs} ((\O_X)^\thicksim_{E^c} / \O_X [-1], -) \;\simeq\; \\
&\;\simeq\; \uhom_{X_\bs} (j^*( (\O_X)^\thicksim_{E^c} / \O_X [-1] ), j^*-) \;\simeq\; \\
&\;\simeq\; \uhom_{X_\bs} ( (\O_U)^\thicksim_{E^c \cap U} / \O_U [-1] , j^*-) \;\simeq\; \\
&\;\simeq\; (j^*-)^\thicksim_{E \cap U} \;.
\end{split}\]
\end{proof}

\sssec{} Given a map of posets \[\xymatrix{
[0; d] \ar[rr] && X^\spc \\
k \ar@{}[u]|\invertin \ar@{|->}[rr] && T_{k-1} \ar@{}[u]|\invertin
}\] satisfying $T_{-1} = \emptyset$,
by composing it with the map (\ref{the-sm-sp-map-global}) and the map (\ref{smashing-spectrum-to-closed-subcat}),
we get a stratification \[
[0;d]^\op \;\arr\; \cls(\qcoh(X_\bs)) \;.
\]
In particular, we get functors \[
L_k(-) \;\simeq\; \bigl( (\O_X)^\thicksim_{T_k} \otimes_{X_\bs} - \bigr)^\thicksim_{T_{k-1}^c}
\] that localize to the corresponding associated graded factors.

\begin{lemma}\label{lem-L-functors-Zar-locally}
Fix an open embedding $j \colon U \hookrightarrow X$ and denote by $L_{U, k}$ the endofunctor of $\qcoh(U_\bs)$ defined by the formula
\[
L_{U,k}(-) \;\simeq\; \bigl( (\O_U)^\thicksim_{U \cap T_k} \otimes_{U_\bs} - \bigr)^\thicksim_{U \cap T_{k-1}^c} \;.
\]
The map \[
j^* L_k \;\arr\; L_{U, k} j^* \;,
\] which we get from the adjunction, is an isomorphism.
\end{lemma}
\begin{proof}
Follows immediately from Lemma \ref{lem-pull-back-preserves-comp} and Lemma \ref{lem-pull-back-preserves-comp-gen}.
\end{proof}

\ssec{The skeletal filtration}
As before, $X$ is a scheme of finite type over $\mathbb{Z}$ of dimension $d$.
Let $[0,d]$ denote the poset $0 < 1 < 2 < \ldots < d$; we define a map of posets
\begin{equation}\label{eq-skeleton-as-spc-closed-subsets-global}\xymatrix{
[0,d] \ar[rr] && X^\spc \\
k \ar@{}[u]|\invertin \ar@{|->}[rr] && \sk_{k-1} \stackrel{\define}{=} \stackrel[\substack{Z \subseteq X \\ \dim(Z) \le k-1}]{}{\cup} Z \ar@{}[u]|\invertin &.
}\end{equation} The \textit{skeletal filtration} \[
[0,d]^\op \;\arr\; \cls(\qcoh(X_\bs))
\] is a stratification of $\qcoh(X_\bs)$ defined by composing the map (\ref{eq-skeleton-as-spc-closed-subsets-global})
with maps (\ref{the-sm-sp-map-global}) and (\ref{smashing-spectrum-to-closed-subcat}). As was explained in section
\ref{ssec-stratifications-and-sm-sp}, we get endofunctors \[
L_k(-) \;\simeq\; \bigl( (\O_X)^\thicksim_{S_k} \otimes_{X_\bs} - \bigr)^\thicksim_{S_{k-1}^c} \;, \qquad k\in\{0, 1, \ldots, d\},
\] whose composites comprise a cubical diagram \[
[0 \le k_0 < \ldots < k_m \le d] \;\mapsto\; L_{k_0} \ldots L_{k_m} \;,
\] and the induced map \[
\id \;\arr\; \stackrel[ \substack{[0 \le k_0 < \ldots < k_m \le d] \\ m \ge 0} ]{}{\lim} \;L_{k_0} \ldots L_{k_m}
\] is an isomorphism by Proposition \ref{resolution-stratified-over-poset}.
Now our goal is to express the vertices in terms of nested completed localizations---as it was done in the affine case.

\begin{theorem}\label{thm-global-vertices}
The endofunctors $L_k \colon \qcoh(X_\bs) \arr \qcoh(X_\bs)$ are continuous, right t-exact, and preserve right t-bounded products and
completions. The value of $L_{k_0} \ldots L_{k_m}$ on a perfect complex $\P$ is given by the formula
\begin{equation}\label{eq-main-thm-BP-adeles-global}
L_{k_0} \ldots L_{k_m} (\P) \;\simeq\; \prod_{\substack{\p_m \in X \\ \dim \overline{\{\p_m\}} = k_m}} \leftcomp\: (\ldots \prod_{\substack{\p_1\supset \p_2\\\dim \overline{\{\p_1\}} = k_1}}\leftcomp\:(\prod_{\substack{\p_0 \supset \p_1\\\dim \overline{\{\p_0\}} = k_0}} \leftcomp\: \P_{\p_0})_{\p_1} \ldots )_{\p_m} \;.
\end{equation}
\end{theorem}
\begin{proof}
First, we note that each functor $L_k$ is continuous by construction and right t-exact by Lemma \ref{lem-L-functors-Zar-locally} and
Theorem \ref{thm-adeles}. 
We pick a finite affine cover $X = \cup_\a U_\a$ such that \[
\qcoh(X_\bs) \;\stackrel{\simeq}{\arr}\;\; \stackrel[\a]{}{\lim}^* \qcoh({U_\a}_\bs) \;,
\] where the limit is computed using the upper-$*$ pull-backs and the map to the limit is also induced by the upper-$*$ pull-backs.
Using the descent and Lemma \ref{lem-L-functors-Zar-locally}, we write \[\begin{split}
L_k &\;\simeq\; \stackrel[\a]{}{\lim} {\phi_\a}_* \phi_\a^* L_k \;\simeq\\
&\;\simeq\; \stackrel[\a]{}{\lim} {\phi_\a}_* L_{U_\a, k} \phi_\a^* \;,
\end{split}\] where $\phi_\a$ denotes the open embedding $U_\a \hookrightarrow X$.
Now the product and completion preservation follow from the corresponding affine statement plus the fact that functors
${\phi_\a}_*$ and $\phi_\a^*$ are t-exact and preserve all limits.

Let us denote the nested product on the right hand side of (\ref{eq-main-thm-BP-adeles-global}) by $L^{\text{BP}}_{k_0\ldots k_m}\P$.
By Theorem \ref{thm-adeles}, and continuity and co-continuity of each $\fakej^*_\a$, for each $\a \in \A$, there is an isomorphism \[
\phi_\a^* \Psi^{\text{BP}}_{k_0\ldots k_m}\P \;\simeq\; L_{U_\a, k_0} \ldots L_{U_\a, k_m} \phi^*_\a \P \;.
\]
On the other hand, by the descent and Lemma \ref{lem-L-functors-Zar-locally},
\begin{equation*}
\begin{split}
L^{\text{BP}}_{k_0\ldots k_m}\P &\;\stackrel{\simeq}{\arr}\; \stackrel[\a]{}{\lim} {\phi_\a}_*\phi^*_\a L^{\text{BP}}_{k_0\ldots k_m}\P \;\simeq\\
&\;\simeq\; \stackrel[\a]{}{\lim} {\phi_\a}_*L_{U_\a, k_0} \ldots L_{U_\a, k_m} \phi^*_\a \P \;\simeq\\
&\;\simeq\; L_{k_0} \ldots L_{k_m} \P \;.
\end{split}
\end{equation*}
\end{proof}

%% file: canonical_colim_of_affs.tex
\section{Solid adelic descent}\label{5th_sec}

Let $X$ be a scheme of finite type over $\Z$.
As a geometric application of our main results,
we will describe $X_\bs$ as a colimit (in the Tannakian sense) of a certain canonically defined cubical
diagram of affinoids, whose terms will be given by the adelic rings with the analytic ring structure inherited from $X_\bs$.
This can be viewed as a generalization or a solid analog of Theorem 3.1 \cite{Groechenig_2017}.
We note that, unlike the decomposition of the category $\qcoh(X_\bs)$ one gets by directly applying the results of \cite{stratified}
to the skeletal filtration, the decomposition we describe in Theorem \ref{adelic-cover} of this section
is geometric---its pieces are quasi-coherent sheaves on certain analytic spaces
(in fact, affinoid, as we prove in Theorem \ref{thm-adeles-are-affinoid}), and the functors between the pieces are symmetric monoidal.

\ssec{Some properties of the completion at the skeleton}\label{ssec-some-prop-of-skeleton-comp}
As usual, we fix a ring $A$ of finite type over $\Z$.
Recall from section \ref{ssec-skeletal-filtration} the definition of the $k$-th skeleton $S_k \in \spec(A)^\spc$.
Also recall from section \ref{sssec-completion-at-a-specialization-closed} that,
for a specialization closed subset $T \in \spec(A)^\spc$, we have a functor of completion at $T$,
which is an idempotent lax-symmetric monoidal endofunctor of $\mod_{A_\bs}$.
We now study in greater detail the functor of completion at $S_k$, which is given by the formula \[\xymatrix{
\mod_{A_\bs} \ar[rrr] &&& \mod_{A_\bs} \\
M \ar@{}[u]|\invertin \ar@{|->}[rr] && M^\thicksim_{S_k} \ar@{}[r]|(0.3){\stackrel{}{\simeq}} & \lim_{\substack{Y \subseteq \spec(A) \\ \dim(Y) \le k}} M^\thicksim_Y \ar@{}[u]|\invertin \;,
}\] where the limit is taken over the poset of all closed subsets $Y \subseteq \spec(A)$ of dimension not greater than $k$.
We postpone some of the proofs until Appendix \ref{app-k-comp}.
\begin{definition}
For an integer $k$, an $A_\bs$-module $M$ is called $S_k$-complete if the map \[
M \arr M^\thicksim_{S_k} \;,
\] is an isomorphism. We note that, since the functor of completion at a specialization closed subset is idempotent,
an object of the form $M^\thicksim_{S_k}$ is $S_k$-complete.
\end{definition}
The following lemma provides a few non-tautological examples of $S_k$-complete modules.
\begin{lemma}\label{lem-some-exms-of-k-complete}
Any module $M \in \mod_{A_\bs}$ complete with respect to a specialization closed subset $T \in \spec(A)^\spc$ satisfying 
$T \cap S_k = T$ is $S_k$ complete. In particular, a module complete with respect to an ideal $I \subset A$ such that $\dim(A/I) \le k$
is $S_k$ complete. Also, an $S_{k'}$-complete module is $S_k$-complete if $k' \le k$.
\end{lemma}
\begin{proof}
We need to prove that the map \[
M \;\arr\; \uhom_{A_\bs} (\Ga_{S_k}, M)
\] is an isomorphism. Using the $T$-completeness, it suffices to prove that the map \[
M^\thicksim_T \;\arr\; \uhom_{A_\bs} (\Ga_{S_k}, M^\thicksim_T)
\] is an isomorphism, which easily follows from Lemma \ref{lem-sm-sp-map-preserves-finite-joins}: \[\begin{split}
\uhom_{A_\bs} (\Ga_{S_k}, M^\thicksim_T) &\;\simeq\; \uhom_{A_\bs} (\Ga_{S_k}, \uhom_{A_\bs}(\Ga_T, M)) \;\simeq\; \\
&\;\simeq\; \uhom_{A_\bs} (\Ga_{S_k}\otimes_{A_\bs} \Ga_T, M) \;\simeq\; \\
&\stackrel{\text{Lemma \ref{lem-sm-sp-map-preserves-finite-joins}}}{\;\simeq\;} \uhom_{A_\bs} (\Ga_{S_k\cap T}, M) \;\simeq\; \uhom_{A_\bs} (\Ga_{T}, M) \;.
\end{split}\]
\end{proof}

\begin{corollary}
The ring $L_{k_0} L_{k_1} \ldots L_{k_m}(A)$ is $S_{k_m}$-complete.
\end{corollary}
\begin{proof}
First, we note that, by Theorem \ref{thm-adeles}, the ring at hand splits as a product: \[
L_{k_0} L_{k_1} \ldots L_{k_m}(A) \;\simeq\;\ \prod_{\substack{\p_m \in \spec(A) \\ \dim(A/\p_m) = k_m}} M(\p_m) \;,
\] where the module $M(\p_m)$ is complete at $p_m$. Therefore, because the subcategory of $S_{k_m}$-complete modules is closed
under all limits, it suffices to prove that each $M(\p_m)$ is $S_{k_m}$-complete, which immediately follows
from Lemma \ref{lem-some-exms-of-k-complete}.
\end{proof}

We finish this section with a couple of useful technical statements.
\begin{proposition}\label{prop-primes-are-conservative-on-k-complete}
The set of functors \[
\{ (-)^\comp_\p \colon \mod_{A_\bs} \arr \mod_{A_\bs}\}_{\substack{\p \in \spec(A) \\ \dim(\p) \le k}}
\] is jointly conservative on the full subcategory of $S_k$-complete modules.
\end{proposition}
\begin{proof}
It follows from Proposition \ref{prop-lim-over-primes} that the natural map \[
M^\thicksim_{S_k} \;\arr\; \stackrel[\substack{\p \in \spec(A) \\ \dim(A/\p) \le k}]{}{\lim} M^\comp_\p
\] is an isomorphism, from which the statement follows immediately.
\end{proof}
We postpone the proof of the following proposition and its corollary until appendix \ref{app-k-comp},
see Proposition \ref{k-comp-is-monoidal} and Corollary \ref{appB-sing-ideal-comp-is-monoidal}.
\begin{proposition}\label{prop-comp-is-monoidal}
The functor of $S_k$-completion (which is a priori lax symmetric monoidal) is symmetric monoidal when restricted to the
subcategory of eventually connective modules, i.e. for any pair $M,N \in \mod^{< +\infty}_{A_\bs}$, the map
\[
M^\thicksim_{S_k} \otimes_{A_\bs} N^\thicksim_{S_k} \;\arr\; \bigl( M\otimes_{A_\bs} N \bigr)^\thicksim_{S_k}
\]
is an isomorphism. In particular, the tensor product of a pair of $S_k$-complete eventually connective $A_\bs$-modules is $S_k$-complete.
\end{proposition}
\begin{corollary}[Proposition A.17 in \cite{Artem}]\label{cor-comp-sing-ideal-is-monoidal}
The functor \[
(-)^\comp_I \colon \mod^{< +\infty}_{A_\bs} \;\arr\; \mod^{< +\infty}_{A_\bs}
\] of completion at an ideal $I\subset A$ is symmetric monoidal. In particular,
the tensor product of a pair of $I$-complete eventually connective $A_\bs$-modules is $I$-complete.
\end{corollary}

\ssec{Solid adelic descent}
In this section we state and prove a descent statement, which can be viewed as a solid analog of Theorem 3.1 in \cite{Groechenig_2017}.
Recall that $X$ denotes a scheme of finite type over $\Z$. First, we establish certain nice properties of the adelic rings in the solid setting.

\begin{proposition}\label{factors-are-tensor-products}
For an ordered tuple $[0 \le k_0 < \ldots < k_m \le \dim(X)]$, the canonical map
\[ L_{k_0}\O_X \otimes_{X_\bs} \ldots \otimes_{X_\bs} L_{k_m}\O_X \arr L_{k_0}\ldots L_{k_m}\O_X \]
is an isomorphism.
\end{proposition}
\begin{proof}
It suffices to prove the statement in the case where $X = \spec(A)$ for a finite type algebra $A$.
Moreover, by induction, it suffices to prove that the map \begin{equation}\label{sec5.2-local-eq1}
L_{k_0}\ldots L_{k_{m-1}}A \otimes_{A_\bs} L_{k_m}A \;\arr\; L_{k_0}\ldots L_{k_m}A
\end{equation}
is an isomorphism.
Now since both the source and the target are $S_{k_m}$-complete by Lemma
\ref{lem-some-exms-of-k-complete} and Proposition \ref{prop-comp-is-monoidal},
by the joint conservativeness of Proposition \ref{prop-primes-are-conservative-on-k-complete},
it suffices to prove that (\ref{sec5.2-local-eq1}) becomes an isomorphism upon application of the functor of $\p$-completion for each prime
$\p \in \spec(A)$ of dimension not greater than $k_m$.

First, we compute the $\p$-completion of the left hand side of (\ref{sec5.2-local-eq1}). By Corollary \ref{cor-comp-sing-ideal-is-monoidal}, 
\[ \bigl(L_{k_0}\ldots L_{k_{m-1}}A \otimes_{A_\bs} L_{k_m}A \bigr)^\comp_\p \;\simeq\; \bigl(L_{k_0}\ldots L_{k_{m-1}}A \bigr)^\comp_\p \otimes_{A_\bs} \bigl( L_{k_m}A \bigr)^\comp_\p \;. \] Using Theorem \ref{thm-adeles}, \[
\bigl( L_{k_m}A \bigr)^\comp_\p \;\simeq\; \prod_{\substack{\p_m \\ \dim(A/\p_m) = k_m}} \bigl( \leftcomp A_{\p_m} \bigr)^\comp_\p \;\simeq\; \left\{ \begin{smallmatrix}\leftcomp A_\p , & \text{ if } \dim(A/\p) = k_m ; \\ 0, & \text{ otherwise.} \end{smallmatrix}\right.
\] Given that the result vanishes unless $\dim(A/\p) = k_m$, while computing the other multiple $\bigl(L_{k_0}\ldots L_{k_{m-1}}A \bigr)^\comp_\p$, we can assume that
$\dim(A/\p) = k_m$. Using Theorem \ref{thm-adeles}, \[\begin{split}
\bigl(L_{k_0}\ldots L_{k_{m-1}}A \bigr)^\comp_\p &\;\simeq\; \Bigl( \prod_{\substack{\p_{m-1} \\ \dim A/\p_{m-1} = k_{m-1}}} \leftcomp\: (\ldots \prod_{\substack{\p_1\supset \p_2\\\dim A/\p_1 = k_1}}\leftcomp\:(\prod_{\substack{\p_0 \supset \p_1\\\dim A/\p_0 = k_0}} \leftcomp A_{\p_0})_{\p_1} \ldots )_{\p_{m-1}} \Bigr)^\comp_\p \\
&\;\simeq\; \prod_{\substack{\p_{m-1} \supset \p \\ \dim A/\p_{m-1} = k_{m-1}}} \leftcomp\: (\ldots \prod_{\substack{\p_1\supset \p_2\\\dim A/\p_1 = k_1}}\leftcomp\:(\prod_{\substack{\p_0 \supset \p_1\\\dim A/\p_0 = k_0}} \leftcomp A_{\p_0})_{\p_1} \ldots )_{\p_{m-1}}
\end{split}\] Putting the last two computations together and using \ref{prop-comp-is-monoidal} one more time, we get \[
\bigl(L_{k_0}\ldots L_{k_{m-1}}A \otimes_{A_\bs} L_{k_m}A \bigr)^\comp_\p \;\simeq\; \leftcomp\:\Bigl( \prod_{\substack{\p_{m-1} \supset \p \\ \dim A/\p_{m-1} = k_{m-1}}} \leftcomp\:(\ldots \prod_{\substack{\p_1\supset \p_2\\\dim A/\p_1 = k_1}}\leftcomp\:(\prod_{\substack{\p_0 \supset \p_1\\\dim A/\p_0 = k_0}} \leftcomp A_{\p_0})_{\p_1} \ldots )_{\p_{m-1}} \Bigr)_\p
\]

Similarly, using Theorem \ref{thm-adeles}, we get that the $\p$-completion $\bigl( L_{k_0}\ldots L_{k_m}A \bigr)^\comp_\p$ of right hand side of
(\ref{sec5.2-local-eq1}) vanishes, unless $\dim(A/\p) = k_m$, and in the latter case is equal to \[
\leftcomp\:\Bigl( \prod_{\substack{\p_{m-1} \supset \p \\ \dim A/\p_{m-1} = k_{m-1}}} \leftcomp\:(\ldots \prod_{\substack{\p_1\supset \p_2\\\dim A/\p_1 = k_1}}\leftcomp\:(\prod_{\substack{\p_0 \supset \p_1\\\dim A/\p_0 = k_0}} \leftcomp A_{\p_0})_{\p_1} \ldots )_{\p_{m-1}} \Bigr)_\p \;.
\]
\end{proof}

\begin{corollary}\label{adelic-rings-are-idempotent}
The ring $L_{k_0}\ldots L_{k_m} \O_X$ is idempotent relative to $\qcoh(X_\bs)$.
\end{corollary}
\begin{proof}
By Proposition \ref{factors-are-tensor-products}, it suffices to prove that $L_k \O_X$ is an idempotent algebra for each $k$.
Moreover, by Zariski descent, it suffices to treat the affine case $X = \spec(A)$.
By Lemma \ref{lem-some-exms-of-k-complete} and Proposition \ref{prop-comp-is-monoidal},
in order to prove that the map \begin{equation}\label{sec5.2-local-eq2}
L_k A \otimes_{A_\bs} L_k A \;\arr\; L_k A
\end{equation}
is an isomorphism, it suffices to prove that it becomes so after the $\p$-completion for each prime $\p$ of dimension not greater than $k$.
Using Theorem \ref{thm-adeles},  Proposition \ref{prop-comp-is-monoidal}, and also the fact that completion preserves products, we get that
both the source and the target of (\ref{sec5.2-local-eq2})
vanish under the $\p$-completion unless $\dim(A/\p) = k$, and in the latter case (\ref{sec5.2-local-eq2})
turns into \[
\leftcomp A_\p \otimes_{A_\bs} \leftcomp A_\p \;\arr\; \leftcomp A_\p \;,
\] which is an isomorphism by Corollary \ref{cor-comp-sing-ideal-is-monoidal} and the fact that $A_\p$ is idempotent over $A_\bs$.
\end{proof}

\sssec{} We are now ready to state and prove the descent statement.
For an ordered tuple $[0 \le k_0 < \ldots < k_m \le \dim(X)]$, we consider a ring \[
L_{k_0}\ldots L_{k_m} \O_X \in \calg(\qcoh(X_\bs)) \;.
\] As before, the unit maps $\id \arr L_k$ produce us a map of rings \[
L_{k_0}\ldots L_{k_m} \O_X \arr L_{l_0}\ldots L_{l_n} \O_X
\] for each inclusion \[
[k_0 < k_1 < \ldots < k_m] \;\subset\; [l_0 < l_1 < \ldots < l_n] \;,
\] which together comprise a cubical diagram---see either Proposition \ref{resolution-stratified-over-poset} or 
Proposition \ref{resolution-filtration}.

\begin{theorem}[Solid adelic descent]\label{adelic-cover}
The map \begin{equation}\label{geometric-functor-1}
\qcoh(X_\bs) \arr \stackrel[\substack{[k_0 < k_1 < \ldots < k_m] \\ 0\le m}]{}{\lim} L_{k_0}\ldots L_{k_m} \O_X - \mod_{\qcoh(X_\bs)}
\end{equation}
is an equivalence.
\end{theorem}

\begin{proof}
Let us denote the functor (\ref{geometric-functor-1}) by $F^*$ and its right adjoint by $F_*$. Likewise, we write
$F^*_{k_0 \ldots k_m}$ for the functor
\[
\qcoh(X_\bs) \arr L_{k_0}\ldots L_{k_m} \O_X - \mod_{\qcoh(X_\bs)}
\]
and ${F_{k_0 \ldots k_m}}_*$ for its right adjoint.

Let us prove first that $F^*$ is fully faithful. The unit $\id \arr F_* F^*$ is given by
\begin{equation*}
\begin{split}
M &\;\mapsto\; \stackrel[\substack{[k_0 < \ldots < k_m]\\0 \le m}]{}{\lim} {F_{k_0 \ldots k_m}}_*F^*_{k_0 \ldots k_m}M \;\simeq\\
&\;\simeq\; \stackrel[\substack{[k_0 < \ldots < k_m]\\0 \le m}]{}{\lim} M\otimes_{X_\bs}{F_{k_0 \ldots k_m}}_*F^*_{k_0 \ldots k_m}\O_X \;\simeq\\
&\;\simeq\; M\otimes_{X_\bs}\stackrel[\substack{[k_0 < \ldots < k_m]\\0 \le m}]{}{\lim} L_{k_0}\ldots L_{k_m}\O_X \;.
\end{split}
\end{equation*}
So it suffices to prove that \[
\O_X \arr \stackrel[\substack{[k_0 < \ldots < k_m]\\0 \le m}]{}{\lim} L_{k_0}\ldots L_{k_m}\O_X
\] is an isomorphism, which is a direct corollary of Theorem \ref{thm-global-vertices}.

To finish the proof of the Theorem, we need to check that the counit $F^* F_* \arr \id$ is an isomorphism, but it follows immideately from Corollary \ref{adelic-rings-are-idempotent}.
\end{proof}

\begin{corollary}
The map \begin{equation*}
\perf(X) \arr \stackrel[\substack{[k_0 < k_1 < \ldots < k_m] \\ 0\le m}]{}{\lim} \perf (L_{k_0}\ldots L_{k_m} \O_X )
\end{equation*}
induced by the map (\ref{geometric-functor-1}) is an equivalence.
\end{corollary}
\begin{proof}
In this proof, we keep the notation introduced in the proof of Theorem \ref{adelic-cover}. We need to prove that, whenever we have
an object
\[\{ \P_{k_0 \ldots k_m}\}_{k_0 < \ldots < k_m} \in \stackrel[\substack{[k_0 < k_1 < \ldots < k_m] \\ 0\le m}]{}{\lim} \perf (L_{k_0}\ldots L_{k_m} \O_X ) \;,\]
its image $F_*\{ \P_{k_0 \ldots k_m}\}_{k_0 < \ldots < k_m}$ under the right adjoint to the equivalence (\ref{geometric-functor-1})
lands inside of the subcategory of perfect objects. In fact, it follows from Proposition 3.7.5 of 
\cite{mann2022padic6functorformalismrigidanalytic} that
the subcategory $\perf(X) \subset \qcoh(X_\bs)$ coincides with the subcategory of dualizable objetcs, and therefore, 
it suffices to prove that the object $F_*\{ \P_{k_0 \ldots k_m}\}_{k_0 < \ldots < k_m}$
is dualizable. The latter easily follows from the fact that
the equivalnce of Theorem \ref{adelic-cover} is symmetric monoidal and
that, being levelwise dualizable, the object
$\{ \P_{k_0 \ldots k_m}\}_{k_0 < \ldots < k_m}$ the object is dualizable as an object of the limit
\[\stackrel[\substack{[k_0 < k_1 < \ldots < k_m] \\ 0\le m}]{}{\lim} L_{k_0}\ldots L_{k_m} \O_X - \mod_{\qcoh(X_\bs)}\;.\]
\end{proof}

\begin{remark}
The above theorem together with Proposition \ref{factors-are-tensor-products}
and Corollary \ref{adelic-rings-are-idempotent} mean that the map \[
\spec_{X_\bs} (L_0\O_X) \sqcup \ldots \sqcup \spec_{X_\bs} (L_{\dim(X)}\O_X) \;\arr\; X_\bs
\] a closed !-covering, whose \v{C}ech nerve is given by that cubical diagram of Theorem \ref{adelic-cover}.
We will prove in section \ref{ssec-adeles-are-affinoid} that each piece of this cover is in fact affinoid.
\end{remark}

\ssec{The adelic opens are affinoid}\label{ssec-adeles-are-affinoid}

Now we would like to prove that each term $L_{k_0}\ldots L_{k_m} \O_X - \mod_{\qcoh(X_\bs)}$ is in fact a category of modules over an analytic ring.
Let $\Ga(X_\bs, -)$ denote the push-forward along $X \arr \spec(\Z)$ fro solid quasi-coherent sheaves. It induces a functor
\begin{equation}\label{geometric-local-eq3}
L_{k_0}\ldots L_{k_m} \O_X - \mod_{\qcoh(X_\bs)} \arr \Ga(X_\bs, L_{k_0}\ldots L_{k_m} \O_X) - \mod_{\Z_\bs} \;.
\end{equation}

\begin{theorem}\label{thm-adeles-are-affinoid}
\begin{enumerate}
\item[a)] The ring $\Ga(X_\bs, L_{k_0}\ldots L_{k_m} \O_X)$ is static, i.e. \[
\Ga(X_\bs, L_{k_0}\ldots L_{k_m} \O_X) \in \mod_{\Z_\bs}^{\heartsuit} \;;
\]
\item[b)] the functor (\ref{geometric-local-eq3}) is t-exact, fully faithful, and admits a symmetric monoidal left adjoint.
\end{enumerate}
\end{theorem}

We would like to split the proof into several steps. First, we prove that the ring $\Ga(X_\bs, L_{k_0}\ldots L_{k_m} \O_X)$
sits in the heart of the t-structure, see Proposition \ref{prop-global-adelic-ring-is-h-discrete}.
\begin{lemma}\label{lem-adelic-ring-is-h-discrete}
We have: \[
L_{k_0}\ldots L_{k_m} \O_X \in \qcoh(X_\bs)^\heartsuit \;.
\]
\end{lemma}
\begin{proof}
We prove the claim by the induction on $m$. The base case $m = 0$ was already treated in
Proposition \ref{prop-ultrasolid-stratum-t-exactness}, so we now prove the induction step.
By Theorem \ref{thm-global-vertices}, the ring $L_{k_0}\ldots L_{k_m} \O_X$ is
connective, and it remains to prove that $L_{k_0}\ldots L_{k_m} \O_X \in \qcoh(X_\bs)^{\ge 0}$.
To keep the notation simple, we will only treat the case $m = 1$; the proof applies to the general case without change.
We pick an affine open $\spec(A) = U \hookrightarrow X$. By Theorem \ref{thm-global-vertices} combined with the fact
that the restriction to $U_\bs$ preserves all limits and colimits, the $A_\bs$-module $L_{k_0}L_{k_1} A$ is a direct summand
in $L_{k_0} L_{k_1} \O_X$, and we can assume that $X = \spec(A)$. By Theorem \ref{thm-adeles}, \[
L_{k_0}L_{k_1} A \;\simeq\; \prod_{\substack{\q \\ \dim(A/\q) = k_1}} \leftcomp\: \bigl( \prod_{\substack{ \p \supset \q \\ \dim(A/\p) = k_0}} \leftcomp A_\p \bigr)_\q \;.
\] By the t-exactness of products, it suffices to prove that \[
\leftcomp\: \bigl( \prod_{\substack{ \p \supset \q \\ \dim(A/\p) = k_0}} \leftcomp A_\p \bigr)_\q \;\in\; \mod^{\ge 0}_{A_\bs}
\] for a fixed $\q$. Here we apply Lemma \cite[\href{https://stacks.math.columbia.edu/tag/0921}{Tag 0921}]{stacks-project}
to the completion at $\q$ and we also use Lemma \ref{lem-tensor-vs-prods} and Lemma \ref{lem-tensor-pres-comp}
to pull the tensor product with $A/\q_N$ inside: \[\begin{split}
\leftcomp\: \bigl( \prod_{\substack{ \p \supset \q \\ \dim(A/\p) = k_0}} \leftcomp A_\p \bigr)_\q &\;\simeq\; \stackrel[N \ge 0]{}{\lim} \; A/\q_N \otimes_{A_\bs} A_\q \otimes_{A_\bs} \prod_{\substack{ \p \supset \q \\ \dim(A/\p) = k_0}} \leftcomp A_\p \;\simeq\; \\ &\;\simeq\; \stackrel[N\ge 0]{}{\lim} \; A_\q \otimes_{A_\bs} \prod_{\substack{ \p \supset \q \\ \dim(A/\p) = k_0}} \leftcomp\: (A/\q_N)_\p \;.
\end{split}\]
Finally, by the induction hypothesis, $\leftcomp\: (A/\q_N)_\p \in \mod^\heartsuit_{(A/\q_N)_\bs}$, and we get the desired
coconnectivity by the left t-exactness of limits over $(\mathbb{N}, \le)^\op$.
\end{proof}

Now that we know that the ring $L_{k_0}\ldots L_{k_m} \O_X$ sits in the heart of the t-structure, the corresponding fact about the
ring of global sections will follow by flasqueness.
\begin{proposition}\label{prop-global-adelic-ring-is-h-discrete}
We have: \[
\Ga(X_\bs, L_{k_0}\ldots L_{k_m} \O_X) \in \mod_{\Z_\bs}^\heartsuit \;.
\]
\end{proposition}
\begin{proof}
For any Zariski open $j\colon U \hookrightarrow X$, the canonical map \[\xymatrix{
L_{k_0}\ldots L_{k_m} \O_X \ar[rr] && j_*j^*(L_{k_0}\ldots L_{k_m} \O_X) \;\simeq\;  j_*(L_{k_0}\ldots L_{k_m} \O_U) \\
\prod_{\substack{\p_m \in X \\ \dim \overline{\{\p_m\}} = k_m}} \leftcomp\: (\ldots \prod_{\substack{\p_0 \supset \p_1\\\dim \overline{\{\p_0\}} = k_0}} \leftcomp\: (\O_X)_{\p_0} \ldots )_{\p_m} \ar@{}[u]|\invsimeq \ar[rr] && \prod_{\substack{\p_m \in U \\ \dim \overline{\{\p_m\}} = k_m}} \leftcomp\: (\ldots \prod_{\substack{\p_0 \in U \\ \p_0 \supset \p_1\\\dim \overline{\{\p_0\}} = k_0}} \leftcomp\: (\O_X)_{\p_0} \ldots )_{\p_m} \ar@{}[u]|\invsimeq
}\] is surjective. Therefore, a standard argument shows that higher global sections $R^{>0}\Ga(X_\bs, -)$
of the object $L_{k_0}\ldots L_{k_m} \O_X \in \qcoh(X_\bs)^\heartsuit$ vanish.
\end{proof}

\sssec{} We now prove part b of Theorem \ref{thm-adeles-are-affinoid}.
The claims about the left adjoint are immediate, the proof of the t-exactness and the fully faithfulness
occupies the rest of this section. To simplify our notation, we will only consider the case $m=0$;
the argument will work in the general case without modification. To simplify out natation further, we will sometimes write $\Ga(-)$ instead
of $\Ga(X_\bs, -)$.

\begin{lemma}\label{adeles-aff-lem1}
For a Zariski open embedding $j \colon U \hookrightarrow X$, the diagram \[\xymatrix{
L_k \O_X - \mod_{\qcoh(X_\bs)} \ar[rr]^{\Ga(X_\bs, -)} \ar[d] && \Ga(L_k\O_X) - \mod_{\Z_\bs} \ar[d]^{\Ga(L_k\O_X|_{U_\bs}) \otimes_{\Ga(L_k \O_X)} -} \\
L_k \O_X|_{U_\bs} - \mod_{\qcoh(X_\bs)} \ar[rr]^{\Ga(X_\bs, -)} && \Ga(L_k\O_X|_{U_\bs}) - \mod_{\Z_\bs} \;,
}\] where \[
L_k \O_X|_{U_\bs} \;\stackrel{\define}{=}\; j_*j^*L_k\O_X \;\simeq\; \prod_{\substack{\p \in U\\ \dim(\overline{\{\p\}}) = k}} \leftcomp\:\O_{X, \p} \;,
\] commutes via the Beck-Chevalley map.
\end{lemma}
\begin{proof}
First, we note that all of the functors are continuous, so it suffices to prove that diagram commutes on modules of the form
$L_k\O_X \otimes_{X_\bs} M$.
Because the ring $L_k\O_X$ splits as a product \[
L_k\O_X \;\simeq\; \prod_{\substack{\p \in X\\ \dim(\overline{\{\p\}}) = k}} \leftcomp\:\O_{X, \p} \;\simeq\; \prod_{\substack{\p \in U\\ \dim(\overline{\{\p\}}) = k}} \leftcomp\:\O_{X, \p} \times \prod_{\substack{\p \in X\setminus U\\ \dim(\overline{\{\p\}}) = k}} \leftcomp\:\O_{X, \p} \;=:\; L_k\O_X|_{U_\bs} \times R \;,
\] we have \[\xymatrix{
L_k\O_X \otimes_{X_\bs} M \simeq (L_k\O_X|_{U_\bs} \otimes M) \times (R \otimes M) \ar@{|->}[rr] \ar@{|->}[d] && \Ga(L_k\O_X|_{U_\bs} \otimes M) \times \Ga(R \otimes M) \ar@{|->}[d]^{\Ga(L_k\O_X|_{U_\bs}) \otimes_{\left(\Ga(L_k\O_X|_{U_\bs}) \times \Ga(R)\right)} -} \\
L_k\O_X|_{U_\bs} \otimes M \ar@{|->}[rr] && \Ga(L_k\O_X|_{U_\bs} \otimes M) \;.
}\]
\end{proof}

We now present $X$ as a finite colimit (in classical schemes) \[
\stackrel[i]{}{\colim} U_i \;\stackrel{\simeq}{\arr}\; X \;,
\] where each $U_i \subseteq X$ is an affine open. According to \cite{condensed}, Theorem 9.8, the cover
$\{U_i{}_\bs \arr X_\bs \}_i$ satisfies the upper-$*$ descent, which means that the functor \[
\qcoh(X_\bs) \;\arr\; {\lim_i}^* \qcoh({U_i}_\bs) \;,
\] where the limit is computed via the upper-$*$ pull-backs and every map $\qcoh(X_\bs) \arr \qcoh({U_i}_\bs)$ is the upper-$*$ pull-back,
is an equivalence.

\begin{lemma}\label{adeles-aff-lem2}
The functor \[
L_k\O_X - \mod_{\qcoh(X_\bs)} \stackrel{\text{upper-}*}{\:\;\arr\;\:} {\lim_i}^* \;\: L_k\O_{X}|_{{U_i}_\bs} - \mod_{\qcoh(X_\bs)} \;,
\] where the limit is computed via the upper-$*$ functors, is fully faithful.
\end{lemma}
\begin{proof}
By the projection formula, it suffices to prove that the map \begin{equation}\label{geometric-adeles-as-limit-on-a-cover}
L_k \O_X \;\arr\; \lim_i L_k\O_{X}|_{{U_i}_\bs}
\end{equation}
is an isomorphism, which is true because $\{U_i{}_\bs \arr X_\bs \}_i$ satisfies the upper-$*$ descent.
\end{proof}
Applying $\Ga(X_\bs, -)$ to the isomorphism (\ref{geometric-adeles-as-limit-on-a-cover})
and using the same argument, we get the following lemma.
\begin{lemma}\label{adeles-aff-lem3}
The functor \[
\Ga(L_k\O_X) - \mod_{\Z_\bs} \stackrel{\text{upper-}*}{\:\;\arr\;\:} {\lim_i}^* \;\: \Ga(L_k\O_{X}|_{{U_i}_\bs}) - \mod_{\Z_\bs}
\] is fully faithful. Again, all the functors are the upper-$*$ functors.
\end{lemma}

\sssec{Proof of the fully faithfulness part of Theorem \ref{thm-adeles-are-affinoid}{\color{red}b}}
Let us take a look at the right adjoint to the upper-$*$ functor \[
{\lim_i}^* \;\: \Ga(L_k\O_{X}|_{{U_i}_\bs}) - \mod_{\Z_\bs} \;\:\arr\;\: {\lim_i}^* \;\: L_k\O_{X}|_{{U_i}_\bs} - \mod_{\qcoh(X_\bs)} \;.
\] It follows from Lemma \ref{adeles-aff-lem1} that
it is given by the limit of the right adjoints and fits into a commutative square \[\xymatrix{
L_k\O_X - \mod_{\qcoh(X_\bs)} \ar[rr]^{\Ga(X_\bs, -)} \ar[d] && \Ga(L_k\O_X) - \mod_{\Z_\bs} \ar[d] \\
{\lim_i}^* \;\: L_k\O_{X}|_{{U_i}_\bs} - \mod_{\qcoh(X_\bs)} \ar[rr] && {\lim_i}^* \;\: \Ga(L_k\O_{X}|_{{U_i}_\bs}) - \mod_{\Z_\bs} \;.
}\]
On the other hand, by affinness of $U_i$, each right adjoint
\[\Ga(X_\bs, -) \colon L_k\O_{X}|_{{U_i}_\bs} - \mod_{\qcoh(X_\bs)} \arr \Ga(L_k\O_{X}|_{{U_i}_\bs}) - \mod_{\Z_\bs}\]
is fully faithful, and therefore,
by Lemma \ref{adeles-aff-lem2} and Lemma \ref{adeles-aff-lem3}, the top horizontal arrow of the square must as well be fully faithful.

\sssec{Proof of the t-exactness part of Theorem \ref{thm-adeles-are-affinoid}{\color{red}b}}
We keep using notation already introduced in the course of proving the theorem. First, we note that our functor is
automatically left t-exact, so it suffices to prove that it is right t-exact.
Using Lemma \ref{adeles-aff-lem1}, we see that,
for each $i$, the composite \[
L_k\O_X - \mod_{\qcoh(X_\bs)} \stackrel{\Ga(X_\bs, -)}{\;\arr\;} \Ga(L_k\O_X) - \mod_{\Z_\bs} \stackrel{\Ga(L_k\O_X |_{{U_i}_\bs}) \otimes_{\Ga(L_k\O_X)} - }{\;\arr\;} \Ga(L_k\O_X |_{{U_i}_\bs}) - \mod_{\Z_\bs}
\] is t-exact, and therefore, in view of Lemma \ref{adeles-aff-lem3}, the claim would immediately follow from the following lemma.
\begin{lemma}
Let $R \in \calg^\heartsuit(\mod_{\Z_\bs})$ be a solid non-derived commutative ring. Let also
$\{R_i\}_{i \in I} \subset \calg^\heartsuit(R - \mod_{\Z_\bs})$ be a finite collection of solid non-derived commutative $R$-algebras
satisfying these conditions:
\begin{enumerate}
\item[a)] each $R_i$ is flat and idempotent over $R$;
\item[b)] the functors $\{R_i \otimes_R - \colon R-\mod_{\Z_\bs} \arr R_i - \mod_{\Z_\bs}\}$ are jointly conservative;
\item[c)] each of the maps $R \arr R_i$ is surjective.
\end{enumerate}
For any subset $I^\prime \subseteq I$, we define a commutative algebra $R_{I^\prime}$ by the formula \[
R_{I^\prime} \;\stackrel{\define}{=}\; \lim_{\{i_1, \ldots, i_m\} \subseteq I^\prime} R_{i_1} \otimes_R \ldots \otimes_R R_{i_m} \;,
\] where the limit is taken over all non-empty subsets of $I^\prime$.
We claim: \begin{enumerate}
\item[1)] the map $R \arr R_I$ is an isomorphism;
\item[2)] for any $I^\prime = I_1^\prime \sqcup I_2^\prime$, the square \[\xymatrix{
R_{I^\prime} \ar[rr]\ar[d] && R_{I_2^\prime} \ar[d] \\
R_{I_1^\prime} \ar[rr] && R_{I_1^\prime} \otimes_R R_{I_2^\prime}
}\] is a pull-back square;
\item[3)] for any $I^\prime \subseteq I$, the algebra $R_{I^\prime}$ is non-derived
(i.e. $R_{I^\prime} \in \calg^\heartsuit(\mod_{\Z_\bs})$), 
idempotent and flat over $R$, and the map $R \arr R_{I^\prime}$ is surjective;
\item[4)] any module $M \in R-\mod_{Z_\bs}$ satisfying $R_i \otimes_R M \in R_i - \mod^{\le 0}_{\Z_\bs}$ is in fact connective---
\[
M \in R-\mod^{\le 0}_{\Z_\bs} \;.
\]
\end{enumerate}
\end{lemma}
\begin{proof}
It follows from the idempotency that, for each $i\in I$, the map $R_i \arr R_i \otimes_R R_I$ is an isomorphism. Then the first
claim follows by the joint conservativity.

In the second claim, all of the terms of the square lie inside the full subcategory of $R-\mod_{\Z_\bs}$ generated by
$\{ R_i-\mod_{\Z_\bs} \subseteq R-\mod_{\Z_\bs}\}_{i \in I^\prime}$. Therefore, it suffices to prove that the square becomes a
pull-back square upon application of $R_i \otimes_R -$ for each $i \in I^\prime$. The latter is true because, for any $I^{\prime\prime}$
and any $i \in I^{\prime\prime}$, the map $R_i \otimes_R R_{I^{\prime\prime}} \arr R_{I^{\prime\prime}}$ is an isomorphism.

In the third claim, the idempotency is clear from the definition, and the rest is proved by induction on the size of the set $I^\prime$
using the second claim to prove the induction step.

In the forth claim, using the previous claims,
one proves by induction that $R_{I^\prime} \otimes_R M \in R_{I^\prime}-\mod^{\le 0}_{\Z_\bs}$.
Namely, pick a decomposition $I^\prime = I^\prime_1 \sqcup I_2^\prime$. By the second claim, we get a pull-back square
\[\xymatrix{
M \otimes_R R_{I^\prime} \ar[rr]\ar[d] && M \otimes_R R_{I_2^\prime} \ar[d] \\
M \otimes_R R_{I_1^\prime} \ar[rr] && M \otimes_R R_{I_1^\prime} \otimes_R R_{I_2^\prime} &.
}\] In particular, we get a long exact sequence \[\begin{split}
\ldots \;\to\; H^0(M \otimes_R (R_{I_1^\prime} \oplus R_{I_2^\prime})) \;\to\; H^0(M \otimes_R R_{I_1^\prime} \otimes_R R_{I_2^\prime}) &\;\to\; H^1(M \otimes_R R_{I^\prime}) \;\to\; \\
\;\to\; 0 \;\to\; 0 &\;\to\; H^2(M \otimes_R R_{I^\prime}) \;\to\; \ldots \;.
\end{split}\] So it suffices to prove that the map \[
H^0(M \otimes_R (R_{I_1^\prime} \oplus R_{I_2^\prime})) \;\to\; H^0(M \otimes_R R_{I_1^\prime} \otimes_R R_{I_2^\prime})
\] is surjective, which immediately follows from the second claim.
\end{proof}

%% file: appendixB.tex
\appendix
\section{Certain properties of the $k$-completion and $k$-complete solid modules}\label{app-k-comp}

In this appendix, we fill the gap that we left in section \ref{ssec-some-prop-of-skeleton-comp} and 
establish certain properties of the completion at the $k$-skeleton \[
S_k \;=\; \cup_{\substack{Z \subseteq \spec(A) \\ Z \text{ is closed}}} Z \;\in\; \spec(A)^\spc
\]
and $S_k$-complete modules.
As usual, in this section, $A$ denotes a commutative ring of finite type over $\Z$.

\ssec{} Here we collect a couple of basic facts
before moving on to the main result of this appendix---Proposition \ref{k-comp-is-monoidal}.

\begin{proposition}\label{appB-conservativity-of-primes-on-k-comp}
The set of functors \[
\{ (-)^\comp_\p \colon \mod_{A_\bs} \arr \mod_{A_\bs}\}_{\substack{\p \in \spec(A) \\ \dim(\p) \le k}}
\] is jointly conservative on the full subcategory of $S_k$-complete modules. Moreover, the set of functors \[
\{ - \otimes_{A_\bs} (A/\p)_\bs \colon \mod_{A_\bs} \arr \mod_{(A/\p)_\bs}\}_{\substack{\p \in \spec(A) \\ \dim(\p) \le k}}
\] is jointly conservative on the full subcategory of $S_k$-complete modules.
\end{proposition}
\begin{proof}
The first assertion follows immediately from
Proposition \ref{prop-lim-over-primes}, the second one follows from the first one---see
\cite[\href{https://stacks.math.columbia.edu/tag/0G1U}{Lemma 0G1U}]{stacks-project}.
\end{proof}

\sssec{}
We now prove---see Propositon \ref{lem-k-comp-restricts-to-closed}---that the completion at a specialization closed subset, such as $S_k$,
is compatible with box restriction to a closed subscheme. This result will help us run the induction
during the proof of Proposition \ref{k-comp-is-monoidal}.
\begin{lemma}\label{lem-local-cohomology-base-change}
Given a map $\phi\colon \spec(B) \arr \spec(A)$ between finite type affine schemes and a specialization closed subset
$T \in \spec(A)^\spc$, the subset $\phi^{-1}(T) \subseteq \spec(B)$ is specialization closed, and the canonical map
\begin{equation}\label{eq-local-cohomology-base-change-map}
B\otimes_A \Ga_T \;\arr\; \Ga_{\phi^{-1}(T)}
\end{equation}
is an isomorphism.
\end{lemma}
\begin{proof}
The first claim is immediate. Therefore, by Proposition \ref{prop-hopkins-neeman}, in order to prove that
the map (\ref{eq-local-cohomology-base-change-map}) is an isomorphism it suffices to prove that the coidempotent objects
have the same support. On other words, it suffices to prove that the map (\ref{eq-local-cohomology-base-change-map})
becomes an isomorphism upon tensoring with $B_\q/\q$ for each $\q \in \spec(B)$.
On the one hand, it follows from the first claim of Proposition \ref{prop-hopkins-neeman} that \[
B_\q/\q \otimes_B \Ga_{\phi^{-1}(T)} \;\simeq\; \left\{ \begin{array}{cc}
B_\q/\q , & \text{ if } \phi(\q) \in T; \\
0, & \text{ otherwise.}
\end{array}\right.
\] On the other hand, \[
B_\q/\q \otimes_A \Ga_{T} \;\simeq\; B_\q/\q \otimes_{A_{\phi(\q)}/\phi(\q)} A_{\phi(\q)}/\phi(\q) \otimes_A \Ga_{T} \;\simeq\; \left\{ \begin{array}{cc}
B_\q/\q, & \text{ if } \phi(\q) \in T; \\
0, & \text{ otherwise.}
\end{array}\right.
\] 
\end{proof}
\begin{proposition}\label{lem-k-comp-restricts-to-closed}
For an eventually connective $A_\bs$-module $M \in \mod^{< +\infty}_{A_\bs}$,
an ideal $I \subset A$, and a specialization closed subset $T \subseteq \spec(A)$,
the canonical map \[
M^\thicksim_{T} \otimes_{A_\bs} (A/I)_\bs \;\arr\; ( M\otimes_{A_\bs} (A/I)_\bs )^\thicksim_{T \cap \spec(A/I)}
\] is an isomorphism. In particular, the map \[
M^\thicksim_{S_k} \otimes_{A_\bs} (A/I)_\bs \;\arr\; ( M\otimes_{A_\bs} (A/I)_\bs )^\thicksim_{S_k} \;,
\] where we, abusing the notation, denote both the $k$-skeleton of $\spec(A)$ and of $\spec(A/I)$ by the same letter $S_k$,
is an isomorphism.
\end{proposition}
\begin{proof}
By Lemma \ref{lem-local-cohomology-base-change}, it suffices to prove that the map \[
\uhom_{A_\bs} (\Ga_T, M) \otimes_{A_\bs} (A/I)_\bs \;\arr\; \uhom_{(A/I)_\bs} (\Ga_T \otimes_A A/I, M\otimes_{A_\bs} (A/I)_\bs)  
\] is an isomorphism. Because the object $\Ga_T \in \mod_A$ has finite projective dimension by
Proposition \ref{prop-comp-alg-is-idmpt-and-cmp},
it suffices to prove that the functor $(A/I)_\bs \otimes_{A_\bs} -$ preserves products t-bounded on the right, which immediately follows
from Lemma \ref{lem-tensor-vs-prods} combined with the fact that $A/I$ is pseudo compact as a $A_\bs$-module.
\end{proof}

\ssec{} The main result of this appendix is the following statement, which says that the functor of completion at
the skeleton is symmetric monoidal on the full subcategory consisting of eventually connective objects.
In particular,
it follows from this statement, that the solid tensor product of a pair of $S_k$-complete eventually connective modules is $S_k$-complete.
For the sake of completeness, we also deduce in Corollary \ref{appB-sing-ideal-comp-is-monoidal}
the corresponding result for the functor of completion at a single ideal, which was previously recordered in Proposition A.17, \cite{Artem}.
\begin{proposition}\label{k-comp-is-monoidal}
The functor of $S_k$-completion (which is a priori lax symmetric monoidal) is symmetric monoidal when restricted to the
subcategory of eventually connective modules, i.e. for any pair $M,N \in \mod^{< +\infty}_{A_\bs}$, the map
\begin{equation}\label{k-comp-is-monoidal-eq}
M^\thicksim_{S_k} \otimes_{A_\bs} N^\thicksim_{S_k} \;\arr\; \bigl( M\otimes_{A_\bs} N \bigr)^\thicksim_{S_k}
\end{equation}
is an isomorphism. In particular, the tensor product of a pair of $S_k$-complete eventually connective $A_\bs$-modules is $S_k$-complete.
\end{proposition}

\sssec{}
The proof of Proposition \ref{k-comp-is-monoidal} occupies the rest of this appendix.
First, we argue that it suffices to prove the proposition in the case
where the ring $A$ is integral. Indeed, if $d$ denotes the dimension of $A$, then the functor of completion at $S_d$ is the identity functor.
In particular, both the source and the target of \ref{k-comp-is-monoidal-eq} are $S_d$-complete. We can now reduce to the integral
case by Corollary \ref{appB-conservativity-of-primes-on-k-comp} and Proposition \ref{lem-k-comp-restricts-to-closed}.
From now on we assume that the ring $A$ is integral.

We will first prove the proposition in the case where $k = d-1$.
Let $S = A\setminus0 \subset A$
denote the multiplicative subset
of non-zerodivisors. It follows from Proposition \ref{prop-sm-sp-map-arbitrary-meets} that the map \[
K^\thicksim_{S_{d-1}} \;\arr\; \lim_{f \in S} K^\comp_f
\] is an isomoprhism for any module $K \in \mod_{A_\bs}$.
Since we are working under the finite type hypothesis, the set $S$ is countable. Let us fix an enumeration \[\xymatrix{
\N \ar[rr]^{\simeq} && S \\
n \ar@{|->}[rr] \ar@{}[u]|\invertin && f_n \ar@{}[u]|\invertin \;,
}\] and let us denote by $\widetilde F_n$ the product $f_1 f_2 \ldots f_n \in S$. We get a string of isomorphisms \[\begin{split}
K^\thicksim_{S_{d-1}} &\;\simeq\; \lim_{f \in S} K^\comp_f \;\simeq\\
&\;\simeq\; \lim_{n} K^\comp_{\widetilde F_n} \;\simeq\\
&\;\simeq\; \lim_{n} \lim_m K / \widetilde F_n^m \;\simeq\\
&\;\simeq\; \lim_n K / \widetilde F_n^n \;\simeq\\
&\;\simeq\; \lim_n K / F_n \;,
\end{split}\]
where $F_n \in A$ denotes the element $\widetilde F_n^n \in A$.
Now let us assume that the module $K$ is eventually connective---$K \in \mod^{< +\infty}_{A_\bs}$. It can be represented as a bounded
above complex of the form \[
\bigl( \ldots \to \bigoplus\prod A \to \bigoplus \prod A \to \bigoplus \prod A\bigr) \;,
\] which we will from now on also denote by $K$. We write the complex $K$ as a direct colimit of its naive truncations:
$K \simeq \colim_{l\ge 0} \sigma_{\ge -l} K$. Above, we expressed the completion $K^\thicksim_{S_{d-1}}$ as a limit over the natural
numbers---\[
K^\thicksim_{S_{d-1}} \;\simeq\; \lim_{n} K/F_n \;\simeq\; \fib \bigl( \prod_{n}K/F_n \to \prod_{n}K/F_n \bigr) \;.
\] It follows from this expression that the completion
distributes with the colimit of naive truncations: \[
K^\thicksim_{S_{d-1}} \;\simeq\; (\stackrel[l\ge 0]{}{\colim} \sigma_{\ge -l} K)^\thicksim_{S_{d-1}} \;\simeq\;\; \stackrel[l\ge 0]{}{\colim} (\sigma_{\ge -l} K)^\thicksim_{S_{d-1}} \;.
\] So we have reduced the proof of the case $k = d-1$ of Proposition \ref{k-comp-is-monoidal} to the case where both $M$ and $N$ have
form $\bigoplus_Q\prod_T A$. Moreover, because the functor of completion commutes with tensoring by $\prod_TA$ by
Lemma \ref{lem-tensor-pres-comp},
we can assume that both $M$ and $N$ have form $\bigoplus A$. Finally, because the functor of completion commutes with
$\omega_1$-filtered colimits (the functor is given by homing out of
$\colim_{f \in A\setminus 0} A[\frac{1}{f}]/A[-1]$,
which can be represented as a countable filtered colimit of perfect $A$-modules and hence $\omega_1$-compact),
we can further reduce to the case where
both $M$ and $N$ have form $\bigoplus_{\N} A$.
Summing up, the proof of Proposition \ref{k-comp-is-monoidal} in the case $k = d-1$ reduces to the following statement.
\begin{lemma}\label{reduced-tensor-product-of-k-complete-is-k-complete}
The canonical map \[
\bigl( \bigoplus_\N A \bigr)^\thicksim_{S_{d-1}} \otimes_{A_\bs} \bigl( \bigoplus_\N A \bigr)^\thicksim_{S_{d-1}} \;\arr\; \bigl( \bigoplus_{\N\times\N} A \bigr)^\thicksim_{S_{d-1}}
\] is an isomorphism.
\end{lemma}
We split the proof of Lemma \ref{reduced-tensor-product-of-k-complete-is-k-complete} into several steps.
\begin{lemma}
The natural map \[
\stackrel[\substack{\phi\colon\N\arr\N\cup\{0\}\\\phi(n)\arr\infty}]{}{\colim} \fib\bigl( \prod_n A^\thicksim_{S_{d-1}} \;\arr\; \prod_n A/F_{\phi(n)} \bigr) \;\arr\; \bigl(\bigoplus_\N A\bigr)^\thicksim_{S_{d-1}} \;,
\] where, for $m \ge 1$, $F_m$ denotes $\widetilde F_m$ raised to the power of $m$, and $F_0$ denotes $1 \in A$, is an isomorphism.
\end{lemma}
\begin{proof}
We present the completed direct sum as a fiber \[
\fib \bigl( \prod_{m\in\N} \bigoplus_\N A/F_m \arr \prod_{m\in\N} \bigoplus_\N A/F_m \bigr) \;.
\] Then, because there is no $R^1\lim$ by the Mittag-Leffler condition, this fiber is actually a kernel 
\begin{equation}\label{geom-lem-single-comp-colim-form-proof-eq1}
\ker \bigl( \prod_{m\in\N} \bigoplus_\N A/F_m \arr \prod_{m\in\N} \bigoplus_\N A/F_m \bigr) \;.
\end{equation}
We now rewrite the source of (\ref{geom-lem-single-comp-colim-form-proof-eq1}) as \[\begin{split}
\prod_{m\in\N} \bigoplus_\N A/F_m \;\simeq\; \prod_{m\in\N} \stackrel[\text{finite } S \subset \N]{}{\colim}\bigoplus_S A/F_m &\;\simeq\; \stackrel[\{\text{finite } S_m \subset \N\}_m]{}{\colim} \prod_{m\in\N} \bigoplus_{S_m} A/F_m \;\simeq \\
&\;\simeq\; \stackrel[\substack{\{\text{finite } S_m \subset \N\}_m \\ S_m \subseteq S_{m+1}, |S_m| \arr \infty }]{}{\colim} \prod_{m\in\N} \bigoplus_{S_m} A/F_m \;,
\end{split}\]
where the second isomorphism uses AB6 and the last isomorphism results from a cofinality argument. For a fixed collection $ S_1 \subseteq S_2 \subseteq \ldots \subseteq \N$ of nested finite subsets, we 
construct a function $\phi \colon \N \arr \N$ having value $m-1$ on $S_{m} \setminus S_{m-1}$.
We observe that this gives us an isomorphism of diagrams \[
\left(\substack{\{\text{finite } S_m \subset \N\}_m \\ S_m \subseteq S_{m+1}, |S_m| \arr \infty }\right) \;\simeq\; \left(\substack{\phi\colon \N \arr \N\cup\{0\} \\ \phi(n) \arr \infty}\right) \;,
\] and it remains to prove that, for a fixed such collection of subsets $\{S_m\}_{m \in \N}$ and the corresponding function
$\phi\colon\N\arr\N\cup\{0\}$, we get an identification \begin{equation}\label{geom-lem-single-comp-colim-form-proof-eq2}
\ker\bigl( \prod_{n\in\N} A^\thicksim_{S_{d-1}} \arr \prod_{n\in\N} A / F_{\phi(n)} \bigr) \;\stackrel{\simeq}{\arr}\; \ker \bigl( \prod_{m\in\N} \bigoplus_{S_m} A/F_m \arr \prod_{m\in\N} \bigoplus_\N A/F_m \bigr) \;.
\end{equation}
We prove that both sides of (\ref{geom-lem-single-comp-colim-form-proof-eq2}) identify as submodules of $\prod_\N A^\comp_{d-1}$. The embedding of the right-hand side
into the product composed with the projection onto the $n$-th factor \[
\prod_\N A^\thicksim_{S_{d-1}} \;\arr\; A^\thicksim_{S_{d-1}} \;\simeq\; \ker \bigl( \prod_{m\in\N} A/F_m \arr \prod_{m\in\N} A/F_m \bigr)
\] is given by projections onto the $n$-th factor \[
\bigoplus_{S_m} A/F_m \;\hookrightarrow\; \bigoplus_{\N} A/F_m \;\arr\; A/F_m \quad\text{and}\quad \bigoplus_{\N} A/F_m \;\arr\; A/F_m \;;
\] in particular, we see that a tuple \[
(a_1, a_2, a_3, \ldots ) \in \prod_\N A^\thicksim_{S_{d-1}}
\] is in the image if and only if, for each $m \in \N$ and $n \in S_m \setminus S_{m-1}$, the element $a_n \in A^\thicksim_{S_{d-1}}$ vanishes modulo $F_{m-1}$. 
\end{proof}

By a similar argument, we get the following variant of the same statement.
\begin{lemma}
The natural map \[
\stackrel[\substack{\phi\colon\N\arr\N\cup\{0\}\\\f(n)\arr\infty}]{}{\colim}\;\stackrel[\substack{\g\colon\N\arr\N\cup\{0\}\\\g(m)\arr\infty}]{}{\colim} \fib\bigl( \prod_{n,m} A^\thicksim_{S_{d-1}} \;\arr\; \prod_{n,m} A/(F_{\f(n)} F_{\g(m)}) \bigr) \;\arr\; \bigl(\bigoplus_{\N\times\N} A\bigr)^\thicksim_{S_{d-1}}
\] is an isomorphism.
\end{lemma}
\begin{proof}
For the completed direct sum, we write \[\begin{split}
\bigl(\bigoplus_{\N\times\N} A\bigr)^\thicksim_{S_{d-1}} &\;\simeq\; \lim_{n} \; \bigl(\bigoplus_{\N\times\N} A/F^2_n \bigr)\;\simeq\\
&\;\simeq\; \lim_{n,m} \; \bigl(\bigoplus_{\N\times\N} A/(F_n F_m) \bigr)\;\simeq\\
&\;\simeq\; \ker \bigl( \prod_{n,m} \bigoplus_{\N\times\N} A/(F_nF_m) \arr \prod_{n,m} \bigoplus_{\N\times\N} A/(F_nF_m) \bigr) \;,
\end{split}\]
and now, by an identical argument, we get \[\begin{split}
\bigl(\bigoplus_{\N\times\N} A\bigr)^\thicksim_{S_{d-1}} \;\simeq\; \stackrel[\substack{\t\colon\N\times\N\arr(\N\cup\{0\})\times(\N\cup\{0\})\\\t(n,m)\arr\infty}]{}{\colim} \fib\bigl( \prod_{n,m} A^\thicksim_{S_{d-1}} \;\arr\; \prod_{n,m} A/(F_{\t_1(n)} F_{\t_2(m)}) \bigr) \;.
\end{split}\]
It remains to observe that the subdiagram \[
\left(\substack{(\f, \g)\colon \N\times\N \arr (\N\cup\{0\})\times(\N\cup\{0\}) \\ \f(n), \g(m) \arr \infty}\right) \subseteq \left( \substack{\t\colon \N\times\N \arr (\N\cup\{0\})\times(\N\cup\{0\}) \\ \t(n,m) \arr \infty} \right)
\] is cofinal. Indeed, since both digrams are filtered, it suffices, for each $\t$,
to construct a function of the form $(\f, \g)$ satisfying $(\f(n), \g(m)) \le \t(n,m)$ for all $n,m$. One easily observes that
setting \[
(\f(l), \g(l)) \;=\; ( \stackrel[n+m \ge l]{}{\mathrm{min}} \t_1(n,m), \stackrel[n+m \ge l]{}{\mathrm{min}} \t_2(n,m) )
\] does the job.
\end{proof}

Putting the last two statements together, we get Lemma \ref{reduced-tensor-product-of-k-complete-is-k-complete}.
\begin{proof}[Proof of Lemma \ref{reduced-tensor-product-of-k-complete-is-k-complete}]
\[\begin{split}
&\bigl( \bigoplus_\N A \bigr)^\thicksim_{S_{d-1}} \otimes_{A_\bs} \bigl( \bigoplus_\N A \bigr)^\thicksim_{S_{d-1}} \;\simeq\;\\
&\stackrel[\substack{\f\colon \N \arr \N\cup\{0\} \\ \f(n) \arr \infty}]{}{\colim}\;\stackrel[\substack{\g\colon \N \arr \N\cup\{0\} \\ \g(m) \arr \infty}]{}{\colim} \Bigl( \prod_{n\in\N}\prod_{m\in\N} A^\thicksim_{S_{d-1}} \arr \prod_{n\in\N}\prod_{m\in\N} A/F_{\f(n)} \times A/F_{\g(m)} \arr \prod_{n\in\N}\prod_{m\in\N} A/F_{\f(n)} \otimes_A A/F_{\g(m)} \Bigr) \;\simeq\\
&\stackrel[\substack{\f\colon \N \arr \N\cup\{0\} \\ \f(n) \arr \infty}]{}{\colim}\;\stackrel[\substack{\g\colon \N \arr \N\cup\{0\} \\ \g(m) \arr \infty}]{}{\colim} \Bigl( \prod_{n\in\N}\prod_{m\in\N} A^\thicksim_{S_{d-1}} \arr \prod_{n\in\N}\prod_{m\in\N} A/(F_{\f(n)}F_{\g(m)}) \Bigr) \;\simeq\\
&\bigl( \bigoplus_{\N\times\N} A \bigr)^\thicksim_{S_{d-1}} \;,
\end{split}\]
where to get from the 2nd row to the 3rd we use the pull-back square \[\xymatrix{
A/(F_{\f(n)}F_{\g(m)}) \ar[d] \ar[rr] && A/F_{\f(n)} \ar[d] \\
A/F_{\g(m)} \ar[rr] && A/F_{\f(n)} \otimes_A A/F_{\g(m)} \;.
}\]
\end{proof}

\sssec{}
We will now finish the proof of Proposition \ref{k-comp-is-monoidal} by the induction on the dimension of $A$.
First, we note that the base cases of dimension $0$ and $1$ are, respectively, vacuous and already dealt with by
Lemma \ref{reduced-tensor-product-of-k-complete-is-k-complete}.
Let us now assume that the statement holds for any integral $A$ of dimension up to $d > 1$
and all $k \in \{0,\dots,d\}$. We will now prove that the map
\begin{equation}\label{eq-the-map-proof-tensor-product-of-k-complete-is-k-complete}
M^\thicksim_{S_k} \otimes_{A_\bs} N^\thicksim_{S_k} \;\arr\; \bigl(M\otimes_{A_\bs}N\bigr)^\thicksim_{S_k}
\end{equation} is an isomorphism for some fixed integral ring $A$ of dimension $d+1$, some $k \in \{0,\dots,d\}$,
and a pair $M, N \in \mod^{< +\infty}_{A_\bs}$ of eventually connective modules.
By Lemma \ref{lem-some-exms-of-k-complete} and Lemma \ref{reduced-tensor-product-of-k-complete-is-k-complete}, both sides of
(\ref{eq-the-map-proof-tensor-product-of-k-complete-is-k-complete}) are $S_d$-complete.
By Proposition \ref{appB-conservativity-of-primes-on-k-comp},
it suffices to prove that the map $(\ref{eq-the-map-proof-tensor-product-of-k-complete-is-k-complete}) \otimes_{A_\bs} (A/\p)_\bs$ is an isomorphism
for each prime $\p \subset A$ of dimension up to $d$, but this follows
from Lemma \ref{lem-k-comp-restricts-to-closed} and the induction hypothesis. It remains to note that the case $k = d+1$ is vacuous.
Proposition \ref{k-comp-is-monoidal} is now proved.

\begin{corollary}[cf. Proposition A.17, \cite{Artem}]\label{appB-sing-ideal-comp-is-monoidal}
The functor \[
(-)^\comp_I \colon \mod^{< +\infty}_{A_\bs} \;\arr\; \mod^{< +\infty}_{A_\bs}
\] of completion at an ideal $I\subset A$ is symmetric monoidal. In particular,
the tensor product of a pair of $I$-complete eventually connective $A_\bs$-modules is $I$-complete.

More generally, given a specialization closed subset $T \subseteq \spec(A)$, the functor \[
(-)^\thicksim_T \colon \mod^{< +\infty}_{A_\bs} \;\arr\; \mod^{< +\infty}_{A_\bs}
\] of completion at $T$ is symmetric monoidal. In particular,
the tensor product of a pair of $T$-complete eventually connective $A_\bs$-modules is $T$-complete.
\end{corollary}
\begin{proof}
Pick a pair $M,N \in \mod_{A_\bs}^{< +\infty}$ of eventually connective $A_\bs$-modules. We prove that the map
\begin{equation}\label{cor-comp-is-monoidal-map}
M^\thicksim_T \otimes_{A_\bs} N^\thicksim_T \;\arr\; (M\otimes N)^\thicksim_T
\end{equation}
is an isomorphism. First, we note that both the source and the target are $S_d$-complete, where $d$ denotes the Krull dimension of $A$.
Therefore, by Proposition \ref{appB-conservativity-of-primes-on-k-comp}, it suffices to prove that the map
$A/\p_\bs \otimes_{A_\bs} (\ref{cor-comp-is-monoidal-map})$ is an isomorphism for each prime ideal $\p \subset A$.
In other words, in the view of Lemma \ref{lem-k-comp-restricts-to-closed}, we can assume that the ring $A$ is integral.

In the case of an integral ring, we either have $T = \spec(A)$ or $T \subseteq S_{d-1}$. The former case is trivial;
in the latter case, both the source and the target of the map (\ref{cor-comp-is-monoidal-map}) are $S_{d-1}$-complete by
Proposition \ref{k-comp-is-monoidal}.
Therefore, by Proposition \ref{appB-conservativity-of-primes-on-k-comp}, it suffices to prove that the map
$A/\p_\bs \otimes_{A_\bs} (\ref{cor-comp-is-monoidal-map})$ is an isomorphism for each prime ideal $\p \subset A$ of codimension
at least one. Now we use Lemma \ref{lem-k-comp-restricts-to-closed} and reduce the claim to smaller dimension, where it holds
by the induction hypothesis.
\end{proof}